# Relational Composition of Physical Systems
## A Categorical Approach

Owen Lynch

Submitted in fulfillment of the requirements for a MSc. in Mathematics from Universiteit Utrecht

*July, 2022*

# PREFACE

*Every mathematician knows that it is impossible to understand any elementary course in thermodynamics.*
Vladimir Arnol'd

This thesis began when I encountered the phenomenon mentioned by Arnol'd above. I was taking a course in statistical mechanics from my adviser Wioletta Ruszel, when I realized that I didn't understand the elementary thermodynamics that the authors of the statistical mechanics textbook assumed as background intuition. This lead me on a quest to formulate some of the most basic assumptions of elementary thermodynamics in a formal way, so that even a mathematician could understand. Moreover, I wished to put thermodynamics on a categorical footing, so as to talk about the composition of thermodynamical systems. I had been thinking about this casually for some time when I happened to spend several weeks in the summer of 2021 at the Topos Institute at the same time as John Baez was visiting. It turned out that he had similar intentions and interests with regards to thermodynamics, as did several other people at Topos, and we all started talking about the subject.

Eventually, this turned into a collaboration between John Baez, his student Joe Moeller, and me, and we wrote a paper on what we called "thermostatics". In this paper, we derived many features of thermodynamical equilibrium from some basic assumptions about entropy functions [1]. However, this was just a beginning, and to keep investigating thermodynamics I asked John Baez if he would like to advise my thesis, which would be a continuation of the work.

For my thesis, I wanted to understand *open* thermo*dynamical* systems, that is thermodynamical systems that could take flows of energy, heat, volume, and particles from the outside world. I was (and still am) very interested in the idea that a open system can seemingly violate the second law of thermodynamics, spontaneously becoming more ordered and complex by continually expelling waste heat. Two works along these lines that particularly spoke to me were Haddad's *A Dynamical Systems Theory of Thermodynamics* [2], which forefronted flows of heat as the basic unit of study, and Schnakenberg's *Thermodynamical Network Analysis of Biological Systems* [3], which talked about some of the phenomena that could happen in open systems, like active transport. For my final paper in a mathematical biology class, I wrote about Schnakenberg, and learned about *bond graphs*, which were a notation that Schnakenberg used to analyze his thermodynamic systems. I deliberately left out all mention of category theory from this paper, as it was a mathematical biology class, but studying Schnakenberg made me convinced that taking seriously the view of flows of energy that I had found within bond graphs could lead to a good compositional theory of thermodynamic systems for my thesis.

It was this line of thinking that lead me to stumble upon the rich theory of port-Hamiltonian systems, and I started reading many papers written by van der Schaft and collaborators. The theory of port-Hamiltonian systems was already quite well-developed, including the theory of how to compose port-Hamiltonian systems. However, there was not yet a comprehensive account of how category theory could be used to structure this composition, and so I sensed an area to which I could make a contribution. This thesis is my attempt to put the composition of port-Hamiltonian systems on a firm categorical footing. Although constraints of time have prevented me from discussing nonequilibrium thermodynamics proper, and my thesis is split between equilibrium thermodynamics and non-equilibrium reversible systems, it is my hope that the techniques, philosophy, and results established here will lay the groundwork for a deeper understanding of open systems in thermodynamics, and additionally perhaps the groundwork for compositional computational simulation of thermodynamical systems for purposes of engineering and control.

Finally, a thesis is a wonderful opportunity not just to expand the frontiers of human knowledge, but also to invite others on the journey to places farther out. Thus, I hope to share with the reader some of the perspective that I have gained while studying this material, and to bring the reader up to date with some of the exciting developments in applied category theory that are the "shoulders of giants" on which I stand.

**Acknowledgments.**

First and foremost, I would like to thank my advisers John Baez and Wioletta Ruszel for supporting me through the process of writing this thesis and stimulating my interest in thermodynamics. I would especially like to thank John Baez for taking me on as his student when he had no obligation to do so, as he had no affiliation with Universiteit Utrecht, and for his careful reading and editing of this work.

My ingress to the world of applied category theory was through the Algebraic Julia project, and all of the contributors to Algebraic Julia have been a wonderful community who has inspired me in many ways. If you have committed any code to a repository under the umbrella of Algebraic Julia, then I am in debt to you for this thesis.

More specifically, I would like to thank Evan Patterson and James Fairbanks for being great mentors and for providing such a grand vision for Catlab and the future of categorical computing. Also, Sophie Libkind opened my mind to the possibilities of operads and operad algebras, putting me on the track to do most of the work in this thesis.

The germ of this thesis was created under the wings of the Topos Institute, while I did a summer internship there, and so I thank Brendan Fong, David Spivak, and the rest of Topos for providing such a wonderful place to do category theory. I would also like to thank the other people that showed up for Summer 2021 for stimulating discussion: Christian Williams, David Jaz Myers, Joshua Myers, Nelson Niu, and Sophie Libkind (again). Finally, I would like to thank the administrative staff at Topos, Nancy Derbish and Juliet Szatko, for the essential work of keeping the gears turning and the lights on.

Markus Lohmayer taught me a great deal about port-Hamiltonian systems, and also read a draft of this thesis in detail and found many embarrassing mistakes which I am very glad to fix.

My undergraduate thesis adviser, Yuri Sulyma, helped me kick off my journey into applied category theory (even though as a homotopy theorist he might have preferred I stuck with pure math), and so I thank him for that.

More personally, I would like to thank my parents, Mary Lee Potter and Theresa Lynch, for supporting (and funding!) my education, and for a great deal of patience on a wide variety of subjects. My sister Claire Lynch I thank for inspiring me, when I get somewhat worn out on environmental issues, to keep the fight alive.

I would also like to thank my friends. Sophie Galowitz got me into illustration on the computer, and many of the figures in the thesis owe her a debt. Emily Barker was patient with my mathematician's takes on chemistry, and I would not have been able to figure out how chemical reactions in thermostatics worked without her. Christiaan van den Brink and Ariane Blok brightened my time in Utrecht when it was darkened by the pandemic. And finally, Taeer Bar-Yam originally got me into category theory, and Haskell, and graphical notation, and a great number of other things I don't think I could really do justice in this acknowledgements section, so I'll leave it at that.

**Funding Statement.** This work was financially supported by an Utrecht Excellence Scholarship, the Topos Institute, and DARPA Award HR001120900067. The views and conclusions contained in this document are those of the authors and should not be interpreted as representing the official policies, either expressed or implied, of the Army Research Office or the U.S. Government. The U.S. Government is authorized to reproduce and distribute reprints for Government purposes notwithstanding any copyright notation herein.

# Notation

We use different fonts to signify various meanings in this thesis. **Bold** font is always used for terms that are being defined rigorously for the first time. *Italics* is used for emphasis, but also for terms that we are referring to before they are rigorously defined. This is so that the reader need not worry they have missed something if it is in italics; it will be defined in due time.

Within math, we use sans serif for the names of categories and serif font for the names of functors.

# Table of contents

# 1. Introduction

## 1.1. Relations and behaviors

> *We view a mathematical model as an exclusion law. A mathematical model expresses the opinion that some things can happen, are possible, while others cannot, are declared impossible. Thus Kepler claims that planetary orbits that do not satisfy his three famous laws are impossible. In particular, he judges nonelliptical orbits as unphysical. The second law of thermodynamics limits the transformation of heat into mechanical work. Certain combinations of heat, work, and temperature histories are declared to be impossible. Economic production functions tell us that certain amounts of raw materials, capital, and labor are needed in order to manufacture a finished product: it prohibits the creation of finished products unless the required resources are available.*
>
> J.W. Polderman and J.C. Willems, *Introduction to Mathematical Systems Theory*, [4]

This thesis is the product of many different intellectual traditions. One of the most prominent, however, is the philosophy that physical systems should be understood not as input-output machines, but rather as objects that have a collection of allowed behaviors. Composition of systems constrains their joint behavior. This is a philosophy that Jan Willems was particularly famous for, and is exposited in many of his works, but particularly in *The Behavioral Approach to Open and Interconnected Systems* [5].

In this thesis, there is a great deal of math, and a specifically a great deal of category theory. But fundamentally, the goal can be stated quite simply: to rigorously develop the Willems point of view. This is obviously a massive undertaking, and thus we narrow our scope to two areas. The first is the area of port-Hamiltonian systems, which, roughly speaking, are dynamical systems that have a notion of energy conservation. The second is the area of thermostatic systems, which is a formalization of thermodynamics that only considers equilibrium.

Without yet getting into any of the math, we present the basic ontology that we use to model systems in this way, pictured in Figure 1.1.

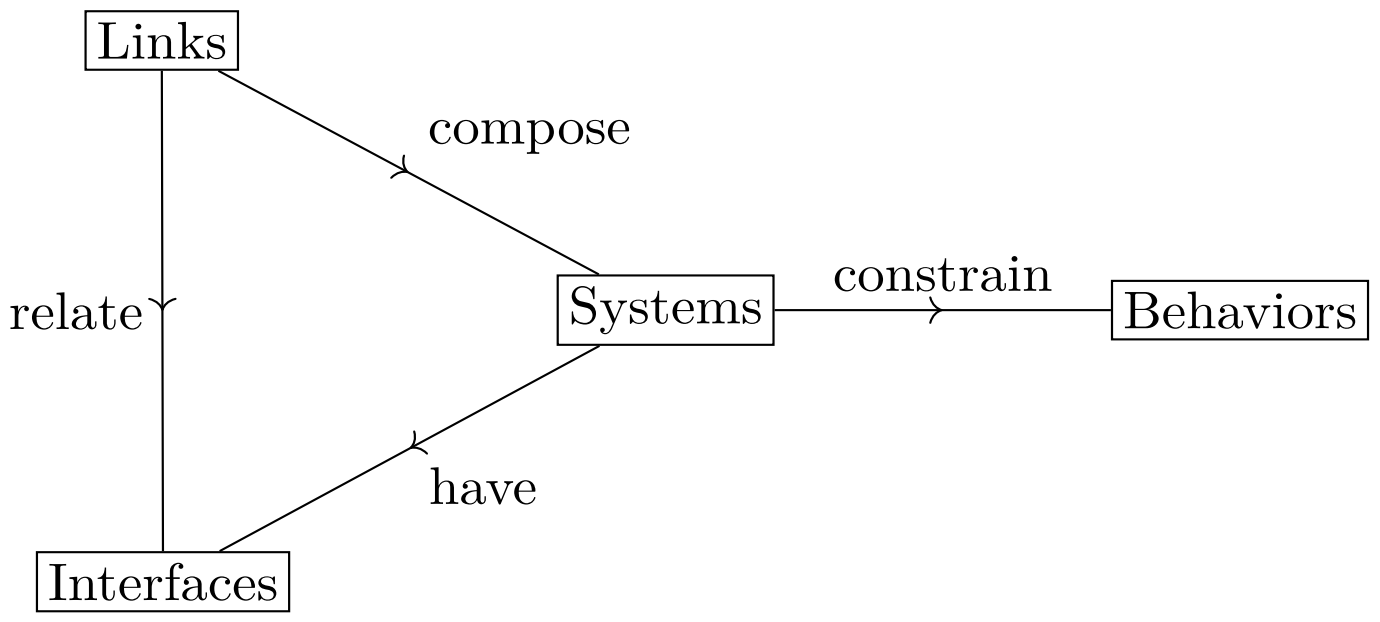

**Figure 1.1.** The ontology of relational composition.

As seen in the figure, our ontology consists of four parts.

1. The central concept is **systems**. In a traditional modeling framework, this would be the *only* concept. These are models of physical phenomena.

2. Systems constrain **behaviors**. In the dynamical setting, these are maps from a time interval to a state space, in the equilibrium setting these are simply elements of a state space. A system only allows certain behaviors.

3. Systems have **interfaces**. Systems interact with other systems through their interfaces, and interfaces allow us to "black box" systems so that all we care about is how they relate to their interface.

4. **Links** relate interfaces, and by relating interfaces, allow the composition of systems. In other words, links mediate how systems interact with one another through their interfaces.

This ontology draws not only from Willems, but also from a variety of sources (some of which were probably themselves inspired by Willems). These sources have provided formalizations for various terms in our ontology.

For instance, there is a thread of research into categorical formalization and application of linear relations that has greatly influenced our treatment of interfaces and links. This thread includes Baez and Erbele [6], Fong [7], and Coya [8].

Our philosophy on the composition of systems via operad algebras that allows us to formalize the statement that "links compose systems", comes mainly from Libkind et al. [9], which itself draws from a long history of using operad algebras to compose systems starting with Spivak [10]. And of course, we also have learned something of this from writing our own paper [1].

A third thread about how to think about behaviors, and specifically time-dependent behaviors, comes from "the Davids" (i.e. Spivak and Myers), and Fong, from [11], [12], and [13], though this thread is unfortunately not covered very much here.

Finally, the meat of this thesis comes from the theory of port-Hamiltonian systems, and this is the source of our "systems". We are indebted to the excellent book written by van der Schaft and Jelsema [14] for our understanding of this subject, and in more generality, to the work of van der Schaft, which has gone much farther than this thesis in many directions, including with respect to thermodynamics, and we hope to catch up to this with the category theory eventually. Port-Hamiltonian systems themselves have an illustrious history reaching back in many directions. We encourage the interested reader to read the introduction to [14] for more details.

The work in this thesis is mainly that of integration. This integration is not easy, especially when each thread requires a good deal of background material to get started. Therefore, in the introduction, we aim to give a preview of the "punchline" of each of the two main areas that we cover: port-Hamiltonian systems and thermostatics, to intrigue the reader and motivate them to plow through the background material. The next two sections are devoted to doing this. In the last section of the introduction, we give an overview of the structure of the rest of the thesis.

## 1.2. Power is the product of effort and flow

The "punchline" of port-Hamiltonian systems has to do with power, effort and flow, and this can all be explained in a fairly elementary way.

Suppose that a physical system has phase space $\mathbb{R}^n$ and $H\colon \mathbb{R}^n \to \mathbb{R}$ is a smooth function that gives the energy of each point in the state space. Moreover, suppose that $\gamma\colon I \to \mathbb{R}^n$ is a path of the physical system through time. Then an application of the chain rule gives

$$\frac{\mathrm{d}}{\mathrm{d}t} H(\gamma(t)) = \nabla H(\gamma(t)) \cdot \gamma'(t)$$

Although this is a simple trick of calculus, we would like to tell a story about this equation that give it more physical meaning. We start by giving examples.

**Example 1.1.** Consider a capacitor, which has phase space $\mathbb{R}$ with coordinate $q$ giving the charge on a plate. Then the energy $H$ is given by $H(q) = \frac{1}{2C} q^2$, where $C$ is the capacitance of the capacitor.

If $q = q(t)$ is a function of time, then the formula for the derivative of $H$ is given by

$$\dot{H} = \frac{1}{C} q \dot{q}$$

The right hand side of this equation is very familiar: $V = \frac{1}{C} q$ is the voltage across the capacitor, and $I = \dot{q}$ is the current flowing into the capacitor. We can thus write

$$\dot{H} = VI$$

That is, the change in the energy of the capacitor is given by the product of the voltage across it and current flowing into it.

We call these two parts of the expression for the derivative of energy “effort” and “flow”: effort is the gradient of the energy function and flow is the derivative of the state function. In this case, voltage takes the role of effort, and current takes the role of flow.

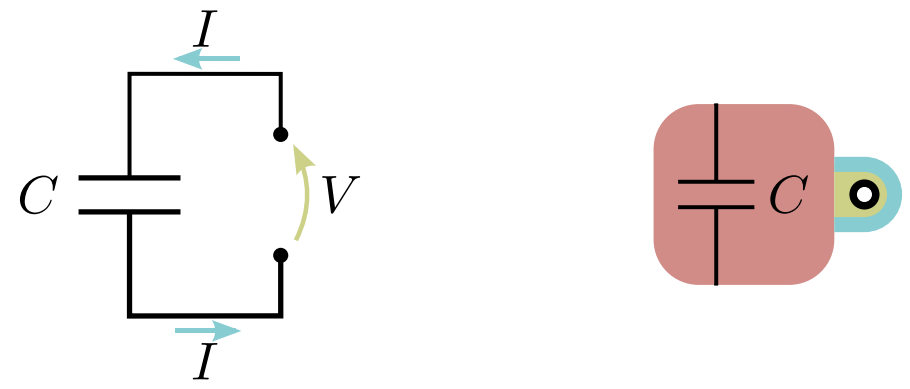

(a) Traditional (b) Dirac diagrams

**Figure 1.2.** Capacitors

In Figure 1.2, we see a capacitor drawn in a traditional style, as well as in the style of Dirac diagrams (original to this thesis). We eventually discuss Dirac diagrams more formally, but for now the reader should parse this diagram as a “system” (red box) with a “power port” that has an effort part (green) and a flow part (blue), which happen to represent voltage and current here. (However, green and blue are not always be voltage and current; in general they represent any pair of quantities whose product is power).

There are many physical systems that are expressed mathematically in the exact same way as a capacitor; see Table 1.1 for some examples. Note that for the systems that can be multi-dimensional (i.e. the particle or the gyroscope), here we are just considering the 1-dimensional variant (a particle constrained to a line or a gyroscope that can only spin on a single axis).

| System | State variable | Constant | Effort | Flow |
|---|---|---|---|---|
| Capacitor | Charge ($q$) | Capacitance ($C$) | Voltage ($V = \frac{1}{C} q$) | Current ($I = \dot{q}$) |
| Lin. Spring | Length ($x$) | Complicance ($C$) | Force ($F = \frac{1}{C} x$) | Velocity ($v = \dot{x}$) |
| Coil Spring | Rotation ($\theta$) | Ang. Compliance ($C$) | Torque ($\tau = \frac{1}{C} \theta$) | Ang. vel. ($\omega = \dot{\theta}$) |
| Inductor | Magnetix flux ($\Phi_B$) | Inductance ($L$) | Current ($I = \frac{1}{L} \Phi_B$) | Voltage ($V = \dot{\Phi}_B$) |
| Particle | Momentum ($p$) | Mass ($m$) | Velocity $\left(v = \frac{1}{m} p\right)$ | Force ($F = \dot{p}$) |
| Gyroscope | Ang. mmtm. ($L$) | Moment of Inertia ($I$) | Ang. vel. ($\omega = \frac{1}{I} L$) | Torque ($\tau = \dot{L}$) |

**Table 1.1.** Systems analogous to capacitors

One confusing thing about Table 1.1 is that in the context of a capacitor, we call voltage effort and current flow, but in the context of an inductor, we call voltage flow and current effort. This is similar to the phenomenon where if Mr. A is standing face-to-face with Ms. B, Mr. A will call “right” the direction that Ms. B calls “left”. In theater, this problem is dealt with by making left and right into immutable directions: “stage-left” and “stage-right”, which are left and right from the perspective of an actor on the stage facing the audience. In the context of bond graphs, which is where this effort and flow terminology comes from, there is a similar convention, where effort and flow are fixed to the “point of view” of a given system in each domain.

The fact that there isn’t a natural choice for which quantity to be “effort” and which quantity to be “flow” in each domain is related to the “Firestone analogy” in electrical circuits, which shows that there is more than one way of analogizing electrical circuits to mechanical systems [15].

The “canonical systems” for the electrical, linear mechanical, and rotational mechanical domains respectively are the capacitor, the linear spring, and the coil spring. The conventions for a variety of domains is listed in Table 1.2.

| Domain | Effort | Flow |
|---|---|---|
| Electrical | Voltage | Current |
| Linear Mechanical | Force | Velocity |
| Rotational Mechanical | Torque | Angular Velocity |
| Hydraulic | Pressure | Volume flow |
| Thermal | Temperature | Entropy flow |
| Chemical | Chemical potential | Mass flow |

**Table 1.2.** Effort and flow in various domains

**Example 1.2.** Consider a free particle of mass $m$ traveling in $\mathbb{R}^3$. The phase space for this particle is $\mathbb{R}^3 \times \mathbb{R}^3$ (position and momentum), and the energy function is $H(q,p) = \frac{1}{2m}\|p\|^2$. If we are only interested in changes of momentum, we can throw away the position, and end up with a phase space of $\mathbb{R}^3$ and an energy function of $H(p) = \frac{1}{2m}\|p\|^2$. Now, suppose that $p$ is a function of time. We end up with the equation

$$\dot{H} = \frac{1}{m} p \cdot \dot{p}$$

We can identify $\frac{1}{m}p$ with the velocity $v$ of this system, and by Newton's second law, $F = \dot{q}$. Thus, we can rewrite this equation as

$$\dot{H} = v \cdot F$$

One way to think about this is that the velocity of a particle measures the resistance of the particle to adding more momentum: it takes more energy to add momentum to something that is already going fast.

When we move to the multidimensional setting, a funny thing happens: we can have non-zero velocity *and* non-zero force, but still no net-power! This is the case when you have a pendulum: the rod of the pendulum exerts force on the mass attached to the pendulum, but this force does no work on the mass (although gravity still does work on the mass).

In the previous two examples, the effort changed based on the state of the system. But this does not always need to be the case.

**Example 1.3.** One can design a "gravitational battery" that stores energy by lifting a bunch of bricks with a crane, and releases energy by letting them fall back to the ground. The state variable here is the height of the bricks, and the energy stored is given by $H(h) = m\,g\,h$, where $m$ is the mass of the bricks and $g$ is the gravity on Earth. Then the flow is $v = \dot{h}$, the velocity that the crane is lifting at, and the effort is simply $F = m\,g$.

Finally, classical thermodynamics can also be thought about in this framework of efforts and flows.

**Example 1.4.** Consider an ideal gas. The fundamental thermodynamic equation is typically written as

$$\mathrm{d}U = T\,\mathrm{d}S - P\,\mathrm{d}V + \mu\,\mathrm{d}N$$

where $U$ is energy, $T$ is temperature, $S$ is entropy, $P$ is pressure, $V$ is volume, $\mu$ is chemical potential, and $N$ is particle number. If $U, T, S, P, V, \mu, N$ are given as functions of time, we can rewrite this as

$$\dot{U} = (T, -P, \mu) \cdot (\dot{S}, \dot{V}, \dot{N})$$

which is of the same form as our earlier examples. Thus, we see that in thermodynamics, intensive variables show up as *efforts* and the derivatives of extensive variables show up as *flows*.

In all of the previous examples, we have considered efforts and flow as they relate to energy functions. However, efforts and flows can exist independently of energy functions, as in the next example.

**Example 1.5.** Consider a network of ideal wires with a collection $P$ of inputs and a collection $Q$ of outputs. Let $\varphi_P, I_P \in \mathbb{R}^P$ be the potentials and currents assigned to the inputs, and $\varphi_Q, I_Q \in \mathbb{R}^Q$ be the potentials and currents assigned to the outputs. Then because no power is lost from current flowing across ideal wires, these potentials and currents must satisfy

$$\varphi_P \cdot I_P = \varphi_Q \cdot I_Q$$

**Example 1.6.** Consider a frictionless pulley system that multiplies the distance pulled by $\alpha$. Then if $v_{\text{in}}$ and $v_{\text{out}}$ are the velocities of the input and output rope respectively, and $F_{\text{in}}$ and $F_{\text{out}}$ are the forces applied to these ropes, we have $v_{\text{out}} = \alpha\, v_{\text{in}}$ and $F_{\text{in}} = \alpha\, F_{\text{out}}$, so that

$$F_{\text{in}}\, v_{\text{in}} = F_{\text{out}}\, v_{\text{out}}$$

It is the work of this thesis to describe a framework in which all of these examples can be discussed formally.

## 1.3. Equilibrium via entropy maximization

The other ground-level subject in this thesis is what we call "thermostatics"; the study of the equilibria of thermodynamical systems. Classically, the first two laws of thermodynamics are

1. Energy is conserved
2. Entropy never decreases

These laws are by no means rigorous, but way of introduction to thermostatics, let us investigate how these laws apply to a very simple scenario. Suppose two thermal bodies are allowed to exchange thermal energy with each other, but are insulated from the environment. For instance, suppose one were to put a hot pie and cold ice cream in an insulated cooler. What might one suppose should happen?

Thermodynamically, each body is described by an energy $U_i$, a temperature $T_i$, and an entropy $S_i$, where $S_i$ is a differentiable function of $U_i$ and

$$\frac{1}{T_i} = \frac{\partial S_i}{\partial U_i}$$

A very simple choice for the function $S_i$ is

$$S_i(U_i) = C_i \log U_i$$

This is because it leads to the equation

$$C_i\, T_i = U_i$$

which describes a system with a constant heat capacity $C_i$; each unit of temperature increase requires $C_i$ units of energy.

Now that we have the basic setup, let us apply the first and second laws. The first law, combined with the assumption that this system is insulated, implies that the total energy $U_1 + U_2$ should not change. Thus, we can only move within the solutions to the equation $U_1 + U_2 = U$ for some fixed $U$. The second law implies that the total entropy $S_1(U_1) + S_2(U_2)$ must never decrease. Given this constraint, we can imagine that the total entropy increases until it hits a local maximum, and then that is the equilibrium.

Now, in fact $S_1(U_1) + S_2(U_2)$ is a *concave* function, so local maximums are global maximums. Thus, we might reasonable assume that the equilbrium is at the global maximum of entropy, with respect to the constraint that energy remains conserved. Let us find this maximum.

This system only has one degree of freedom because of the $U_1 + U_2 = U$ constraint. Thus, we reparameterize it with $U_2 = U - U_1$, and we get

$$S_{\text{tot}}(U_1) = S_1(U_1) + S_2(U - U_1)$$

We now set the derivative of this to 0 in order to find the local maximum, and we get

$$S_1'(U_1) - S_2'(U - U_1) = 0$$

We can rewrite this as

$$\frac{\partial S_1}{\partial U_1} = \frac{\partial S_2}{\partial U_2}$$

Then, using our equation from before, we find

$$\frac{1}{T_1} = \frac{1}{T_2}$$

Thus, at the point of equilibrium, the two temperatures are equal to each other. This is precisely what we would expect; we would expect that after a long time, the ice cream and the pie would roughly come to the same temperature. Note that this is true for any concave differentiable functions $S_i$; we have not used our assumption that heat capacity is constant yet.

Now, suppose that we considered both thermal bodies as a single system, for instance by assuming that the time that it takes for the thermal bodies to come to equilibrium with each other is small compared to the time it takes for them to come to equilibrium with external systems. Then, we can parameterize the thermal bodies by their total energy $U$. It then turns out that we can write down an explicit form for their entropy (when at equilibrium with each other) as a function of $U$, and it is

$$S(U) = (C_1 + C_2) \log U + K$$

where $K$ is a constant not depending on $U$. Thus, the two bodies together act as a single body with heat capacity $C_1 + C_2$, which is also quite intuitive.

Let us now look over what we have done in this section, and see if we can extract a general pattern. We started by defining state variables for each system, which in this case were just energy. Then we wrote down entropy for each of these systems as a function of the state variables. Built-in to the entropy function were constants parameterizing the system (i.e. heat capacity). Then the entropy function captured the relationship between the state variables and other variables of the system (i.e. temperature), via the derivative of the entropy function.

We then wrote down a *constraint* on the systems, derived from a *conservation law*. This time we used conservation of energy, but we could also have more constraints based on conservation of matter, or conservation of other quantities. Once we did this, we could derive the equilibria by maximizing entropy with respect to the constraints. Finally, we used all of this setup to think about the *composed* system made out of parts that had equilibriated with eachother.

We formalize this whole process that we have just outlined, using a similar categorical method to port-Hamiltonian systems. Along the way, we also find that our process applies to more than just traditional classical thermodynamical systems, and we end up with a framework where classical, statistical, and quantum systems can all live under the same roof.

As of yet, we do not have a direct connection between the port-Hamiltonian systems and the thermostatic systems. However, we hope that presenting both illustrates that our categorical approach is quite broad, and thus can be lifted to many different types of systems. And hopefully the similarity of the approaches will make it possible for one day a connection to be found.

## 1.4. Overview

As said before, the last part of this introduction lays out the structure of this thesis, and serves as a guide for what is to come. The basic structure is given in Figure 1.3.

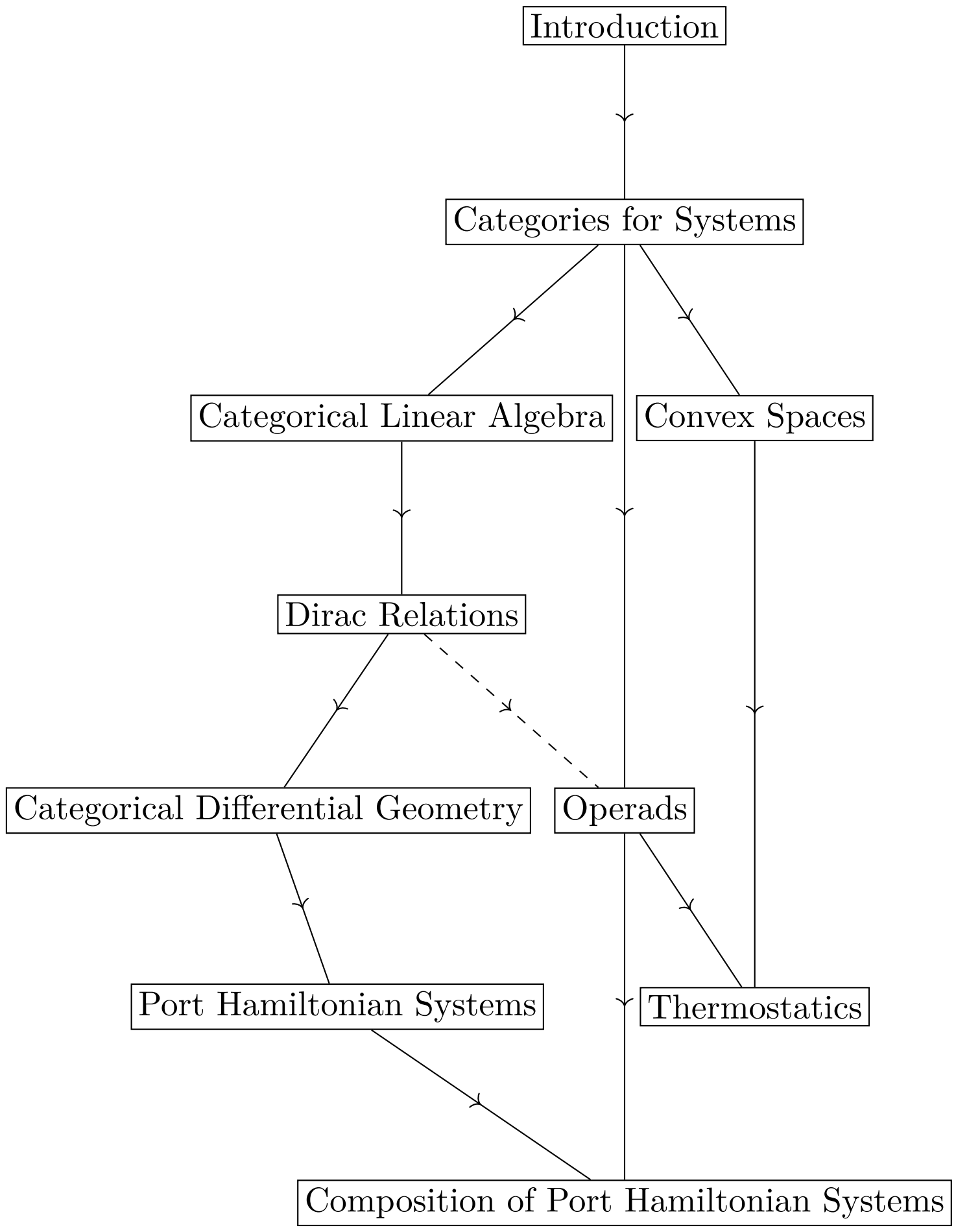

**Figure 1.3.** Dependencies between the chapters of the thesis. Dotted arrows mean that a chapter uses another chapter for examples, but not for the main logical development.

After this introduction, we start with Chapter 2, which is more or less pure category theory. In this chaper, we cover the principle of equivalence, monoidal category theory, and regular category theory. This chapter forms the basis for all that is to come in the thesis, and is purely review.

We then go on in Chapter 3 to cover linear algebra from a categorical viewpoint. There is a lot of rich structure in linear algebra that can be treated categorically; the category of vector spaces admits two interesting monoidal structures (direct product and tensor product), and additionally the category of linear relations is well-treated by category theory. Later on in the thesis, there are several categories that in spirit are quite similar to the category of linear relations, so this is a useful (and simple) preview. This chapter is likewise review.

In Chapter 4, we start treating some new material. Specifically, we introduce "Dirac relations", which formalize the idea of power-preserving interconnection. We also introduce "Dirac diagrams", which are a convenient and mnemonic way of picturing Dirac relations. We use Dirac diagrams and Dirac relations in order to compose our physical systems, in Chapter 8.

Next, Chapter 5 takes the material from Chapter 3 and 4 and "lifts" it to now be talking about vector bundles over manifolds. This is an essential generalization for the development of categorical port-Hamiltonian systems. For this "lifting" we use a theory of smooth functors, which is a neat construction of general interest. Most of this material is not original to this thesis, but the last section which lifts Dirac relations is original.

In Chapter 6, we get to port-Hamiltonian systems. We start with a review of classical mechanics, using Poisson manifolds, and then show how this classical mechanics viewpoint generalizes to port-Hamiltonian systems. Most of this material is known, but we also discuss morphisms between port-Hamiltonian systems, which to our knowledge is a concept original to this thesis.

Chapter 7 takes a break from the development of port-Hamiltonian systems and goes back to category theory; the only official prerequisite for this chapter is Chapter 2. This chapter introduces our formalism for composing systems: operads and their algebras. We use the operad of undirected wiring diagrams as an example, and also show how to make operads of Dirac relations. However, knowledge of Dirac relations is not necessary for this chapter, as they are used only in the examples. This chapter is wholly review of known material.

Finally, in Chapter 8, we prove that port-Hamiltonian system form an operad algebra of the operad of Dirac relations. This is the big result of this thesis, and is quite technical (though not at all surprising, and not too difficult). After proving this, we give examples of using it to compose port-Hamiltonian systems.

The rest of the thesis applies similar methods to thermostatic systems. This material appeared earlier in Baez, Lynch and Moeller [1], and so here we present a condensed version, with the proofs omitted. Our hope is that this version shows the essential simplicity of the thermostatic approach, and also demonstrates the versatility of the operadic framework. Chapter 9 does a fairly standalone treatment of convex spaces, which are the state spaces of thermostatic systems. Then Chapter 10 presents the operad algebra of thermostatic systems, along with several examples. Most of this material was already covered in [1], but our treatment of chemical reactions is novel to this thesis.

There are two broad "tracks" through this thesis. The first is the "port-Hamiltonian track." This track goes straight through from Chapter 2 to Chapter 8. The second is the "thermostatics track." This track is shorter, and goes 2, 9, 7, 10.

However, the "logical" path through the material is not always the most intuitive. The reader with strong grasp of differential geometry might skip directly to Chapter 6 to get an idea of what a port-Hamiltonian system is, before returning to Chapter 2 and doing the material in order. Additionally, it might make sense to look at some of Chapter 3 before doing Chapter 2.

We close the introduction by stressing that the core of this thesis is simply the juxtaposition of several well-understood theories, as discused in Section 1.1. Thus, we hope that the reader will take this thesis as an invitation to open the door for themselves, read the literature that we have drawn from, and apply these methods to their own favorite subject. It is perhaps not easy, but we believe that the fruit of these methods has just begun to be picked.

# 2. Categories for Systems

## 2.1. Principle of equivalence

> *Several considerable paradoxes follow from this, amongst others that it is never true that two substances are entirely alike, differing only in being two rather than one.*
>
> G.W. Leibniz, *Discourse on Metaphysics*, 1686, [16]

> *In folklore we find the following witty slogan: if there is a contradiction between physics and logic, you must change the logic.*
>
> René Lavendhomme, *Basic Concepts of Synthetic Differential Geometry*, 1996, [17]

The task of the mathematical physicist is to take what is intuitively obvious to physicists and render it rigorously clear to mathematicians. In this quest, the study of logic can be a great boon, providing a foundation of rock on which to build houses. However, sometimes it is necessary to remember that the physics must come first.

For instance, for a great many purposes, physicists found it useful to work with Dirac deltas and their ilk. The onus is then on the mathematician to loosen their definition of function as a simple map of sets and invent new theories to deal with this problem, such as measure theory, distribution theory, and Fourier theory. Of course, this can still be accomplished within classical foundations, but the point is that the fact that "function" is defined in a certain way within classical foundations does not mean that people who use function in a different way are wrong. A more thorough account of this story can be found in Jaynes [18, Appendix B].

Leibniz's "Identity of Indiscernibles", part of which is quoted above, states that two things which have all of the same properties in common are in fact one thing [19]. In mathematics, we use this idea when we use the definite article "the" to refer to "the real numbers." In set theory, there are any number of different sets which could be used for $\mathbb{R}$, and yet this fact never matters in the practice of mathematics. This is because the “properties” of all of these sets are the same. The rest of this section is devoted to discussion and formalization of this idea. This may seem like an embarrasing triviality, but as we delve into the categorical structures needed later on in this thesis a firm grasp on this bit of philosophy is very helpful.

The mathematical history of this formalization is long and illustrous; we can see that the idea goes all the way back to Leibniz. Tracking the entire story is beyond the scope of this section, but some modern takes on the subject are found in Awodey [20] and Coquand [21]. The author’s thinking on this matter was particularly influenced by Riehl’s talk [22], and more generally the homotopy type theory project [23].

Formalization of the identity of indiscernibles hinges on the question of what “property” means in mathematics? In set theory, there is a very specifc form of identity of indiscernibles which says that if two sets have all of their elements in common, then they are equal. However, this is no good: we could have two sets that had *no* elements in common, but each could serve just fine as $\mathbb{R}$. In set theory, we also have a concept of bijection. But this is also no good: $\mathbb{R}$ and $\mathbb{C}$ have a bijection because they are of the same cardinality, yet intuitively these sets have quite different properties!

The answer is that we must specify which properties we care about, and only consider bijections that preserve these properties. For instance, the reals are a complete ordered field, and we can show that there is a canonical field isomorphism between any two complete ordered fields [24]. This notion of “isomorphism” is saying that all the relevant properties of two objects are the same.

The process of taking this answer seriously leads to category theory. Of necessity, we assume that the reader has had some exposure to category theory before, as a full tutorial on category theory is beyond the scope of this thesis. Works that the reader might consult for introduction to the subject include Leinster [25], Riehl [26] and Maclane [27]. However, we reserve the right to refresh the reader’s memory on certain definitions for the sake of pedagogy.

Definition 2.1. *Two objects $X, Y$ in a category $\mathsf{C}$ are said to be isomorphic if there exists $f\colon X \to Y$, $g\colon Y \to X$ such that $1_X = g \circ f$ and $1_Y = f \circ g$*

In mathematics, we generally hold the following opinion.

**Principle of Equivalence.** All mathematical constructions should be invariant under isomorphism in the relevant category. That is, we should not notice if someone were to replace all of our objects with isomorphic objects.

However, when we try to apply the principle of equivalence to categories themselves, we end up with something that breaks the principle of equivalence! That is, when we assert that $G \circ F = 1_{\mathsf{C}}$ for functors $F\colon \mathsf{C} \to \mathsf{D}$ and $G\colon \mathsf{D} \to \mathsf{C}$, we are asserting that $G(F(X)) = X$, which is not invariant under isomorphism. Thus, we make the following definition.

Definition 2.2. *An **equivalence** of categories $\mathsf{C}$ and $\mathsf{D}$ consists of*

1. *A pair of functors $F\colon \mathsf{C} \to \mathsf{D}$, $G\colon \mathsf{D} \to \mathsf{C}$*
2. *A pair of natural isomorphisms $1_{\mathsf{C}} \cong G \circ F$, $1_{\mathsf{D}} \cong F \circ G$*

We then rephrase the principle of equivalence.

**Categorified Principle of Equivalence.** All categorical constructions should be invariant under equivalence of categories. That is, we should not notice if someone were to replace all of our categories with equivalent categories.

Equivalence can be surprisingly weak. Let $\mathsf{C}$ be the category with one object and one morphism (the identity), and let $\mathsf{D}$ be a category with $2^{\mathbb{N}}$ objects, but precisely one morphism between any two (these morphisms are then isomorphisms). Then $\mathsf{C}$ and $\mathsf{D}$ are equivalent. We say that a category is **contractible** if it is equivalent to the one object, one morphism category.

**Example 2.3.** The subcategory of the category of fields consisting of all complete, ordered fields is contractible.

If a category is contractible, according to the principle of equivalence, we might as well treat it like it has one object. This leads to the following interpretation of Leibniz’ Identity of Indiscernibles principle.

**Identity of Indiscernibles.** The English phrase "the object with property $P$" is understood to mean in a mathematical context that the category of objects that satisfy property $P$ along with $P$-preserving morphisms is contractible.

**Example 2.4.** If $\mathsf{C}$ is category, then we say that 1 is a **terminal object** of $\mathsf{C}$ if for all $X \in \mathsf{C}$, there exists a unique map $f: X \to 1$. Then the subcategory consisting of all terminal objects of $\mathsf{C}$ is contractible, because there is a unique map between any two terminal objects. Thus, we are justified by identity of indiscernibles in saying "the" terminal object.

The same can be said for **initial objects**, 0 is an initial object of $\mathsf{C}$ if for all $X \in \mathsf{C}$ there exists a unique map $f: 0 \to X$. The subcategory of initial objects is contractible, so we say "the" initial object.

The subtle part of this definition is the idea that not only are all of the objects isomorphic but also they are *canonically* isomorphic; i.e. we have exactly one isomorphism between any two objects. It would sound strange to say "the vector space of dimension $n$" because given two such vector spaces, there is not a canonical isomorphism between them. However, it is natural to say "the vector space $\mathbb{R}^n$", because $\mathbb{R}^n$ comes with a canonical ordered basis, so any two instances of $\mathbb{R}^n$ are canonically isomorphic.

Finally, note that we do not take this as dogmatic. If it is convenient to use a stricter form of equality than isomorphism, we do not hesitate to do so. It is just that we must be aware that strict equality should not necessarily be the default. Full application of the principle of equivalence takes us down a long road whose end is as of yet unknown, though many have travelled along it much farther than we do here and found a wide variety of interesting topics; see [23].

## 2.2. Monoidal categories

Monoidal categories are a foundational object in category theory and especially applied category theory. They were introduced independently in 1963 by Benabou [28] and Maclane [29]. Prominent early uses in applied category theory are [30] and [6], both of which are excellent introductions to the subject that the reader is encouraged to refer to as a complement to this section.

This section rigorously defines what a monoidal category is. However, the reader is encouraged to skip back and forth between this section and section 2.3 to develop intuition on how to think about monoidal categories.

Monoidal categories are a *categorification* of monoids [31]. Roughly speaking, categorification refers to the process of taking a definition with sets, and replacing the sets with categories. So we start our introduction to monoidal categories by reviewing what a monoid is.

DEFINITION 2.5. *A* ***monoid*** *consists of*

- *a set $M$*
- *a function $\mu: M \times M \to M$*
- *an element $e \in M$*

*such that*

- *$\mu(a, \mu(b, c)) = \mu(\mu(a, b), c)$ for all $a, b, c \in M$*
- *$\mu(e, a) = a = \mu(a, e)$ for all $a \in M$*

The most straightforward way of categorifying this is a strict monoidal category.

DEFINITION 2.6. *A* ***strict monoidal category*** *consist of*

- *a category $\mathsf{C}$*
- *a functor $\otimes: \mathsf{C} \times \mathsf{C} \to \mathsf{C}$*
- *an object $I \in \mathsf{C}$*

*such that*

- *$- \otimes (- \otimes -) = (- \otimes -) \otimes -$ as functors. That is, for all objects $X, Y, Z \in \mathsf{C}_0$,*

$$X \otimes (Y \otimes Z) = (X \otimes Y) \otimes Z$$

*and for all morphisms $f: X \to X'$, $g: Y \to Y'$, $h: Z \to Z'$*

$$f \otimes (g \otimes h) = (f \otimes g) \otimes h$$

- $(-\otimes I) = (-) = (I \otimes -)$ *as functors. That is, for all objects* $X \in \mathsf{C}_0$,

  $$X \otimes I = X = I \otimes X$$

  *and for all morphisms* $f\colon X \to X'$

  $$f \otimes 1_I = f = 1_I \otimes f$$

There are also additional laws "hidden" in the functorality of $\otimes$, namely the fact that for

$$X_1 \xrightarrow{f_1} Y_1 \xrightarrow{g_1} Z_1, X_2 \xrightarrow{f_2} Y_2 \xrightarrow{g_2} Z_2$$

we have

$$(g_1 \otimes g_2) \circ (f_1 \otimes f_2) = (g_1 \circ f_1) \otimes (g_2 \circ f_2)\colon X_1 \otimes X_2 \to Z_1 \otimes Z_2$$

and additionally

$$1_{X_1} \otimes 1_{X_2} = 1_{X_1 \otimes X_2}$$

**Example 2.7.** Let $\mathbb{N}$ be the category where the objects are natural numbers, and where a morphism $f \in \mathrm{Hom}_{\mathbb{N}}(m, n)$ is a function $f\colon \{1, \ldots, m\} \to \{1, \ldots, n\}$. (A category theorist might know this as the skeleton of $\mathsf{FinSet}$). Then $(\mathbb{N}, +, 0)$ forms a monoidal category, where $f + g$: $m_1 + m_2 \to n_1 + n_2$ is defined by

$$(f + g)(k) = \begin{cases} f(k) & \text{if } k \leq m_1 \\ g(k - m_1) + n_1 & \text{if } k > m_1 \end{cases}$$

An example of this monoidal composition of morphisms is pictured in Figure 2.1.

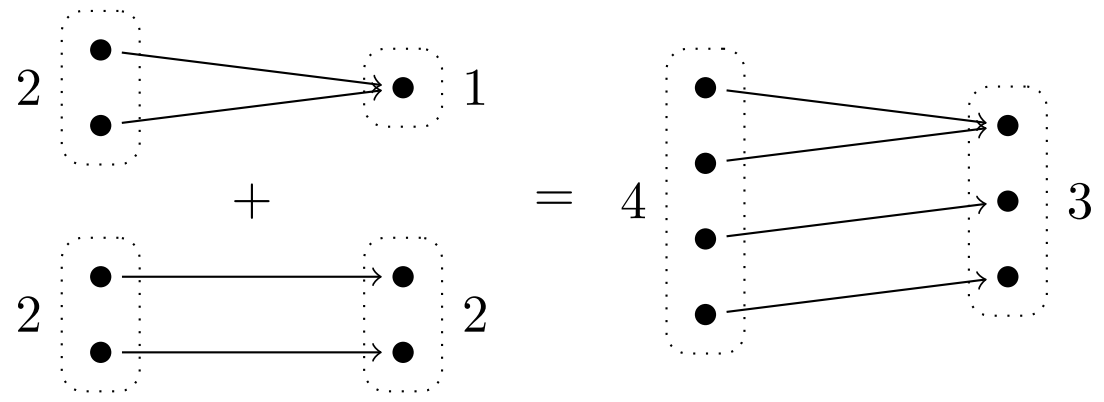

**Figure 2.1.** Picturing Monoidal Composition in $\mathbb{N}$

**Example 2.8.** Let $\mathbb{N}^2 = \mathbb{N} \times \mathbb{N}$ be the category where the objects are tuples $(n_1, n_2)$, and a morphism $f \in \mathrm{Hom}_{\mathbb{N}^2}((m_1, m_2), (n_1, n_2))$ is a tuple $(f_1\colon n_1 \to m_1, f_2\colon n_2 \to m_2)$.

We can think of this as a "two-colored" version of Example 2.7, as is pictured in Figure 2.2.

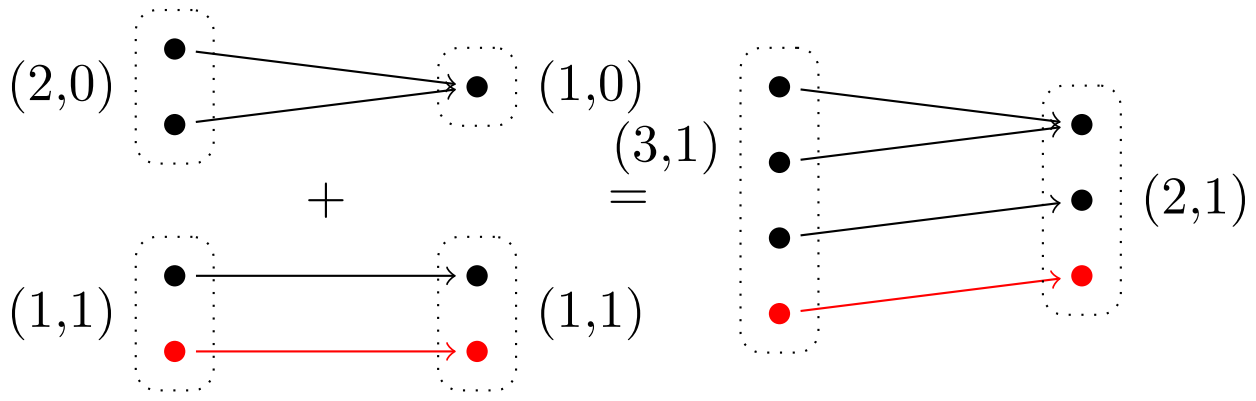

**Figure 2.2.** Picturing Monoidal Composition in $\mathbb{N}^2$

This can be generalized to $\mathbb{N}^S$, where $S$ is a category and $\mathbb{N}^S$ is the category of functors $F$: $S \to \mathbb{N}$. These categories are useful for modeling a variety of data structures, see [32] for details.

**Example 2.9.** If $\mathsf{C}$ is any category, then let $\mathrm{End}(\mathsf{C})$ be the category where the objects are functors $F\colon \mathsf{C} \to \mathsf{C}$ and the morphisms are natural transformations $\alpha\colon F \Rightarrow G$. Then $(\mathrm{End}(\mathsf{C}), \circ, 1_C)$ is a strict monoidal category.

Although strict monoidal categories are convenient (especially for computation), they violate the principle of equivalence by requiring that $X \otimes (Y \otimes Z) = (X \otimes Y) \otimes Z$; in accordance with the principle of equivalence we should instead have $X \otimes (Y \otimes Z) \cong (X \otimes Y) \otimes Z$. For instance, if

$$X \times Y = \{(x, y) \mid x \in X, y \in Y\}$$

then set-theoretically we do not have $X \times (Y \times Z) = (X \times Y) \times Z$. However there is a canonical isomorphism $X \times (Y \times Z) \cong (X \times Y) \times Z$.

Moreover, if $(\mathsf{C}, \otimes, I)$ is a strict monoidal category, and $\mathsf{C}$ is equivalent to another category $\mathsf{D}$, then $\mathsf{D}$ is not necessarily given a strict monoidal structure by this equivalence. Rather, $\mathsf{D}$ is given what we call a *monoidal structure*, defined as follows. In other words, monoidal structures are the generalization of strict monoidal structures that are invariant under equivalence of categories.

Definition 2.10. *A **monoidal category** is a triple $(\mathsf{C}, \otimes, I)$ of a category $\mathsf{C}$, a functor $\otimes\colon \mathsf{C} \times \mathsf{C} \to \mathsf{C}$, and an object $I \in \mathsf{C}_0$, along with a natural isomorphism called the **associator***

$$a_{X,Y,Z}\colon (X \otimes Y) \otimes Z \xrightarrow{\sim} X \otimes (Y \otimes Z)$$

*and two natural isomorphisms*

$$l_X\colon I \otimes X \xrightarrow{\sim} X$$

$$r_X\colon X \otimes I \xrightarrow{\sim} X$$

*called **left** and **right unitors**, such that the following coherence conditions hold:*

- *for all $X, Y \in C$ the **triangle equation**:*

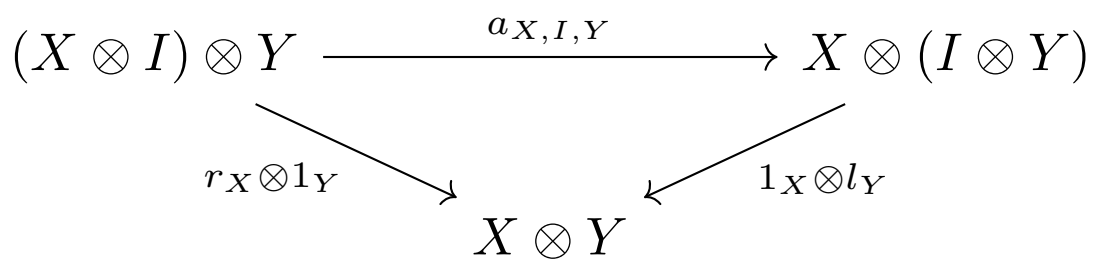

- *for all $W, X, Y, Z \in C$ the **pentagon equation**:*

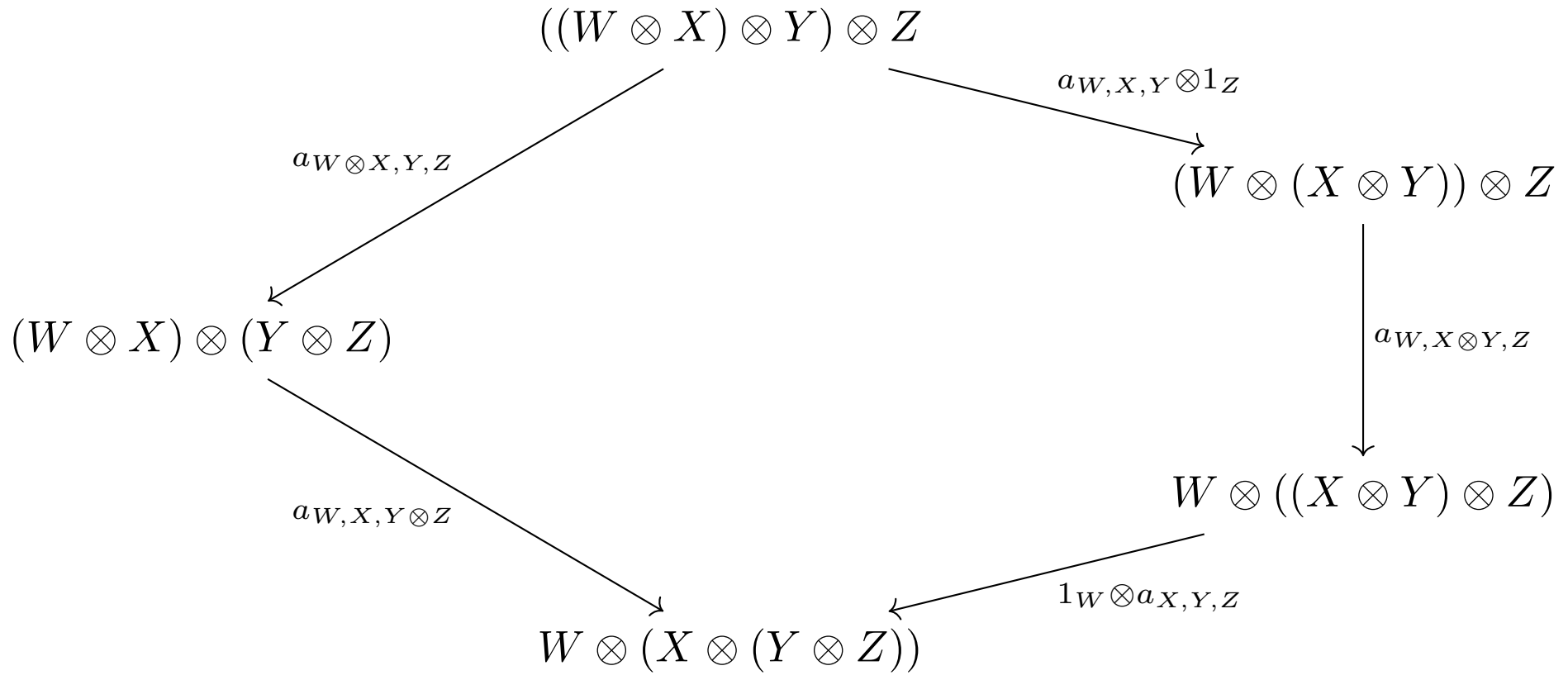

These coherence conditions can look intimidating at first. However, essentially what they are doing is providing sufficient conditions for:

- All of the parenthizations of $X_1 \otimes \cdots \otimes X_n$ to be *canonically* isomorphic. This allows us to speak of "the" object $X_1 \otimes \cdots \otimes X_n$.
- Any way we can insert or drop copies of $I$ into $X_1 \otimes \cdots \otimes X_n$ to form, for example $X_1 \otimes I \otimes X_2 \otimes X_3 \otimes I$, to be canonically isomorphic to just $X_1 \otimes \cdots \otimes X_n$ (along with any parenthization).

The fact that the two triangles and the pentagon are sufficient to prove this is a celebrated theorem of Mac Lane, found in [27].

**Example 2.11.** We would like the following sentence to be true. If $\mathsf{C}$ is a category with products $\times$ and a terminal object 1, then $(\mathsf{C}, \times, 1)$ is a monoidal category, and $(\times, 1)$ is called the **cartesian monoidal structure** on $\mathsf{C}$. However, there is a problem with this definition. $\mathsf{C}$ "having products" does not mean that we have chosen a specific object $X \times Y$ for every $X$ and $Y$; $X \times Y$ is only specified up to canonical isomorphism. There are various ways to get around this; one is by picking an arbitrary representative for each product.

Note that the fact that all parenthizations of $X_1 \times \cdots \times X_n$ are canonically isomorphic is simply a result of the universal property of the product of $n$ objects in a category.

**Example 2.12.** The same goes for if $\mathsf{C}$ is a category with coproducts $+$ and an initial object 0: $(\mathsf{C}, +, 0)$ is a monoidal category, and $(+, 0)$ is called the **cocartesian monoidal structure** on $\mathsf{C}$.

**Example 2.13.** The category FinSet of finite sets and functions between them has a cocartesian monoidal structure that behaves very similarly to $(\mathbb{N}, +, 0)$.

**Example 2.14.** The category of vector spaces and linear maps has a monoidal structure given by direct product $V \oplus W$, which is both the cartesian and cocartesian monoidal structure. We discuss this more in Chapter 3.

In all of our examples except for Example 2.9, it was true that $X \otimes Y \cong Y \otimes X$. The following definition formalizes this structure, and additionally gives coherence conditions.

Definition 2.15. *A* ***symmetric monoidal category*** *is a monoidal category* $(\mathsf{C}, \otimes, I)$ *along with a natural isomorphism* $B_{X,Y}: X \otimes Y \xrightarrow{\sim} Y \otimes X$*, called the* ***braiding****, such that*

- *The hexagon identity holds:*

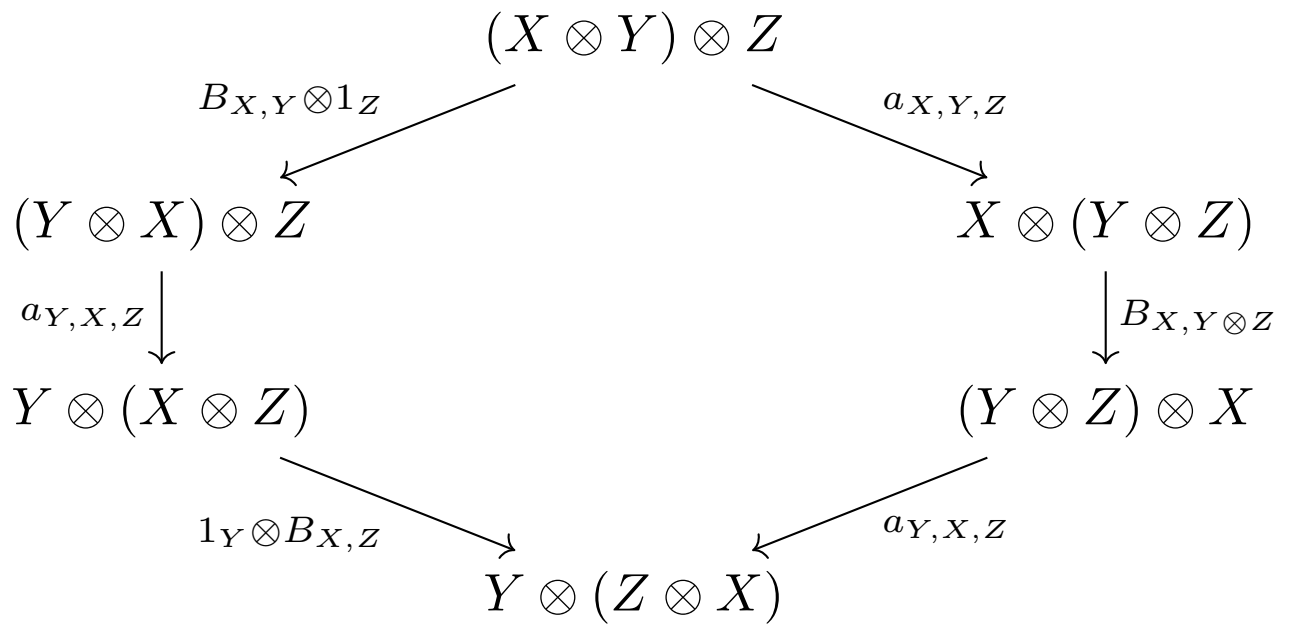

- $B_{Y,X} \circ B_{X,Y} = 1_{X \otimes Y}$

These coherence conditions are sufficient to show that there is a unique natural isomorphism $X_1 \otimes \cdots \otimes X_n \to X_{\sigma(1)} \otimes \cdots \otimes X_{\sigma(n)}$ built out of the braiding for any permutation $\sigma \in S(n)$.

**Example 2.16.** Any cartesian monoidal category is a symmetric monoidal category, as there is a natural isomorphism $B_{X,Y}: X \times Y \to Y \times X$ that is given by applying the universal property of products to the two projections $\pi_2: X \times Y \to Y$ and $\pi_1: X \times Y \to X$. Dually, any cocartesian monoidal category is also symmetric.

Finally, there is an even nicer class of symmetric monoidal categories that we often use: compact closed categories. We start with the definition of closed monoidal categories, originally due to Eilenberg and Kelly [33].

Definition 2.17. *A monoidal category* $(\mathsf{C}, \otimes, I)$ *is called a* ***closed monoidal category*** *if* $-\otimes Y$ *has a right adjoint* $[Y, -]$ *for every* $Y$*. We call* $[Y, Z]$ *the* ***hom-object*** *corresponding to* $Y$ *and* $Z$*.*

**Example 2.18.** In the cartesian monoidal category $(\mathsf{Set}, \times, 1)$, for two sets $X$ and $Y$ we define $[Y, Z]$ to be the set of functions $f: Y \to Z$. It is well-known that

$$\mathrm{Hom}(X, [Y, Z]) \cong \mathrm{Hom}(X \times Y, Z)$$

The process of taking a function $f: X \times Y \to Z$ and producing a function $f: X \to [Y, Z]$ is known as "currying" and the reverse process is known as "uncurrying".

Compact closed categories are closed monoidal categories where the hom-objects are computed in a particular way.

Definition 2.19. *A* ***dual object*** *for an object* $A$ *in a symmetric monoidal category* $(\mathsf{C}, \otimes, I)$ *is an object* $A^*$ *along with*

- *a morphism* $\mathrm{ev}_A: A^* \otimes A \to I$
- *a morphism* $i_A: I \to A \otimes A^*$

*such that the following diagram commutes:*

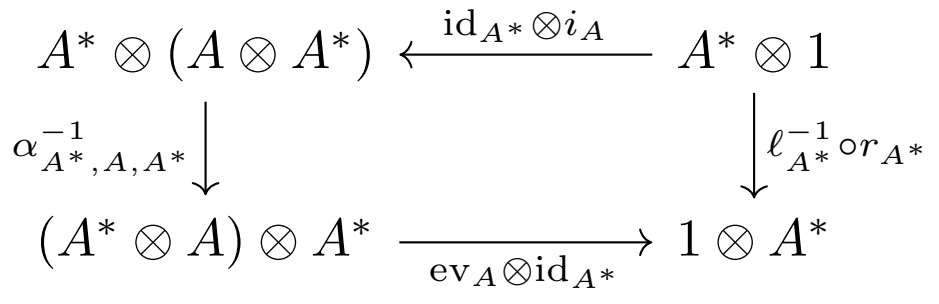

*See Ponto [34] for more details on dual objects.*

The diagram for dual objects is also known as the "zigzag" relation, for reasons that become apparent in the next chapter. In a general (i.e., not-necessarily symmetric) monoidal category, there are two diagrams that must commute; in the symmetric setting only one is necessary.

DEFINITION 2.20. *A* ***compact closed category*** *is a symmetric monoidal category where each object $A$ is equipped with a dual object $A^*$.*

It is not immediately clear from the above definition that a compact closed category is in fact even a closed monoidal category. This is the content of the next proposition, whose proof is delayed until the next section.

PROPOSITION 2.21. *A compact closed category is a closed monoidal category with hom-object $[Y, Z] = Z \otimes Y^*$.*

**Proof.** We must display a natural isomorphism

$$\mathrm{Hom}(X, Z \otimes Y^*) \cong \mathrm{Hom}(X \otimes Y, Z)$$

We prove this with string diagrams in the next section. □

## 2.3. String diagrams

String diagrams are a way of picturing morphisms in monoidal category in a two-dimensional syntax. Roughly speaking, the reason why we have a two-dimensional syntax for monoidal categories and a one-dimensional syntax

$$X \xrightarrow{f} Y \xrightarrow{g} Z$$

for regular categories is that for monoidal categories we have two types of composition for morphisms (i.e. categorical composition with $\circ$, and monoidal composition with $\otimes$), and we represent each type of composition with juxtaposition along a different axis.

The basic building block of string diagrams is picturing a morphism $f: X_1 \otimes \cdots \otimes X_n \to Y_1 \otimes \cdots \otimes Y_m$ as a bead with several wires coming in and out. So for instance, $f: X_1 \otimes X_2 \to Y_1 \otimes Y_2 \otimes Y_3$ would be pictured as

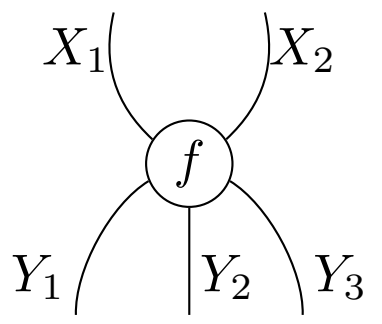

Note that the coherence condition for associativity is essential in interpreting this picture: it is not possible to parenthesize in string diagrams, so we need to be able to say "the" object $Y_1 \otimes Y_2 \otimes Y_3$.

We represent categorical composition of morphisms with vertical juxtaposition:

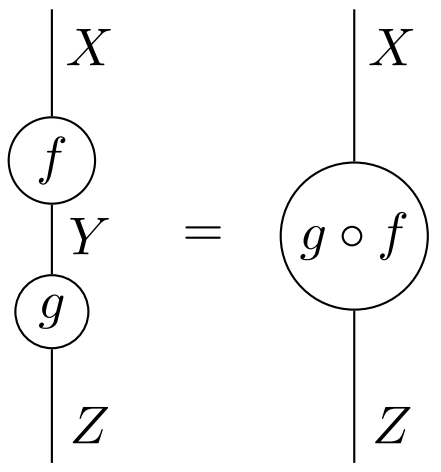

and we represent monoidal composition with horizontal juxtaposition:

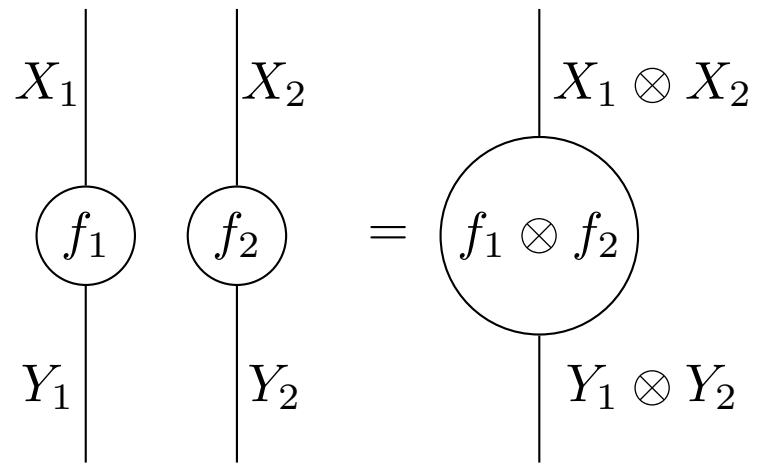

Notice that there is another "ambiguity in parenthesization" inherent in string diagrams; we can parse the following string diagram in two ways.

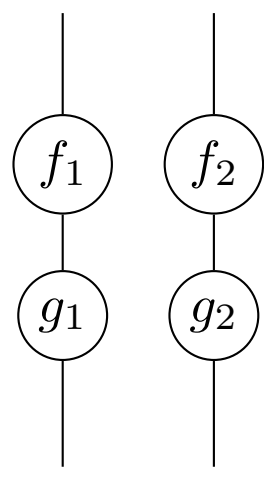

The first way is by letting $f = f_1 \otimes f_2$, $g = g_1 \otimes g_2$, and then taking

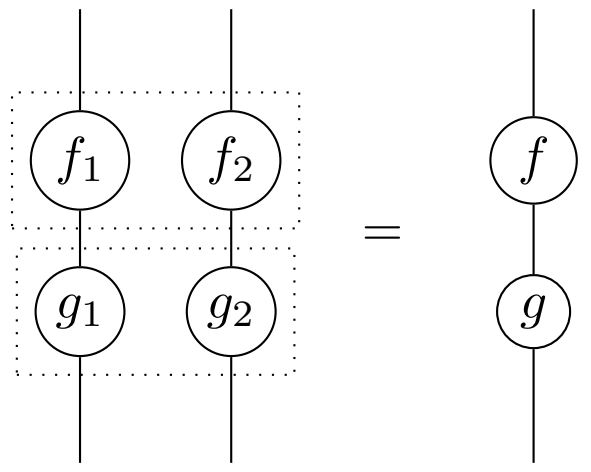

and the second way is by letting $h_1 = g_1 \circ f_1$, $h_2 = g_2 \circ f_2$, and then taking

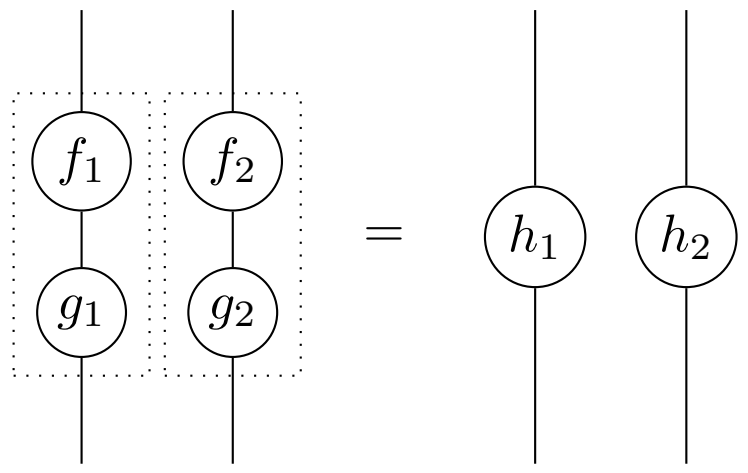

By the fact that $\otimes$ is a functor, these end up being equal, so it does not matter in what order we interpret the string diagrams.

Identities are particularly simple in string diagrams too. The categorical identity $1_X : X \to X$ is pictured as a plain string,

and the monoidal identity $I$ is pictured as a simple lack of a wire, so that a function $f : I \to X \otimes Y$ would be represented as

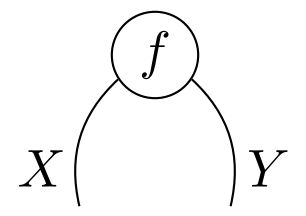

Thus, the laws of monoidal categories are made implicit in our notation; it would be impossible to interpret our notation unambiguously without those laws. For instance, without the boxes to indicate the categorical and monoidal identities, one would not be able to tell the following two diagrams apart.

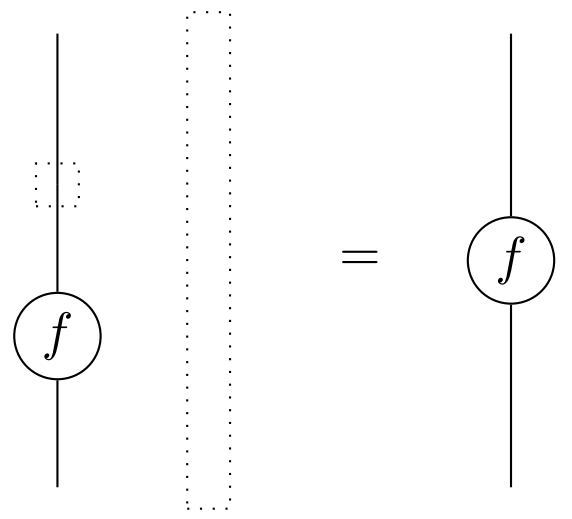

Finally, we represent the braiding isomorphism in symmetric monoidal categories by a simple crossing of wires.

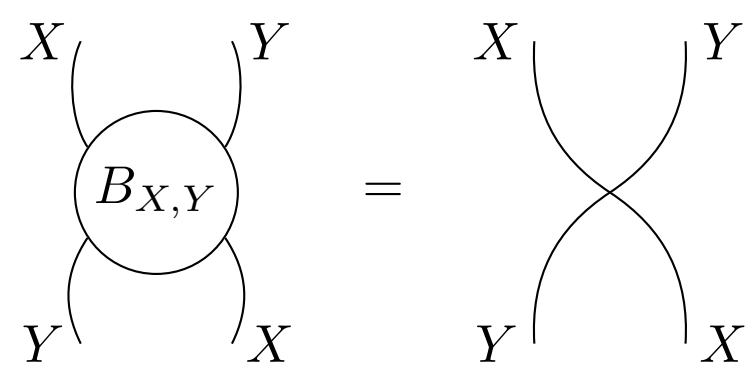

The coherence conditions for symmetric monoidal categories imply that we can treat the strings in string diagrams "topologically"; that is we can not worry so much about the exact pattern of string crossings, and only care where the strings start and end. For instance, coherence implies the following equality.

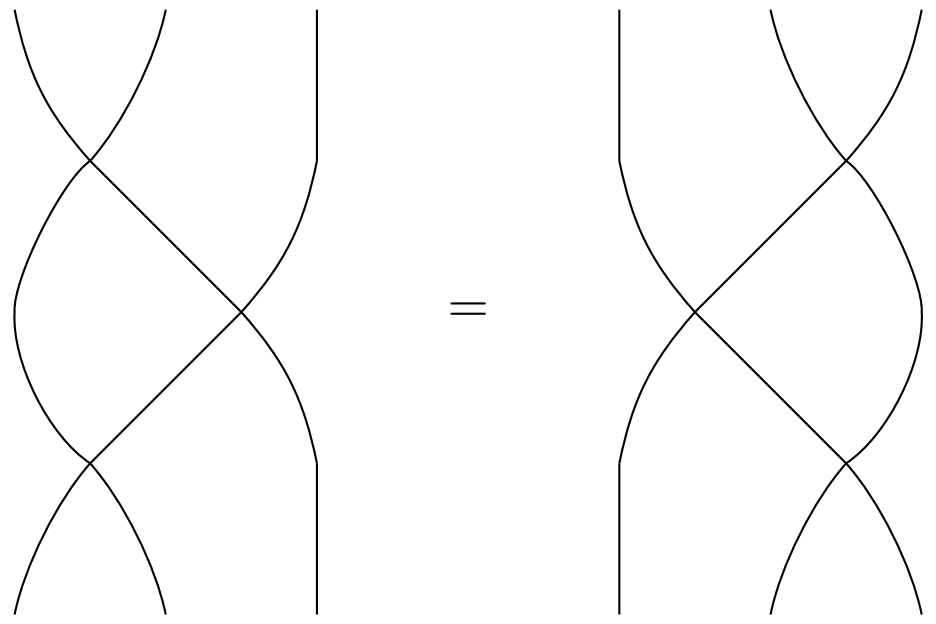

Finally, we can also represent the compact closed structure with string diagrams. In this representation, we draw $i_A\colon I \to A \otimes A^*$ with a downwards-facing curve:

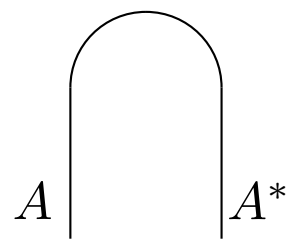

and $\mathrm{ev}_A\colon A^* \otimes A \to I$ with an upwards-facing curve:

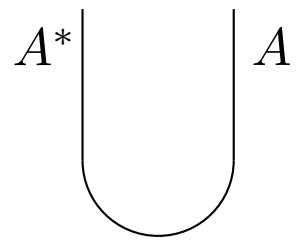

We can then rewrite the commutative diagram for duals as the following equality, which we think of as "pulling straight" the wires. One can now see where the term "zigzag" relation comes from.

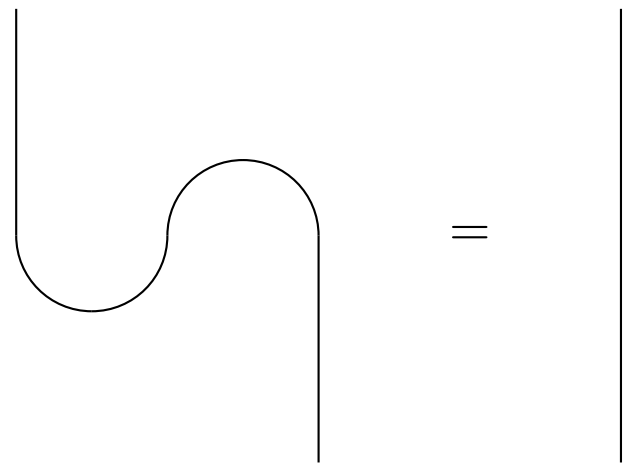

The identification of $[A, B]$ with $B \otimes A^*$ in a compact closed category is (morally speaking) captured in the following diagram, although technically the following diagram pictures $\mathrm{Hom}(A, B) \cong \mathrm{Hom}(I, B \otimes A^*)$.

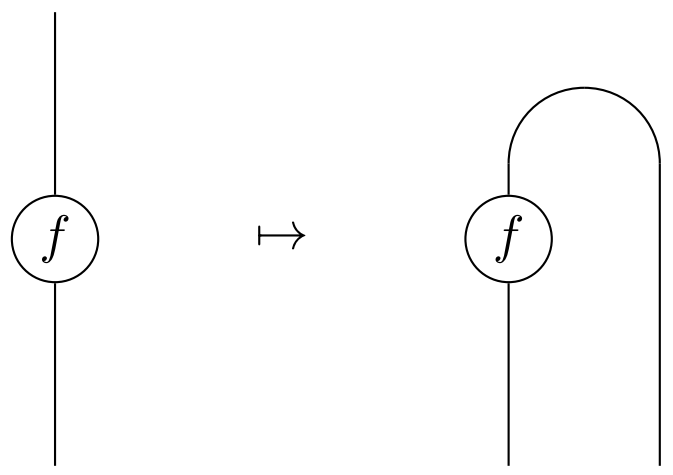

To actually show that $B \otimes A^*$ is a valid hom-object, we now prove Proposition 2.21.

**Proof.** (of Proposition 2.21) We sketch the proof using string diagrams. Recall that we must show

$$\mathrm{Hom}(X, Z \otimes Y^*) \cong \mathrm{Hom}(X \otimes Y, Z)$$

In one direction, this map is given by

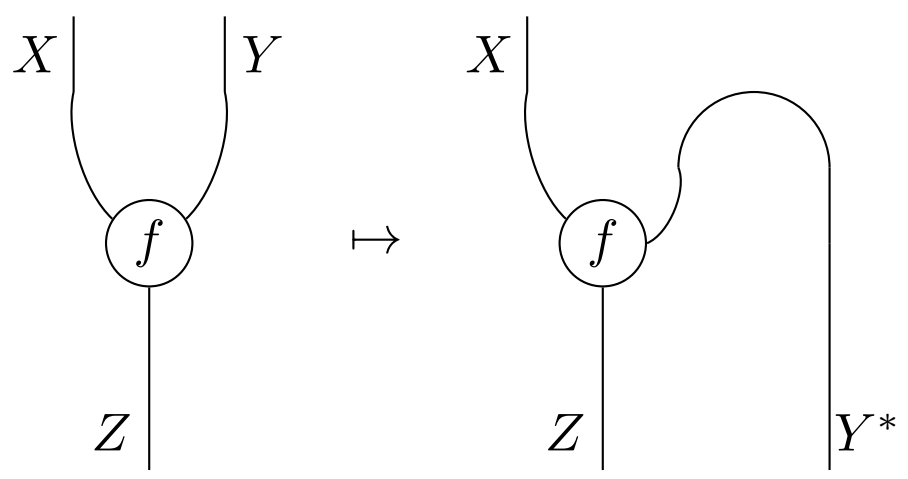

and in the other direction, this map is given by

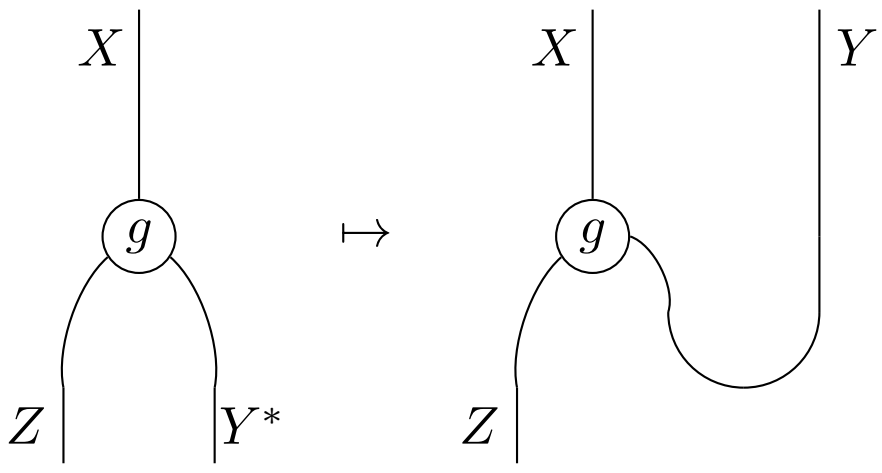

The "pulling straight" identity shows that these two operations are inverse, i.e.

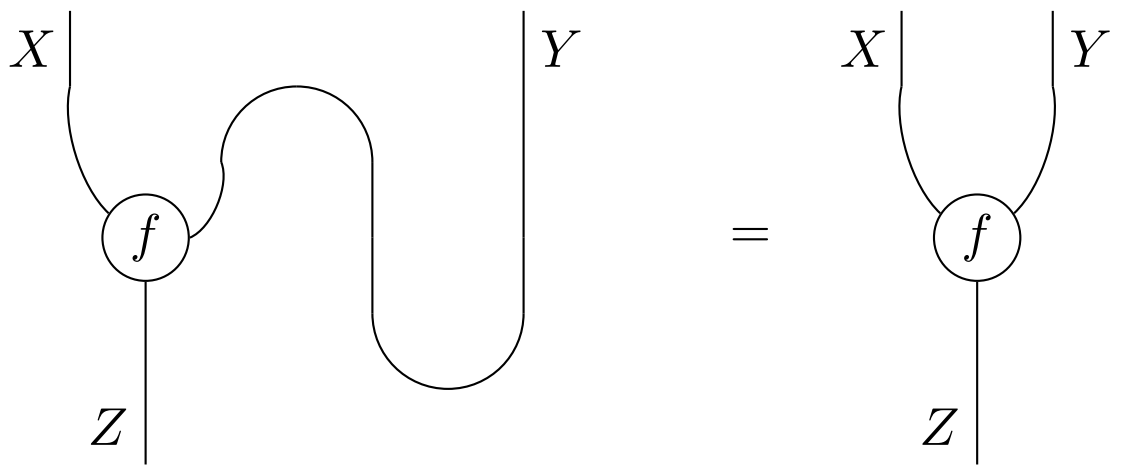

We are done. □

It should now be clear why we waited until we had developed string diagrams to give this proof; what would have been a great deal of tedious symbol manipulation can be reduced to a couple of intuitively clear pictures.

In general, in a compact closed category we have many options for moving wires around. Namely, we can cross wires using the symmetric monoidal structure, and we can bend wires using the closed structure. This idea is formalized by Joyal and Street [35], [36], where string diagrams are treated as topological spaces, and it is shown that if string diagram $A$ can be deformed into string diagram $B$, then they must have naturally isomorphic interpretations as morphisms in a monoidal category (possibly with symmetric or compact closed structure too).

We see many examples of string diagrams in the upcoming Chapter 3, but before we get there, we would be remiss if we did not discuss what the relevant morphisms between monoidal categories are.

## 2.4. Monoidal functors

We now discuss functors between monoidal categories. As is often the case in higher category theory, there are choices to make about how strict to be. We present three "gradations" of strictness, from most to least strict. Each of these variants have scenarios in which it is useful. These definitions can all be found in [30].

Definition 2.22. *A* ***strict monoidal functor*** *from a strict monoidal category* $(\mathsf{C}, \otimes_{\mathsf{C}}, I_{\mathsf{C}})$ *to a strict monoidal category* $(\mathsf{D}, \otimes_{\mathsf{D}}, I_{\mathsf{D}})$ *consists of a functor* $F\colon \mathsf{C} \to \mathsf{D}$ *such that*

- $F(X \otimes_{\mathsf{C}} Y) = F(X) \otimes_{\mathsf{D}} F(Y)$
- $F(I_{\mathsf{C}}) = I_{\mathsf{D}}$

**Example 2.23.** There is a "black projection" from $(\mathbb{N}^2, +, (0,0))$ to $(\mathbb{N}, +, 0)$, that sends $(n_1, n_2)$ to $n_1$. This sends a colored diagram like Figure 2.2 to just its black part. There is also a corresponding "black inclusion" from $\mathbb{N}$ to $\mathbb{N}^2$ that sends $n$ to $(n, 0)$.

As always, expecting these equalities to hold is in general not going to be true. So we use a weaker definition.

DEFINITION 2.24. *A **strong monoidal functor** between monoidal categories* $(\mathsf{C}, \otimes_{\mathsf{C}}, I_{\mathsf{C}})$ *and* $(\mathsf{D}, \otimes_{\mathsf{D}}, I_{\mathsf{D}})$ *consists of*

- *a functor* $F\colon \mathsf{C} \to \mathsf{D}$
- *an isomorphism* $\epsilon\colon I_{\mathsf{D}} \to F(I_{\mathsf{C}})$
- *a natural isomorphism* $\mu_{X,Y}\colon F(X) \otimes_{\mathsf{D}} F(Y) \xrightarrow{\sim} F(X \otimes_{\mathsf{C}} Y)$

*such that* $\mu$ *and* $\epsilon$ *interact with the unitor and associator in the appropriate way; namely the following diagrams commute. The diagram for the associator is*

$$\begin{array}{ccc}
(F(X) \otimes_{\mathsf{D}} F(Y)) \otimes_{\mathsf{D}} F(Z) & \xrightarrow{a_{F(X),F(Y),F(Z)}} & F(X) \otimes_{\mathsf{D}} (F(Y) \otimes_{\mathsf{D}} F(Z)) \\
{\scriptstyle \mu_{X,Y} \otimes_{\mathsf{D}} 1_{F(z)}} \downarrow & & \downarrow {\scriptstyle 1_{F(X)} \otimes_{\mathsf{D}} \mu_{Y,Z}} \\
F(X \otimes_{\mathsf{C}} Y) \otimes_{\mathsf{D}} F(Z) & & F(X) \otimes_{\mathsf{D}} F(Y \otimes_{\mathsf{C}} Z) \\
{\scriptstyle \mu_{X \otimes_{\mathsf{C}} Y, Z}} \downarrow & & \downarrow {\scriptstyle \mu_{X, Y \otimes_{\mathsf{C}} Z}} \\
F((X \otimes_{\mathsf{C}} Y) \otimes_{\mathsf{C}} Z) & \xrightarrow[F(a_{X,Y,Z})]{} & F(X \otimes_{\mathsf{C}} (Y \otimes_{\mathsf{C}} Z))
\end{array}$$

*and the diagrams for the left and right unitors are*

$$\begin{array}{ccc}
I_{\mathsf{D}} \otimes_{\mathsf{D}} F(X) & \xrightarrow{\epsilon \otimes_{\mathsf{D}} 1_{F(X)}} & F(I_{\mathsf{C}}) \otimes_{\mathsf{D}} F(X) \\
{\scriptstyle l_{F(X)}} \downarrow & & \downarrow {\scriptstyle \mu_{I_{\mathsf{C}}, X}} \\
F(X) & \xleftarrow[F(l_X)]{} & F(I_{\mathsf{C}} \otimes_{\mathsf{C}} X)
\end{array}$$

*and*

$$\begin{array}{ccc}
F(X) \otimes_{\mathsf{D}} I_{\mathsf{D}} & \xrightarrow{1_F(X) \otimes_{\mathsf{D}} \epsilon} & F(X) \otimes_{\mathsf{D}} F(I_C) \\
{\scriptstyle r_{F(X)}} \downarrow & & \downarrow {\scriptstyle \mu_{X, I_{\mathsf{C}}}} \\
F(X) & \xleftarrow[F(r_X)]{} & F(X \otimes_{\mathsf{C}} I_{\mathsf{C}})
\end{array}$$

These coherence diagrams suffice to show that there is a unique natural isomorphism $F(X_1) \otimes_{\mathsf{D}} \cdots \otimes_{\mathsf{D}} F(X_n) \to F(X_1 \otimes_{\mathsf{C}} \cdots \otimes_{\mathsf{C}} X_n)$ constructed out of $\mu$, which is reasonable.

The *first* step in categorification is to replace equalities with isomorphisms, but often a further step is taken, where the isomorphisms are replaced with morphisms that are not necessarily invertible. And so we have a third definition.

DEFINITION 2.25. *A **lax monoidal functor** between monoidal categories* $(\mathsf{C}, \otimes_{\mathsf{C}}, I_{\mathsf{C}})$ *and* $(\mathsf{D}, \otimes_{\mathsf{D}}, I_{\mathsf{D}})$ *consists of*

- *a functor* $F\colon \mathsf{C} \to \mathsf{D}$
- *a morphism* $\epsilon\colon I_{\mathsf{D}} \to F(I_{\mathsf{C}})$
- *a natural transformation* $\mu_{X,Y}\colon F(X) \otimes_{\mathsf{D}} F(Y) \xrightarrow{\sim} F(X \otimes_{\mathsf{C}} Y)$

*satisfying the commutative diagrams of Definition 2.24. Thus, we have a unique morphism* $F(X_1) \otimes_{\mathsf{D}} \cdots \otimes_{\mathsf{D}} F(X_n) \to F(X_1 \otimes_{\mathsf{C}} \cdots \otimes_{\mathsf{C}} X_n)$ *for every* $X_1, \ldots, X_n$.

**Example 2.26.** Consider the category 1 with just one object, say $*$, and its identity morphism. This category has a unique monoidal structure, which happens to be the cartesian monoidal structure. A lax monoidal functor $(F, \mu, \varepsilon)$ from 1 to $(\mathsf{Set}, \times, 1)$ consists of a set $F(*)$, along with a map

$$\mu\colon F(*) \times F(*) \to F(* \otimes *) = F(*)$$

and a map

$$\varepsilon\colon 1 \to F(*)$$

such that $(F(*), \mu, \varepsilon)$ is a monoid. Associativity of $\mu$:

$$\mu(a, \mu(b, c)) = \mu(\mu(a, b), c)$$

and unitality of $\mu$ and $\varepsilon$:

$$\mu(a, \varepsilon) = a = \mu(\varepsilon, a)$$

is implied by the commutative diagrams in Definition 2.24. In fact, all monoids can be represented as lax monoidal functors from 1 to $(\mathsf{Set}, \times, 1)$.

More generally, a lax monoidal functor from 1 into any monoidal category $(\mathsf{C}, \otimes, I)$ is a monoid-like structure internal to $\mathsf{C}$, and we call this a **monoid in the monoidal category $\mathsf{C}$**; see Maclane [27, Section III.6].

Lax monoidal functors are also important for defining operad algebras in Chapter 7.

Finally, in the case that we are working with symmetric monoidal categories, we additionally ask that functors between them respect the braiding. We only give the definition in the lax case: the strict and strong cases can be produced by requiring that the maps $\mu$ and $\epsilon$ are equalities, and isomorphisms, respectively.

Definition 2.27. *A* ***lax symmetric monoidal functor*** *between symmetric monoidal categories* $(\mathsf{C}, \otimes_{\mathsf{C}}, I_{\mathsf{C}})$ *and* $(\mathsf{D}, \otimes_{\mathsf{D}}, I_{\mathsf{D}})$ *is a lax monoidal functor* $(F, \mu, \epsilon)$ *such that the following diagram commutes.*

$$\begin{array}{ccc} F(X) \otimes_{\mathsf{D}} F(Y) & \xrightarrow{B_{F(X),F(Y)}} & F(Y) \otimes_{\mathsf{D}} F(X) \\ \downarrow{\scriptstyle \mu_{X,Y}} & & \downarrow{\scriptstyle \mu_{Y,X}} \\ F(X \otimes_{\mathsf{C}} Y) & \xrightarrow[F(B_{X,Y})]{} & F(Y \otimes_{\mathsf{C}} X) \end{array}$$

## 2.5. Categories of relations

We end by discussing a particularly useful construction of symmetric monoidal categories: that of symmetric monoidal categories of *relations*. This is because, as discussed before, the Willems approach to systems is essentially *relational*. That is, systems are not input-output machines, but rather should be thought of as imposing relations on the things that they are connected to, without a preordained sense of causality.

In order to get a handle on how this is formalized within category theory, we begin with a simple category and then we discuss how to generalize it.

Definition 2.28. *Let* $\mathsf{Rel}$ *be the following category. The objects of* $\mathsf{Rel}$ *are sets, and a morphism from* $X$ *to* $Y$ *is a subset*

$$R \subset X \times Y$$

*which we write as* $R: X \nrightarrow Y$ *and call a* ***relation from X to Y*** *(the slashed arrow represents relations as opposed to functions). Composition of a relation* $R: X \nrightarrow Y$ *and a relation* $S: Y \nrightarrow Z$ *is given by*

$$S \circ R = \{(x, z) \in X \times Z \mid \exists y \in Y, (x, y) \in R \wedge (y, z) \in S\}$$

We want to think of $\mathsf{Rel}$ as somehow "built out of" $\mathsf{Set}$. In this thesis, we encounter several categories that are similar to $\mathsf{Rel}$, but "built out of" categories other than $\mathsf{Set}$. Thus, it is worth thinking about what kind of category has the right properties in order to construct something like $\mathsf{Rel}$.

This right kind of category is a *regular category*. Regular categories have a long history, beginning with Freyd and Schedrov [37] and independently with Carboni and Walters [38] (though note that Carboni and Walters use the terminology "cartesian bicategory"). For a comprehensive textbook level treatment of the subject, we refer the reader to Borceux [39], or Johnstone [40]. Finally, a more modern take on the material can be found in Fong and Spivak [41], which emphasizes the graphical nature of regular categories.

Regular categories can be somewhat abstruse at first glance, so we attempt to explain them by showing how they are the natural structure in order to have a category of relations. A morphism $f: A \to B$ is **monic** (or alternatively, a **monomorphism**) if for any $g_1, g_2: X \rightrightarrows A$, $f \circ g_1 = f \circ g_2$ implies that $g_1 = g_2$. This generalizes injectivity. Conversely, $f: A \to B$ is **epic** (or alternatively, an **epimorphism**) if for any $g_1, g_2: B \rightrightarrows X$, $g_1 \circ f = g_2 \circ f$ implies that $g_1 = g_2$. This generalizes surjectivity.

We use monomorphisms to generalize the idea of subset.

DEFINITION 2.29. *For a category* $\mathsf{C}$*, and for an object* $X \in \mathsf{C}$*, we define* $\mathrm{Mono}(X)$ *to be the subcategory of* $\mathsf{C}/X$ *consisting of monomorphisms* $f: A \to X$*.*

PROPOSITION 2.30. *In* $\mathrm{Mono}(X)$*, two monomorphisms* $f: A \to X$ *and* $g: B \to X$ *are isomorphic in at most one way.*

**Proof.** Suppose that $h, h': A \to B$ are both isomorphisms such that the following diagrams commute:

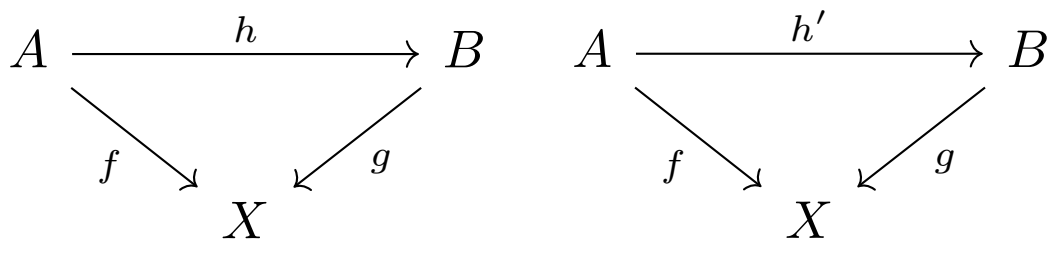

Then $g \circ h = f = g \circ h'$. Thus, by the fact that $g$ is a monomorphism, $h' = h$, and we are done. □

This implies that $\mathrm{Mono}(X)$ is equivalent to a category where isomorphic objects are equal.

DEFINITION 2.31. *Let* $\mathrm{Sub}(X)$ *be the category where the objects are equivalence classes of objects in* $\mathrm{Mono}(X)$ *under the equivalence relation of isomorphism. An object of* $\mathrm{Sub}(X)$ *is known as a* ***subobject*** *of* $X$*.*

**Example 2.32.** In $\mathsf{Set}$, we can identify the equivalence class of monomorphisms into a set $A$ with a literal subset of $A$.

DEFINITION 2.33. *If* $\mathsf{C}$ *is a category with products, a* ***relation*** *in* $\mathsf{C}$ *between objects* $X$ *and* $Y$ *is a subobject* $R$ *of* $X \times Y$*.*

We can also think of a relation as a **span**

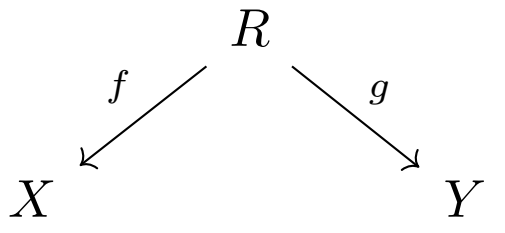

such that $f$ and $g$ are *jointly monic*, that is $\langle f, g\rangle: R \to X \times Y$ is monic. Before we investigate relationships, we first investigate spans.

DEFINITION 2.34. *If* $\mathsf{C}$ *is a category with pullbacks, the category* $\mathrm{Span}(\mathsf{C})$ *is defined in the following way. The objects are objects of* $\mathsf{C}$*, and a morphism from* $X$ *to* $Y$ *is an* equivalence class of spans

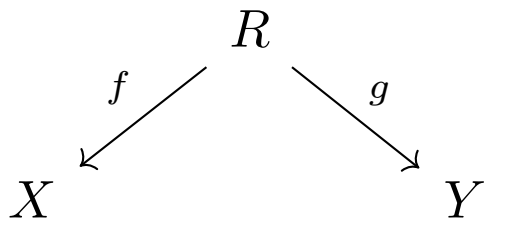

*where two spans* $(f: R \to X, g: R \to Y)$ *and* $(f': R' \to X, g': R' \to Y)$ *are equivalent if there is an isomorphism* $h: R \to R'$ *making the following diagram commute*

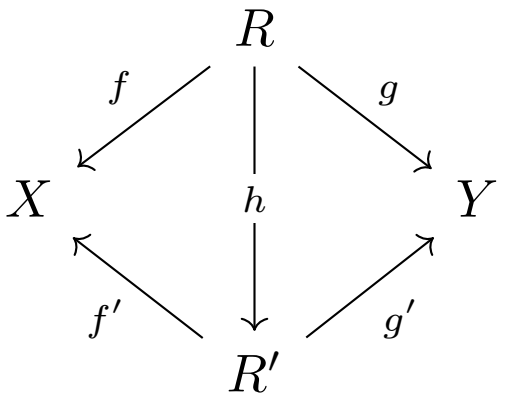

*Composition of spans is done via pullback, as in the following diagram.*

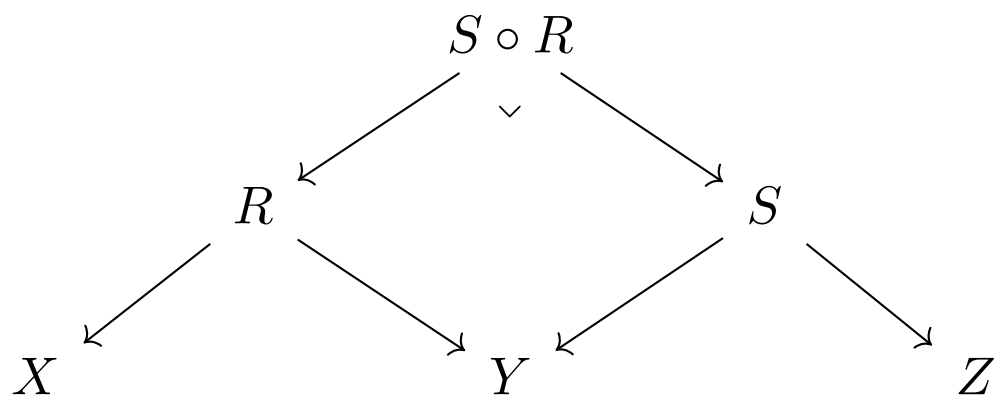

PROPOSITION 2.35. *Definition 2.34 specifies a well-defined category.*

**Proof.** For the original proof, see Benabou [42], and for a more modern take, see Baez [43]. However, we go over some of the details for showing that the composition is well-defined, for the sake of completeness.

Well-definedness comes in two parts. First of all, limits in a category are not defined uniquely; they are only unique up to unique isomorphism. This is not a problem however, because the equivalence class of the limit is uniquely specified.

Secondly, we must take into account the fact that the spans that we are composing are only specified by their equivalence class; it must be true that composing different representatives of the same equivalence classes produces the same result. This is true because the limits of isomorphic diagrams are isomorphic. We are done. □

However, we cannot simply compose relations via pullback, because even if $R \hookrightarrow X \times Y$ and $S \hookrightarrow Y \times Z$ are monic, $S \circ R \to X \times Z$ (as defined above) might not be. In Set we can fix this by taking the *image* of the map $S \circ R \to X \times Z$. In order to do this in a general category, we need what are known as *image factorizations*.

Definition 2.36. *Let $f: A \to B$ be a morphism in a category* C*. We say that*

$$A \xrightarrow{e} C \xrightarrow{m} B$$

*is an* ***image factorization*** *of $f$ if $e$ is epic, $m$ is mono, and given any other factorization $f = m' \circ e'$ into an epimorphism $e$ and a monomorphism $m$, $m = m' \circ k$ for some $k$.*

Proposition 2.37. *Image factorizations are unique up to canonical isomorphism. That is if $A \xrightarrow{e} C \xrightarrow{m} B$ and $A \xrightarrow{e'} C' \xrightarrow{m'} B$ are both image factorizations of $f: A \to B$, then there exists a unique isomorphism $\varphi: C \to C'$ such that the following diagram commutes:*

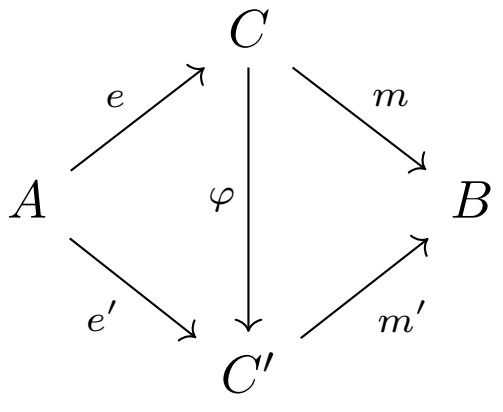

Corollary 2.38. *If $f$ has an image factorization $A \xrightarrow{e} C \xrightarrow{m} B$, then there is a unique subobject of $B$ corresponding to $C \xrightarrow{m} B$. We call this subobject the* ***image*** *of $f$, and write it as $\mathrm{im}(f)$.*

Note that $\mathrm{im}(f)$ is not an object of C, it is a subobject of $B$, which is an equivalence class of monomorphisms into $B$. That being said, we will immediately abuse notation and treat $\mathrm{im}(f)$ as an object of C.

**Example 2.39.** In Set, every function $f: A \to B$ has a image factorization with

$$\mathrm{im}(f) = \{b \in B \mid \exists a \in A, f(a) = b\}$$

which coincides with the traditional notion of image. We can see in this example that the idea of "image" is very much connected to the idea of "existentials"; recall that we needed an existential quantifier to define composition of relations in Definition 2.28.

We can also define images in Set in a different way. If $f: A \to B$ is any function, let

$$\mathrm{im}(f) = A / \sim \tag{2.1}$$

where the equivalence relation $\sim$ is given by $a_1 \sim a_2$ if and only if $f(a_1) = f(a_2)$.

This can be replicated in an arbitrary category in the following way. Suppose that $f: A \to B$ is any morphism. Then take the pullback

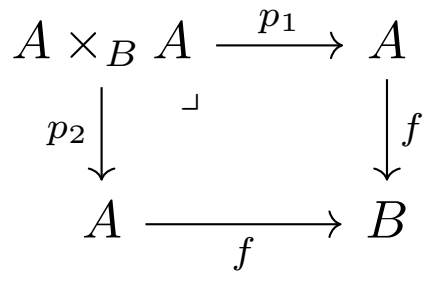

The pair of morphisms $p_1, p_2\colon A \times_B A \rightrightarrows B$ is called the **kernel pair** of $f$. If we now take the coequalizer

$$A \times_B A \underset{p_2}{\overset{p_1}{\rightrightarrows}} A \longrightarrow \mathrm{coeq}(p_1, p_2)$$

then the universal property of coequalizers gives us that $f$ factors through $\mathrm{coeq}(p_1, p_2)$. Moreover, the map from $A$ to $\mathrm{coeq}(p_1, p_2)$ is an epimorphism because it is a coequalizer map. Later, we show that the map from $\mathrm{coeq}(p_1, p_2)$ to $B$ is a monomorphism, so this is indeed an image factorization.

If one thinks carefully about this construction, one sees that it is precisely generalizes the construction of $\mathrm{im}(f) = A/\sim$. This generalization serves us to define the right structure to use to build a category of relations.

Definition 2.40. *A **regular epimorphism** in a category* $\mathsf{C}$ *is a morphism* $f\colon A \to B$ *that is the coequalizer* $X \rightrightarrows A \overset{f}{\longrightarrow} B$ *of some pair of morphisms.*

**Example 2.41.** Any epimorphism in $\mathsf{Set}$ is a regular epimorphism.

**Example 2.42.** In $\mathsf{CMon}$, the category of commutative monoids, the inclusion $f\colon \mathbb{N} \to \mathbb{Z}$ is an epimorphism, as a monoid homomorphism $\varphi\colon \mathbb{Z} \to M$ is totally determined by $\varphi(1)$. However, $f$ is not a regular epimorphism.

A *regular category* has precisely the properties needed to form images in the manner of Example 2.39, and have them be well-behaved enough to form a good category of relations.

Definition 2.43. *A **regular category** is a category* $\mathsf{C}$ *such that*

- $\mathsf{C}$ *has all finite limits (so in particular, products and pullbacks)*
- *every kernel pair has a coequalizer*
- *the pullback of a regular epimorphism along any morphism is again a regular epimorphism*

This last condition is a bit dense; what it means is that if $f\colon A \to C$ is a regular epimorphism, and $g\colon B \to C$ is any morphism, then $g^*(f)$ in the following diagram is also a regular epimorphism.

$$\begin{array}{ccc} A \times_C B & \longrightarrow & A \\ {\scriptstyle g^*(f)}\downarrow & \lrcorner & \downarrow{\scriptstyle f} \\ B & \underset{g}{\longrightarrow} & C \end{array}$$

**Example 2.44.** $\mathsf{Set}$ is a regular category.

**Example 2.45.** If $\mathsf{D}$ is a regular category, and $\mathsf{C}$ is any category, $\mathsf{D}^{\mathsf{C}}$ is a regular category.

**Example 2.46.** The category of models of any "finitary algebraic theory" is regular. Precisely, this means the category models of any Lawvere theory is regular; see Barr and Wells [44, Chapter 4] for more details on Lawvere theories. For the reader unfamiliar with Lawvere theories, this class of categories includes the categories of

- (commutative) monoids
- (abelian) groups
- rings (and commutative rings, rings with unit, etc.)
- modules over a ring
- vector spaces over a field
- convex spaces (as defined in Chapter 9)

**Example 2.47.** A "convenient" category of topological spaces, e.g. the category $\mathsf{cgHaus}$ of compactly-generated Hausdorff spaces, is regular. This is one of the reasons that algebraic topologists prefer working with categories other than the category of traditional topological spaces. The use of such convenient categories goes back to Steenrod [45].

Proposition 2.48. *If* $\mathsf{C}$ *is a regular category, and* $f: A \to B$ *is any morphism, then* $f$ *factors into*

$$A \xrightarrow{e} \mathrm{im}(f) \xrightarrow{m} B$$

*where* $e$ *is a regular epimorphism and* $m$ *is a monomorphism, and moreover this factorization is unique up to unique isomorphism.*

**Proof.** We follow the construction of Example 2.39, i.e. we take the coequalizer of the kernel pair. The full proof can be found in Borceux [39, Theorem 2.1.3]. □

We now finally give the general construction for a category of relations.

Construction 2.49. *Suppose that* $\mathsf{C}$ *is a regular category. Then define the category* $\mathrm{Rel}(\mathsf{C})$ *in the following way. The objects are the objects of* $\mathsf{C}$*, and a morphism from* $A$ *to* $B$ *is a subobject* $R \hookrightarrow A \times B$.

*To define composition of morphisms, suppose that* $R \hookrightarrow A \times B$ *and* $S \hookrightarrow B \times C$ *are relations. Then let* $f: T \to A \times C$ *be constructed by pullback in the following diagram*

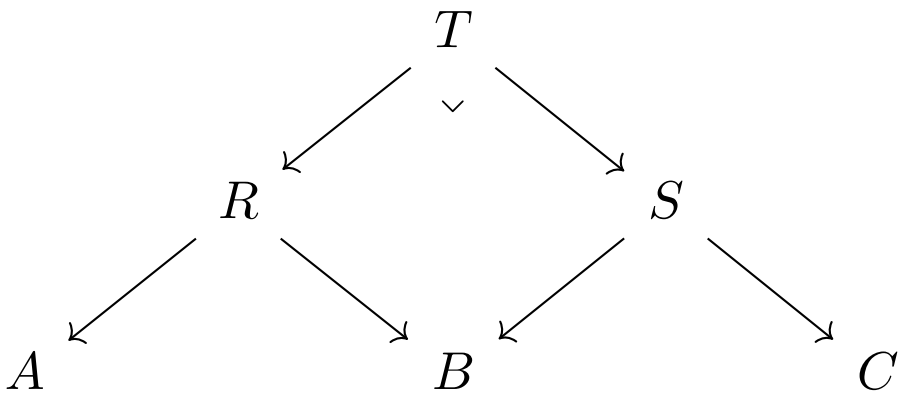

*We then let* $S \circ R = \mathrm{im}(f)$ *with its natural monomorphism into* $A \times C$.

*This is all well-defined because both image factorizations and pullbacks are defined up to isomorphism so the same reasoning as in the definition of* $\mathrm{Span}(\mathsf{C})$ *applies. We call* $\mathrm{Rel}(\mathsf{C})$ *the* ***category of relations*** *for* $\mathsf{C}$.

The category of relations of a regular category has a natural compact closed structure. We exposit this with a series of propositions, given without proof.

Proposition 2.50. *If* $\mathsf{C}$ *is a regular category, then there is a functor* $\mathrm{graph}: \mathsf{C} \to \mathrm{Rel}(\mathsf{C})$ *that is the identity on objects, and sends a morphism* $f: A \to B$ *to the relation*

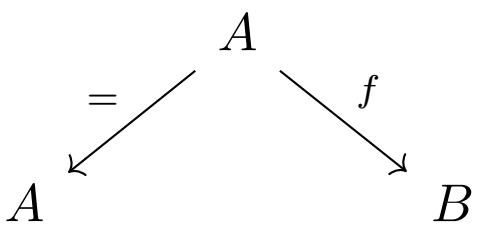

*There is also a functor* $\mathrm{cograph}: \mathsf{C}^{\mathrm{op}} \to \mathrm{Rel}(\mathsf{C})$ *defined in the obvious complementary way. Moreover, both of these functors are injective on morphisms.*

The next Proposition can be found in Borceux [39, Theorem 2.8.4].

Proposition 2.51. *If* $\mathsf{C}$ *is a regular category, there is a natural symmetric monoidal structure on* $\mathrm{Rel}(\mathsf{C})$ *defined in the following way. On objects of* $\mathrm{Rel}(\mathsf{C})$ *(which are objects of* $\mathsf{C}$*), we define* $A \otimes B$ *to be the categorical product in* $\mathsf{C}$ *(note that this is not the categorical product in* $\mathrm{Rel}(\mathsf{C})$*; this is why we use* $\otimes$ *instead of* $\times$*). Then if* $R_1: A_1 \nrightarrow B_1$ *and* $R_2: A_2 \nrightarrow B_2$*, we define* $R_1 \otimes R_2: A_1 \otimes A_2 \nrightarrow B_1 \otimes B_2$ *to be the subobject of* $(A_1 \otimes A_2) \times (B_1 \otimes B_2)$ *given by applying the isomorphism*

$$(A_1 \times B_1) \times (A_2 \times B_2) \cong (A_1 \times A_2) \times (B_1 \times B_2) = (A_1 \otimes A_2) \times (B_1 \otimes B_2)$$

*to the categorical product*

$$R_1 \times R_2 \hookrightarrow (A_1 \times B_1) \times (A_2 \times B_2)$$

*The monoidal unit* $I$ *is the terminal object in* $\mathsf{C}$.

Finally, note that subobjects of $(A \times B) \times C$ are in natural bijection with subobjects $A \times (C \times B)$, so $\mathrm{Hom}_{\mathrm{Rel}(\mathsf{C})}(A \otimes B, C) \cong \mathrm{Hom}_{\mathrm{Rel}(\mathsf{C})}(A, C \otimes B)$. Thus, we expect that $B$ is its own dual. We record this well-known fact in the following proposition.

Proposition 2.52. *If* $\mathsf{C}$ *is a regular category, then every object in* $\mathrm{Rel}(\mathsf{C})$ *is its own dual, with maps* $1 \to A \otimes A$ *and* $A \otimes A \to 1$ *given by the diagonal relation* $A \hookrightarrow A \times A$*. Thus,* $\mathrm{Rel}(\mathsf{C})$ *is a compact closed category.*

We would give examples at this point, but in fact the entire next chapter consists of examples of this structure, so we simply move on to that.

# 3. Categorical Linear Algebra

## 3.1. Linear maps

Linear algebra in monoidal category theory has a long history. String diagrams were first developed as a notation for working with tensors by Penrose [46], and are often known as "Penrose graphical notation" in this context. However, the use of a different monoidal structure for string diagrams, the direct product monoidal structure, is more recent and is independently due to [6] and [47].

For this chapter, we work over an arbitrary field $K$ of characteristic 0, typically $\mathbb{R}$, but perhaps also $\mathbb{C}$ or $\mathbb{R}(s)$ (the field of fractions of the ring of polynomials $\mathbb{R}[s]$). Moreover, all our vector spaces are finite-dimensional.

Definition 3.1. *The category* $\mathsf{Lin}$ *has vector spaces as objects, and linear transformations as morphisms.*

Proposition 3.2. $\mathsf{Lin}$ *has a symmetric monoidal structure given by the direct product* $\oplus$

$$V \oplus W = \{(v, w) \mid v \in V, w \in W\}$$

*and unit vector space* $0 = K^0$.

The direct product $\oplus$ is both the categorical product, with projection maps

- $\pi_1: V \oplus W \to V$, $\pi_1(v, w) = v$
- $\pi_2: V \oplus W \to W$, $\pi_2(v, w) = w$

and the categorical coproduct, with injection maps

- $\iota_1: V \to V \oplus W$, $\iota_1(v) = (v, 0)$
- $\iota_2: W \to V \oplus W$, $\iota_2(w) = (0, w)$

Thus, $\oplus$ is a **biproduct**, which is a word for an operation that is both product and coproduct; moreover, 0 is both a terminal and initial object in $\mathsf{Lin}$. Consequently, $(\mathsf{Lin}, \oplus, 0)$ is both a cartesian monoidal category and a cocartesian monoidal category.

**Example 3.3.** For every object $V \in \mathsf{Lin}$ and every $n \in \mathbb{N}$, there is a map $\bullet: V^{\oplus n} \to V$ given by

$$(v_1, \ldots, v_n) \mapsto v_1 + \cdots + v_n$$

Note that this is given by the universal property of the coproduct $V^{\oplus n}$, applied to the maps $\mathrm{id}_V, \ldots, \mathrm{id}_V$. Thus, this map would exist in any cocartesian monoidal category; the reader should think what this map would be for $\mathsf{Set}$. Moreover, when $n = 0$, this is the unique linear map $0 \to V$. We draw this map using a black dot with $n$ inputs and one output:

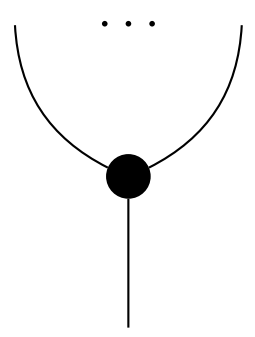

This has the curious property that a map composed only of instances of $\bullet$ is determined solely by its connection structure as a string diagram. That is, suppose that $f: V^{\oplus m} \to V^{\oplus n}$ is comprised solely of instances of $\bullet$. Then

$$f(x)_j = \sum_{i \text{ connected to } j} x_i$$

This implies that, for instance

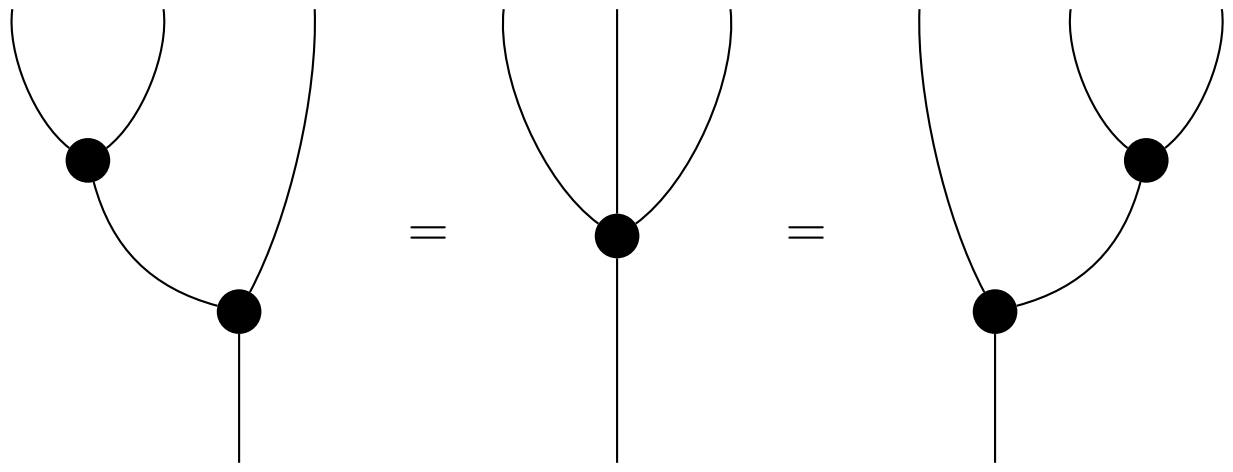

which can be seen as a type of "associativity law".

If we just consider the nullary and binary versions of $\bullet$, they form a monoid in $(\mathsf{Lin}, \oplus, 0)$ (see Example 2.26 for more about monoids in monoidal categories).

**Example 3.4.** For every object $V \in \mathsf{Lin}$ and every $n \in \mathbb{N}$, there is a map $\circ: V \to V^{\oplus n}$ given by

$$v \mapsto (v, \ldots, v)$$

Note that this is given by the universal property of the product $V^{\oplus n}$, applied to the maps $\mathrm{id}_V, \ldots, \mathrm{id}_V$. Thus, this map would exist in any cartesian monoidal category. Moreover, when $n = 0$ this is the unique linear map $V \to 0$. We draw this map using an unfilled dot with one input and $n$ outputs.

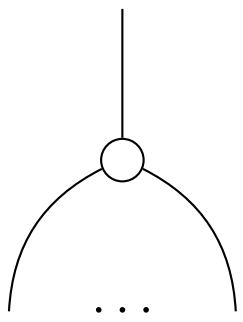

This has a similar property to $\bullet$. That is, in a string diagram comprised only of copies of $\circ$, all that matters is the connection structure. Structures like $\circ$ are known by the name of **comonoid in a monoidal category** [47].

**Example 3.5.** For every $V$, and every $\lambda \in K$, where $K$ is the underlying field, there is a map $m_\lambda$: $V \to V$ given by $v \mapsto \lambda v$. We draw this as a triangle with $\lambda$ in it, like the following:

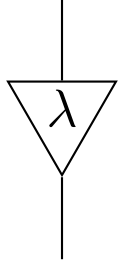

Proposition 3.6. *Any linear map $K^{\oplus m} \to K^{\oplus n}$ can be represented as a composite of $\bullet$, $\circ$, and $m_\lambda$, using both normal and monoidal composition.*

**Proof.** Any linear map can be represented as a matrix. We can then express a matrix by successive steps of $\circ$, $m_\lambda$, and $\bullet$, as illustrated in Figure 3.1. □

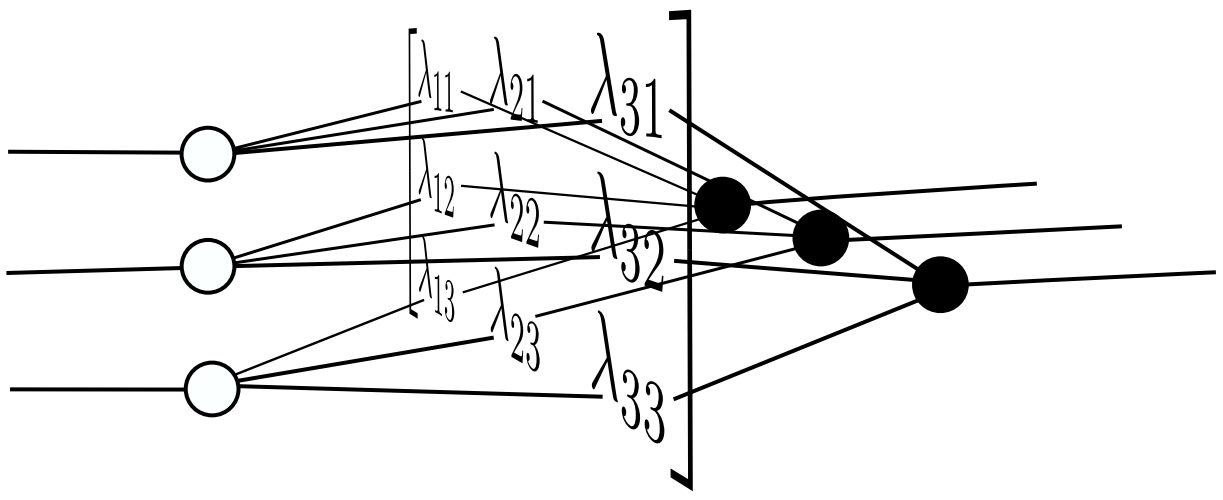

**Figure 3.1.** Proof by picture that all morphisms in $\mathsf{Lin}$ are be generated by $\bullet$, $\circ$, and $m_\lambda$

## 3.2. Linear relations

We can apply the machinery that we developed in section 2.5 to $\mathsf{Lin}$ to develop a *category of linear relations*. This works because $\mathsf{Lin}$ is a regular category (and this is because $\mathsf{Lin}$ is the category of models of a Lawvere theory; see Lawvere [48] or Barr and Wells [44]).

DEFINITION 3.7. *A **linear relation** from a vector space $V$ to a vector space $W$ is a linear subspace $R \subset V \oplus W$. We write this as $R: V \nrightarrow W$. If $R: U \nrightarrow V$ and $S: V \nrightarrow W$, then*

$$S \circ R = \{(u, w) \mid \exists v \in V.(u, v) \in R, (v, w) \in S\} \subset U \oplus W$$

*Then let* LinRel *be the category with vector spaces as objects and linear relations as morphisms.*

As discussed in Proposition 2.50, we can take a linear map $L: V \to W$ and turn it into a linear relation

$$\mathrm{graph}(L) = \{(v, Lv) \mid v \in V\} \subset V \oplus W$$

We can do the same for a map $L: W \to V$:

$$\mathrm{cograph}(L) = \{(Lw, w) \mid w \in W\} \subset V \oplus W$$

Moreover, we have for $M: W \to X$, $L: V \to W$,

$$\mathrm{graph}(M) \circ \mathrm{graph}(L) = \mathrm{graph}(M \circ L)$$

and

$$\mathrm{cograph}(L) \circ \mathrm{cograph}(M) = \mathrm{cograph}(M \circ L)$$

Thus, both Lin and $\mathsf{Lin}^{\mathrm{op}}$ are naturally subcategories of LinRel.

The direct product $\oplus$ is also a symmetric monoidal product on LinRel, even though it is neither the categorical product nor coproduct in LinRel. All of the necessary natural isomorphisms come from Lin when thought of as a subcategory of LinRel as discussed before. This is a consequence of Proposition 2.51. Finally, LinRel is a compact closed monoidal category, by Proposition 2.52.

The maps $\bullet$, $\circ$, and $\lambda$ from Lin all extend to more general relations in LinRel.

**Example 3.8.** For every vector space $V$ and $m, n \in \mathbb{N}$, there is a linear relation $\bullet: V^{\oplus m} \nrightarrow V^{\oplus n}$ defined by

$$((v_1, \ldots, v_m), (v'_1, \ldots, v'_n)) \in \bullet \text{ iff } \sum_{i=1}^{m} v_i = \sum_{j=1}^{n} v'_j$$

This is the composite of

$$\mathrm{graph}(\bullet: V^{\oplus m} \to V): V^{\oplus m} \nrightarrow V$$

and

$$\mathrm{cograph}(\bullet: V^{\oplus n} \to V): V \nrightarrow V^{\oplus n}$$

We draw this in the obvious way, as a black dot with $m$ wires coming in and $n$ wires coming out.

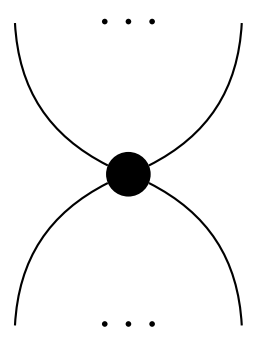

We call this the **summing junction**. This models a situation where there is some conserved quantity, like electrical current, flowing through a junction. The sum of the currents going in is equal to the sum of the currents coming out.

Just like with $\bullet$ in Lin, if we have a diagram that is only composed of $\bullet$, the only thing that matters is the connection structure of the diagram. These properties make the collection of morphisms collectively referred to by $\bullet$ into a so-called "Frobenius monoid" [6].

**Example 3.9.** We also make a Frobenius monoid (or "junction") out of $\circ$, which we call the **matching junction**.

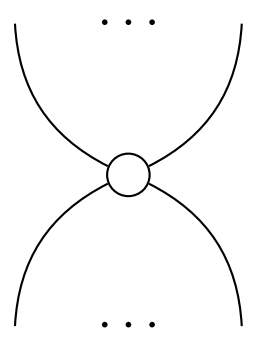

This is defined by

$$\circ = \{((v, \ldots v), (v, \ldots, v)) \mid v \in V\}: V^{\oplus m} \nrightarrow V^{\oplus n}$$

That is, all of the wires going into ∘ must have the same value in $V$. This models some quantity, like electrical potential, that is held to be equal for all of the wires.

Note that in the real world wires carry a voltage and a potential; keep this in mind for Chapter 4.

**Remark 3.10.** The story of "scalars" in $\mathsf{LinRel}$ is quite interesting, because the morphism $m_\lambda$: $V \to V$ does not simplify as easily in the relational setting as • and ∘. Thus, we use a downwards triangle for graph($m_\lambda$) and an upwards triangle for cograph($m_\lambda$). When $\lambda$ and $\gamma$ are non-zero, then all the following are equal, where $\bar{\gamma} = \gamma^{-1}$ and $\bar{\lambda} = \lambda^{-1}$ (for compactness).

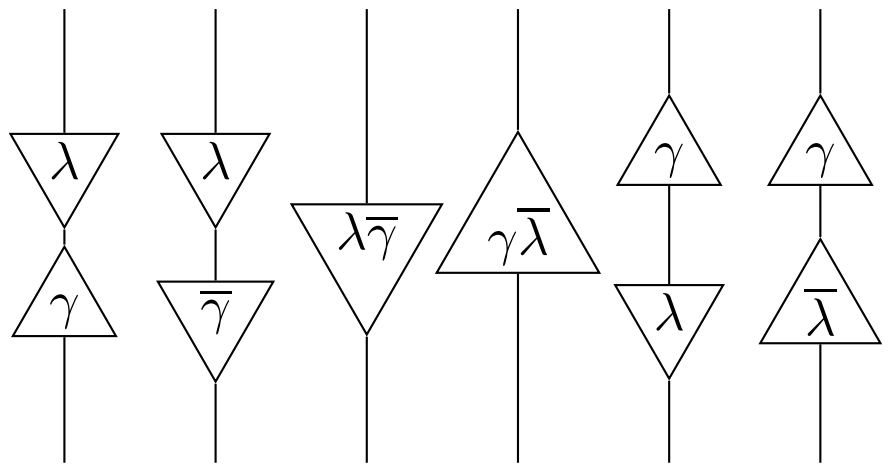

The subspace represented by all of these is $\{(v, v') | v, v' \in V, \lambda v = \gamma v'\}$. So the group of units (invertible elements) of $K$ injects into the monoid of linear relations $R: V \nrightarrow V$ with composition. However, the following four morphisms are all different.

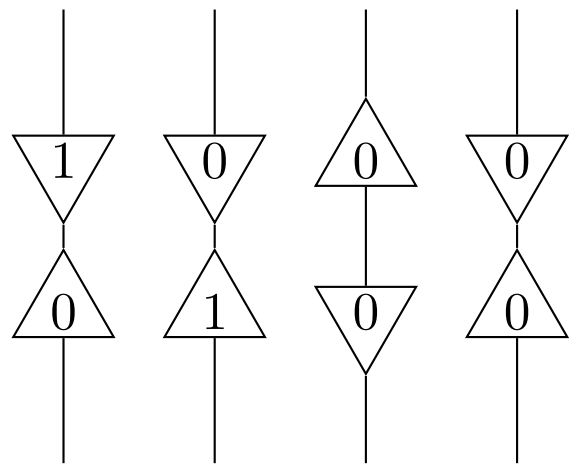

From left to right, their respective subspaces are $\{(0, v) | v \in V\}$, $\{(v, 0) \,|\, v \in V\}$, $\{(0, 0)\}$, $\{(v, v') \,|\, v, v' \in V\}$. This gives an interesting perspective on the question of "dividing by zero". Namely, we actually *can* give an non-trivial semantics to something similar to the inverse of 0 (i.e. cograph($m_0$)), as long as we are willing for our monoid to be non-commutative (as witnessed by the non-commutativity of graph($m_0$) and cograph($m_0$)); see Sobocinski [49].

## 3.3. Multilinear algebra

The category $\mathsf{Lin}$ has another monoidal product on it: the tensor product $\otimes$. This is the more traditional monoidal product on $\mathsf{Lin}$ (as opposed to $\oplus$), and it is the tensor product that lends monoidal categories their symbol for a generic monoidal product.

The tensor product is neither the categorical product nor the categorical coproduct, however it does have some interesting categorical structure. The reader is likely somewhat familiar with the tensor product, however we hope that this section can deepen the reader's knowledge. For an elementary reference on multilinear algebra, the reader should consult Treil [50, Chapter 8], and for an in-depth exposition of string diagrams and monoidal categories in multilinear algebra, the reader should consult Coecke [51].

A slick way in category theory to define something is by saying that it is a *representation* of a certain functor.

Definition 3.11. *Let $\mathsf{C}$ be any category, and let $F: \mathsf{C} \to \mathsf{Set}$ be a functor. We say that $X \in \mathsf{C}$ is a* ***representation*** *for $F$ if $F$ is naturally isomorphic to* $\mathrm{Hom}(X, -)$.

Proposition 3.12. *If a representation for $F$ exists, it is unique up to canonical isomorphism.*

**Proof.** This is a corollary of the Yoneda lemma, but we will give an explicit proof for pedagogical purposes. To prove this, we show that the map $y: \mathsf{C}^{\mathrm{op}} \to \mathsf{Set}^{\mathsf{C}}$ defined by $X \mapsto \mathrm{Hom}(X, -)$ is a full and faithful functor. Having shown this, if $F \cong \mathrm{Hom}(X, -)$ and we have any isomorphism $F \cong \mathrm{Hom}(Y, -)$, we can compose these isomorphisms to get $\mathrm{Hom}(X, -) \cong \mathrm{Hom}(Y, -)$, and thus $X \cong Y$.

The map $y$ is often called the "Yoneda embedding"; this proof provides the reason why it is called an embedding.

It is not so hard to show that $y$ is faithful; if $f, g: Y \to X$ are distinct morphisms, then $f^*, g^*$: $\mathrm{Hom}(X, -) \to \mathrm{Hom}(Y, -)$ given by postcomposition are clearly distinct, as $\mathrm{id}_X \circ f \neq \mathrm{id}_X \circ g$.

To show that $y$ is full, suppose that $\alpha$ is a natural transformation between $\mathrm{Hom}(X, -)$ and $\mathrm{Hom}(Y, -)$. Then we claim that $\alpha = \alpha(\mathrm{id}_X)^*$. To make sense of this, first note that $\alpha(\mathrm{id}_X) \in \mathrm{Hom}(Y, X)$, and recall that $\alpha(\mathrm{id}_X)^*$ is the postcomposition functor. We now claim that for any $f: X \to Z$, $\alpha_Z(f) = \alpha_X(\mathrm{id}_X)^* f = \alpha_X(\mathrm{id}_X) \circ f$. This follows from the following naturality diagram:

$$\begin{array}{ccc} \mathrm{Hom}(X, X) & \xrightarrow{\alpha_X} & \mathrm{Hom}(Y, X) \\ {\scriptstyle f^*}\downarrow & & \downarrow{\scriptstyle f^*} \\ \mathrm{Hom}(X, Z) & \xrightarrow[\alpha_Z]{} & \mathrm{Hom}(Y, Z) \end{array}$$

If we start with $\mathrm{id}_X \in \mathrm{Hom}(X, X)$, following the top we get $\alpha_X(\mathrm{id}_X) \circ f$, and following the bottom we get $\alpha_Z(f)$. Thus $\alpha = \alpha(\mathrm{id}_X)^*$, and we are done. □

DEFINITION 3.13. *A* ***bilinear map*** $f: X \times Y \to Z$ *is a function that satisfies*

- $f(\lambda_1 x_1 + \lambda_2 x_2, y) = \lambda_1 f(x_1, y) + \lambda_2 f(x_2, y)$
- $f(x, \lambda_1 y_1 + \lambda_2 y_2) = \lambda_1 f(x, y_1) + \lambda_2 f(x, y_2)$

*For any* $X, Y \in \mathsf{Lin}$, *we can define a functor* $\mathrm{Bilin}_{X,Y}: \mathsf{Lin} \to \mathsf{Set}$ *that sends* $Z$ *to the set of bilinear functions* $X \times Y \to Z$, *and sends a map* $h: Z \to Z'$ *to the map*

$$\mathrm{Bilin}_{X,Y}(h)(f: X \times Y \to Z) = h \circ f: X \times Y \to Z'$$

*We often call a bilinear map* $X \times X \to K$ *a* ***bilinear form****.*

DEFINITION 3.14. *A* ***tensor product*** $X \otimes Y$ *is a representation of* $\mathrm{Bilin}_{X,Y}$. *That is, the set* $\mathsf{Lin}(X \otimes Y, Z)$ *is naturally isomorphic to* $\mathrm{Bilin}_{X,Y}(Z)$.

PROPOSITION 3.15. *There exists a tensor product for every* $X, Y \in \mathsf{Lin}$.

**Proof.** This is a standard construction, see [52, Chapter VIII.2]. □

This means that we can say "the" tensor product, in accordance with the principle of equivalence, as it exists and is unique up to unique isomorphism.

DEFINITION 3.16. *The identity* $X \otimes Y \to X \otimes Y$ *corresponds to a bilinear map* $X \times Y \to X \otimes Y$ *that we denote by* $(x, y) \mapsto x \otimes y$.

It can be shown that in general, $X_1 \otimes \cdots \otimes X_n$ is a representation for the functor $Z \mapsto \{n\text{-linear maps } X_1 \times \cdots \times X_n \to Z\}$, where $n$-linear is the natural generalization of bilinear. Using this fact, and the fact that $\times$ is symmetric, associative, and unital, we can show that $\otimes$ is a symmetric monoidal product with $I = K$ being the monoidal unit, so $(\mathsf{Lin}, \otimes, K)$ is a symmetric monoidal category. We record this well-known fact in the following Proposition, which can be found in [30].

PROPOSITION 3.17. $(\mathsf{Lin}, \otimes, K)$ *is a symmetric monoidal category.*

This symmetric monoidal category has additional structure on it, because the set of linear maps $\mathsf{Lin}(X, Y)$ has a natural vector-space structure. To emphasize this, we use the notation $[X, Y]$ to mean the *vector space* of linear maps from $X$ to $Y$. Moreover, $[-, -]$ is a functor from $\mathsf{Lin}^{\mathrm{op}} \times \mathsf{Lin}$ to $\mathsf{Lin}$; this is because we can define its action on morphisms in $\mathsf{Lin}^{\mathrm{op}} \times \mathsf{Lin}$ to be pre- and post-composition.

Now, bilinear functions $f: X \times Y \to Z$ are in one-to-one correspondence with linear functions $\hat{f}: X \to [Y, Z]$, via the definition

$$\hat{f}(x) = f(x, -)$$

This gives a natural isomorphism $\mathsf{Lin}(X \otimes Y, Z) \cong \mathsf{Lin}(X, [Y, Z])$. Thus, $- \otimes Y$ is a left adjoint to $[Y, -]$. This results in the following proposition.

Proposition 3.18. *$(\mathsf{Lin}, \otimes, I)$ is a closed monoidal category, where the right adjoint to $-\otimes Y$ is $[Y, -]$.*

The tensor product starts to really get interesting when we combine it with duals.

Definition 3.19. *For a vector space V, we define the **dual space** $V^*$ to be $[V, K]$. Moreover, we define the contravariant functor $(-)^* = [-, K]$, so that if $L: V \to W$, then $L^*: W^* \to V^*$.*

Definition 3.20. *There is a canonical bilinear map $V^* \times V \to K$ (or equivalently, a linear map $V^* \otimes V \to K$) that we denote by*

$$(\varphi, v) \mapsto \varphi(v) = \langle \varphi, v \rangle$$

*and call the **natural pairing**.*

Proposition 3.21. *If $e_1, \ldots, e_n$ is a basis for V, then $f_1, \ldots, f_n$ defined by*

$$f_i(e_j) = \delta_{ij}$$

*is a basis for $V^*$, called the **dual basis**. As a consequence, $\dim V \cong \dim V^*$.*

Proposition 3.22. *There is a natural map $V \to V^{**}$ defined by*

$$v \mapsto \langle -, v \rangle$$

*This map is injective, and thus an isomorphism because $\dim V = \dim V^*$ (recall that we work only with finite-dimensional vector spaces). This is known as the "double-dual isomorphism."*

Proposition 3.23. *$(V \otimes W)^* \cong V^* \otimes W^*$.*

**Proof.** For every $(\varphi, \eta) \in V^* \times W^*$, there is a bilinear map $V \times W \to K$ defined by

$$(v, w) \mapsto \varphi(v)\, \eta(w)$$

This construction gives a natural map $V^* \otimes W^* \to (V \otimes W)^*$, and because these vector spaces both have the same dimension ($\dim V \otimes W = \dim V \dim W$) and this map is injective, this must be an isomorphism. □

Finally, as one might guess, it turns out that linear algebra dual is the same as the categorical dual (see Definition 2.20) in the context of the symmetric monoidal category $(\mathsf{Lin}, \otimes, K)$. We state this more formally in a proposition.

Proposition 3.24. *The categorical dual of a vector space V in the symmetric monoidal category $(\mathsf{Lin}, \otimes, K)$ is $V^* = [V, K]$. Specifically, consider the maps $\mathrm{id}_V: K \to V \otimes V^*$ and $\mathrm{ev}_V: V^* \otimes V \to K$ defined in the following way. Let $v_1, \ldots, v_n$ be a basis of V with corresponding dual basis $\varphi_1, \ldots, \varphi_n$ of $V^*$, and define*

$$\mathrm{id}_V(\lambda) = \lambda(v_1 \otimes \varphi_1 + \cdots + v_n \otimes \varphi_n)$$

*(It can be shown that this is actually basis-independent). Also, define $\mathrm{ev}_V$ to be the map equivalent to the bilinear map*

$$(\varphi, v) \mapsto \langle \varphi, v \rangle$$

*Then the two maps $\mathrm{id}_V$ and $\mathrm{ev}_V$ satisfy the zig-zag identities, proving that $V^*$ is the categorical dual of V.*

Corollary 3.25. *$(\mathsf{Lin}, \otimes, K)$ is a compact closed category.*

In fact, $(\mathsf{Lin}, \otimes, K)$ was the original compact closed category; compact closed categories were invented in order to generalize $(\mathsf{Lin}, \otimes, K)$.

The practical upshot of all of this category theory is that we can move vector spaces from the domain to the codomain of a morphism by dualizing them, i.e.

$$\mathrm{Hom}(A \otimes B, C) \cong \mathrm{Hom}(A, C \otimes B^*)$$

This is an extremely common operation in multilinear algebra, and thus satisfying that it has a categorical basis.

## 3.4. Quadratic forms

We now make some definitions of a less abstract kind that are useful in the next chapter. We assume that we work in a field $K$ of characteristic 0.

Definition 3.26. *A **symmetric bilinear form** $f: V \times V \to K$ is a bilinear form such that*

$$f(u,v) = f(v,u)$$

*for all $u, v \in V$.*

Definition 3.27. *A **quadratic form** is a map $q: V \to K$ such that for all $a \in K$, $v \in V$*

$$q(a\,v) = a^2\,q(v)$$

*and the function $(u,v) \mapsto q(u+v) - q(u) - q(v)$ is bilinear.*

Definition 3.28. *Given a quadratic form, its **associated symmetric bilinear form** $\langle -,- \rangle_q$ is defined by*

$$\langle u,v \rangle_q = \frac{1}{2}(q(u+v) - q(u) - q(v))$$

Proposition 3.29. *Given a symmetric bilinear form $f: V \times V \to K$, define $q_f(v) = f(v,v)$. Then $q_f$ is a quadratic form, and $\langle -,- \rangle_{q_f} = f$.*

**Proof.**

$$\begin{aligned}
\langle u,v \rangle_{q_f} &= \frac{1}{2}(q(u+v) - q(u) - q(v)) \\
&= \frac{1}{2}(f(u+v,u+v) - f(u,u) - f(u,u)) \\
&= \frac{1}{2}(f(u,u) + f(u,v) + f(v,u) + f(v,v) - f(u,u) - f(v,v)) \\
&= \frac{1}{2}(f(u,v) + f(v,u)) \\
&= \frac{1}{2}(f(u,v) + f(u,v)) \\
&= f(u,v)
\end{aligned}$$

□

The other direction holds to, so the operations of taking the quadratic form associated with a symmetric bilinear form and taking a symmetric bilinear form associated with a quadratic form are inverses.

Definition 3.30. *Suppose that $q: V \to \mathbb{R}$ is a quadradic form, and suppose that $e_1, \ldots, e_n$ is a basis of $V$ such that*

$$q\left(\sum_{i=1}^{n} a_i\,e_i\right) = \sum_{i=1}^{n} \kappa_i\,a_i^2$$

*where $\kappa_i \in \{0,-1,+1\}$. Then we say that $e_1,\ldots,e_n$ is a **diagonalizing basis** of **signature** $(k,n,m)$ if $k$ of the $\kappa_i$'s are $0$, $n$ of the $\kappa_i$'s are $+1$ and $m$ of the $\kappa_i$'s are $-1$. If $k=0$, then we say $q$ is **nondegenerate**, and we write the signature as $(n,m)$.*

Proposition 3.31. *Every quadratic form has a unique signature (although, there is not a unique diagonalizing basis).*

**Proof.** This can be found in many texts on linear algebra, for instance Shafarevich [53, Chapter 6]. □

Proposition 3.32. *If $e_1, \ldots, e_n$ is a diagonalizing basis for the quadratic form $q: V \to \mathbb{R}$ with coefficients $\kappa_1, \ldots, \kappa_n$ then*

$$\langle e_i, e_j \rangle_q = \begin{cases} \kappa_i & \text{if} \quad i = j \\ 0 & \text{otherwise} \end{cases}$$

The fact that each quadratic form has a unique signature despite the diagonalizing basis non-unique is analogous to how each finite-dimensional vector space has a unique dimension, although the basis that proves that the vector space has a given dimension is non-unique.

## 3.5. Exterior algebra

We end this chapter with a section on exterior algebra, which is an important subject for differential geometry. This material can be found in Aluffi [52, VIII.4.1].

Definition 3.33. *An **alternating bilinear map** $\lambda: V \times V \to W$ is a bilinear map such that $\lambda(v, v) = 0$ for all $v \in V$. An **alternating n-linear map** $\lambda: V^n \to W$ is a $n$-linear map such that $\lambda(v_1, \ldots, v_n) = 0$ if $v_i = v_j$ for $i \neq j$.*

Proposition 3.34. *If $\lambda: V \times V \to W$ is an alternating bilinear map, then*

$$\lambda(v_1, v_2) = -\lambda(v_2, v_1)$$

*More generally if $\lambda: V^n \to W$ is an alternating $n$-linear map, then for any permutation $\sigma \in S_n$,*

$$\lambda(v_1, \ldots, v_n) = \operatorname{sgn}(\sigma)\, \lambda(v_{\sigma(1)}, \ldots, v_{\sigma(n)})$$

**Proof.** For $v_1, v_2 \in V$, if $\lambda$ is an alternating bilinear map, then

$$\begin{aligned} 0 &= \lambda(v_1 + v_2, v_1 + v_2) \\ &= \lambda(v_1, v_1) + \lambda(v_2, v_2) + \lambda(v_1, v_2) + \lambda(v_2, v_1) \\ &= 0 + 0 + \lambda(v_1, v_2) + \lambda(v_2, v_1) \end{aligned}$$

Thus

$$\lambda(v_1, v_2) + \lambda(v_2, v_1) = 0$$

To prove the $n$-linear case, we can factor a permutation into a composition of swaps, and then apply a similar argument. □

Definition 3.35. *The **second exterior power** $\Lambda^2 V$ is a representative for the functor $W \mapsto \{$alternating bilinear maps $V \times V \to W\}$. That is, a linear map $f: \Lambda^2 V \to W$ is equivalent to an alternating bilinear map $V \times V \to W$. In general, the **nth exterior power** is a representative for the functor $W \mapsto \{$alternating $n$-linear maps $V^n \to W\}$.*

**Remark 3.36.** The identity $V \times V \to V \times V$ is equivalent to a map $V \times V \to \Lambda^2 V$; we call this the wedge product and denote it $(v_1, v_2) \mapsto v_1 \wedge v_2$.

**Remark 3.37.** The vector space $\Lambda^n V$ is related to the vector space $V^{\otimes n}$ in two different ways. First of all, we have a natural inclusion of $\Lambda^n V$ into $V^{\otimes n}$ given by the inclusion of alternating $n$-linear maps into $n$-linear maps. This inclusion is important because it identifies $\Lambda^2 V$ with a subspace of $V \otimes V \cong [V^*, V]$. This subspace consists of linear maps $J: V^* \to V$ such that for all $v \in V$,

$$\langle v, J v \rangle = 0$$

This identification is important later on for discussing Poisson structures and Poisson manifolds.

Secondly, we can identify $\Lambda^n V$ with the quotient of $V^{\otimes n}$ by the subspace spanned by all $v_1 \otimes \cdots \otimes v_n$ such that $v_i = v_j$ for some $i \neq j$. This identification is useful because it gives a natural quotient map

$$V^{\otimes n} \to \Lambda^n V$$

Definition 3.38. *The **exterior algebra** of a vector space $V$ is the vector space*

$$\Lambda V := \bigoplus_{n=0}^{\infty} \Lambda^n V$$

*This has the structure of an algebra over $K$ when we equip it with the wedge product*

$$\wedge \colon \Lambda^k V \times \Lambda^l V \to \Lambda^{k+l} V$$

*that is defined via*

$$\Lambda^k V \times \Lambda^l V \to V^{\otimes k} \times V^{\otimes l} \to V^{\otimes (k+l)} \to \Lambda^{k+l} V$$

*where we have used the inclusion of $\Lambda^n V$ into $V^{\otimes n}$, the natural bilinear map $V^{\otimes k} \times V^{\otimes l} \to V^{\otimes (k+l)}$, and the quotient map from $V^{\otimes n}$ into $\Lambda^n V$.*

The exterior algebra is important when discussing *differential forms* in Chapter 5.

# 4. Dirac Relations

## 4.1. Bond graphs and Dirac diagrams

In this thesis, we do not attempt to formally define what a bond graph is, because the formalism that we developed does not exactly fit with the practice of bond graphs. Thus, rather than presenting an account of bond graphs that does not mesh with standard usage, we develop our own notation, which we call "Dirac diagrams" and which is inspired by both bond graphs and string diagrams. Consequently, for our purposes, a bond graph is simply a type of diagram used to represent systems made out of parts that interconnect via power-preserving relations. We do not formalize bond graphs beyond this brief description.

Bond graphs have a long and rich history going back to 1959. An overview of this history can be found in a brief document by H.M. Paytner [54], the man who coined the term "bond graph" and started their use. Many textbooks have been written on bond graphs over the years, from Thomas [55] in 1975 to Cunha and Machado [56] in 2021.

As bond graphs were originally used in an engineering context, it took a while for mathematicians to notice them and attempt formalization. A full history of the mathematics sparked by bond graphs is beyond the scope of this thesis; in lieue of a full history we simply list some of the historical works that have illuminated the path of the present work. The book that introduced the author to the subject of bond graphs was *Thermodynamic Network Analysis of Biological Systems*, by Schnakenberg [3]. However, the main thread of this thesis comes from the work of van der Schaft and collaborators in developing port-Hamiltonian systems. A full history of this development can be found in the first chapter of van der Schaft and Jeltsema [14].

Within applied category theory, Coya [8] attempted to formalize the syntax of bond graphs with a generators and relations point of view, and their semantics with Lagrangian relations. Coya's work differs from the present work in that it does not consider state, i.e. it only covers bond graphs and not port-Hamiltonian systems.

More recently, Lohmayer and Leyendecker [57] have worked on representing bond graphs with undirected wiring diagrams (which we cover in Chapter 7), with the eventual goal of doing a similar formalization to this thesis. However, the categorical details of that formalization are not quite developed; we hope that this thesis ends up being helpful to the eventual full treatment within category theory of their approach. Exergetic port-Hamiltonian systems go farther than this thesis in terms of capturing thermodynamics, so we look forward to a future synthesis.

In this chapter we describe the syntax and semantics of Dirac diagrams, but we defer description of the components that these Dirac diagrams are composing in Chapter 6. The reader who finds this backwards and who has a familiarity with differential geometry is encouraged to skip forwards and then come back here afterwards; some of Chapter 6 can be read independently from this one.

Finally, for this chapter, when we say vector space, we mean finite-dimensional vector space over $\mathbb{R}$.

## 4.2. Bond spaces

A power-transfering interaction between physical systems is represented by a **bond** in Dirac diagrams, as pictured in Figure 4.1.

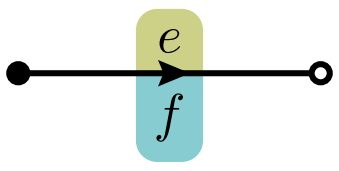

**Figure 4.1.** A single bond

A bond connects an **output port** (filled in) to an **input port** (empty) with a blue flow ($f$) and a green effort ($e$). The meaning of this directionality is that a *positive* power $e \cdot f$ represents energy is flowing from the filled port (output port) to the unfilled port (input port) (i.e., left to right), and a *negative* power represents flow in the other direction. In a full-fledged Dirac diagram, each end of the bond is connected to something; Figure 4.1 is simply a fragment of a Dirac diagram.

In Figure 4.2 we see larger fragments of Dirac diagrams, where bonds are used to connect components with a *state* that changes over time (like charge, flux linkage, position, or momentum). We do not model state until Chapter 6, so a precise characterization of these diagrams is left until then, but it is helpful to know that the flows and efforts are "going somewhere."

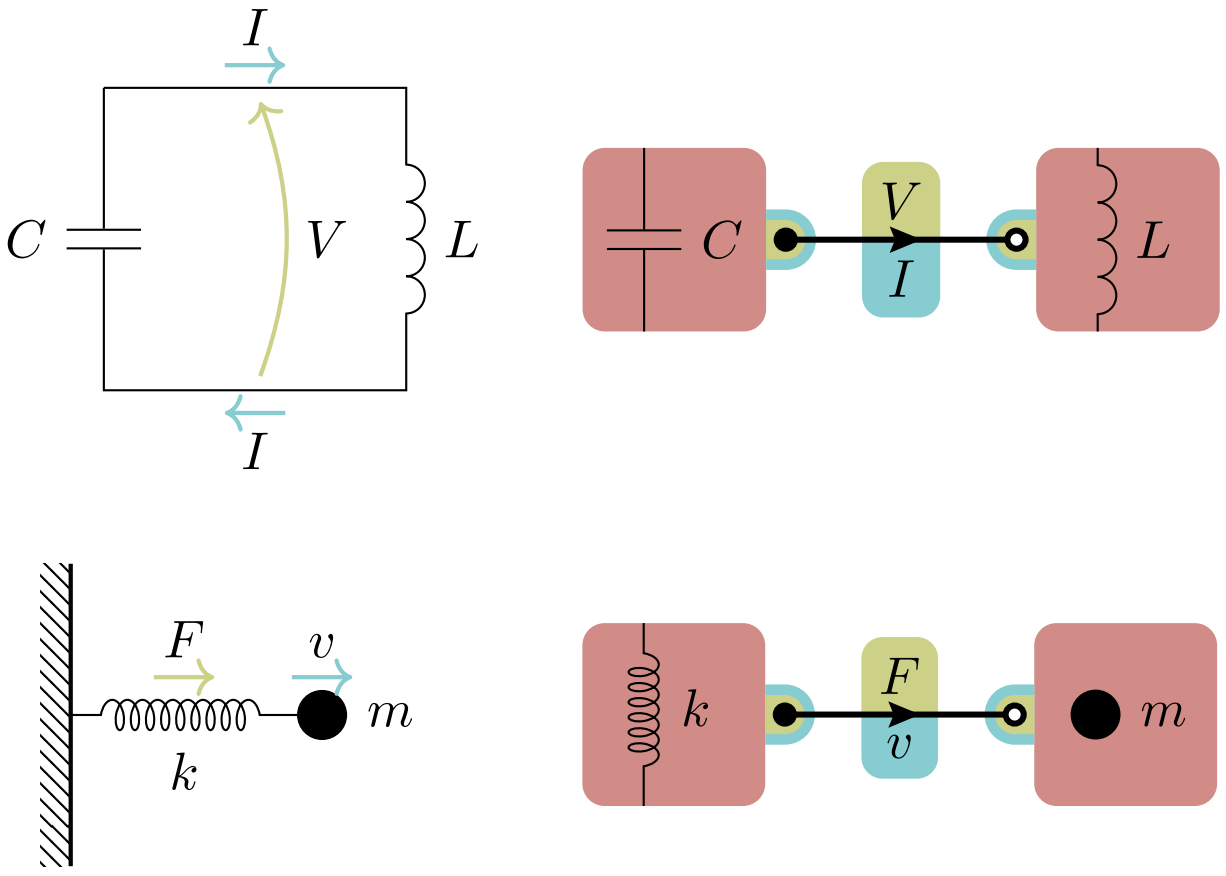

**Figure 4.2.** Traditional notation vs. Dirac diagram. A single bond can represent a pair of wires in an electrical context, or a physical connection between objects in the mechanical context.

Flows and efforts can be vectorial, as systems are often connected by more than one pair of variables. If the space of flows is modeled by the vector space $\mathcal{F}$, then we model the space of efforts and flows with a vector space $\mathcal{E} \oplus \mathcal{F}$, where $\mathcal{E} = \mathcal{F}^*$. The total power is given by a quadratic form $\pi_{\mathcal{E},\mathcal{F}}$ on $\mathcal{E} \oplus \mathcal{F}$ called the **power form**, defined by

$$\pi_{\mathcal{E},\mathcal{F}}(e + f) = \langle e, f \rangle$$

where $e \in \mathcal{E}$, $f \in \mathcal{F}$, and $\langle e, f \rangle \in \mathbb{R}$ is the natural pairing. Note that $\langle -, - \rangle : \mathcal{E} \times \mathcal{F} \to \mathbb{R}$ is a bilinear map when viewed as a two-argument function, but a quadratic form when viewed as a one-argument function.

As there is not a canonical way of saying which variable is flow and which variable is effort (there is only a convention), we seek an "invariant" version of $\mathcal{E} \oplus \mathcal{F}$. This consists of a vector space $V$ and a quadratic form $\pi: V \to \mathbb{R}$ that is isomorphic to $\pi_{\mathcal{E},\mathcal{F}}$, in the sense that it has the same *signature* (see Definition 3.30).

Proposition 4.1. *$\pi_{\mathcal{E},\mathcal{F}}$ is a nondegenerate quadratic form of signature $(n, n)$, where $n = \dim \mathcal{E} = \dim \mathcal{F}$.*

**Proof.** Let $f_1, \ldots, f_n$ be a basis for $\mathcal{F}$ and $e_1, \ldots, e_n$ be the dual basis. Then

$$e_1 + f_1, \ldots, e_n + f_n, e_1 - f_1, \ldots, e_n - f_n$$

diagonalizes $\pi_{\mathcal{E},\mathcal{F}}$ with coefficients $+1,\ldots,+1,-1,\ldots,-1$. □

DEFINITION 4.2. *A **split quadratic form** on a vector space $V$ of dimension $2n$ is a nondegenerate quadratic form on $V$ of signature $(n,n)$.*

DEFINITION 4.3. *A **bond space** is a vector space $V$ along with a split quadratic form $\pi: V \to \mathbb{R}$.*

**Example 4.4.** If $\mathcal{E}$ is a vector space and $\mathcal{F} = \mathcal{E}^*$, then $\mathcal{E} \oplus \mathcal{F}$ with $\pi_{\mathcal{E},\mathcal{F}}(e+f) = \langle e, f\rangle$ is a bond space.

DEFINITION 4.5. *We define the **standard bond space** $\mathbb{B}$ to be the bond space $\mathbb{R}^* \oplus \mathbb{R}$ with power form $\pi_{\mathbb{B}}(e+f) = e\,f$. In general, $\mathbb{B}^n$ is the bond space $(\mathbb{R}^n)^* \oplus \mathbb{R}^n$, with the power form $\pi_{\mathbb{B}^n}(e+f) = \langle e, f\rangle$. The first component is the effort and the second component is the flow.*

DEFINITION 4.6. *If $(V,\pi)$ is a bond space, then $(V,-\pi)$ is also a bond space. When we elide the $\pi$, we write $(V,-\pi)$ as $\bar{V}$.*

$(V,-\pi)$ represents the same variables, but with the opposite convention for direction of power flow. I.e., if $\pi(e+f)$ is the power coming *in*, then $-\pi(e+f)$ is the power going *out*. Eventually, the operation $V \mapsto \bar{V}$ makes a certain category compact closed (see Proposition 4.24), so $V \mapsto \bar{V}$ can be thought of a kind of dual operation.

**Example 4.7.** $\bar{\mathbb{B}}^n$ is the bond space with underlying vector space $(\mathbb{R}^n)^* \oplus \mathbb{R}^n$ and power form $\pi_{\bar{\mathbb{B}}^n}(e+f) = -\langle e, f\rangle$.

We now develop the category of bond spaces and *power-preserving maps*. The function of this category is to provide the symmetric monoidal structure to a category of relations defined later, just as Lin and the biproduct provides a symmetric monoidal structure for LinRel. The rest of this section may be safely skimmed, although it is important from a technical standpoint it does not directly bear on the development of Dirac diagrams.

DEFINITION 4.8. *A power-preserving map from $(V,\pi)$ to $(V',\pi')$ consists of a linear map $L: V \to V'$ such that $\pi = \pi' \circ L$.*

**Example 4.9.** Let $\lambda \in \mathbb{R}$ such that $\lambda \neq 0$, and let $\mathcal{F}$ be any vector space and $\mathcal{E} = \mathcal{F}^*$. Then $e+f \mapsto \lambda e + \frac{1}{\lambda} f$ is a power-preserving map from $\mathcal{E} \oplus \mathcal{F}$ to $\mathcal{E} \oplus \mathcal{F}$. This is an abstract representation of mechanical advantage; devices like levers and pulleys allow one to multiply force while dividing velocity and vice versa.

PROPOSITION 4.10. *The composition of two power-preserving maps is power-preserving.*

DEFINITION 4.11. *The category* Power *has bond spaces as objects and power-preserving maps as morphisms.*

PROPOSITION 4.12. *All bond spaces are isomorphic to $\mathcal{E} \oplus \mathcal{F}$ for some $\mathcal{E} = \mathcal{F}^*$. More specifically, all bond spaces are isomorphic to $\mathbb{B}^n$ for some $n$.*

**Proof.** Let $(V,\pi)$ be a bond space, and let $v_1,\ldots,v_n,w_1,\ldots,w_n$ be a diagonalizing basis for $\pi$. Then map this basis to the basis $e_1+f_1,\ldots,e_n+f_n,e_1-f_1,\ldots,e_n-f_n$ for $(\mathbb{R}^n)^* \oplus \mathbb{R}^n$ where $f_1,\ldots,f_n$ is the standard basis for $\mathbb{R}^n$ and $e_1,\ldots,e_n$ is the dual basis for $(\mathbb{R}^n)^*$. This produces a bijective power-preserving map, as required. □

Although each bond space is isomorphic to one of the form $\mathcal{E} \oplus \mathcal{F}$, there is not a *canonical* isomorphism between a given bond space $(V,\pi)$ and a bond space of the form $\mathcal{E} \oplus \mathcal{F}$. Later on we do constructions that depend on the bond space being in the particular form $\mathcal{E} \oplus \mathcal{F}$.

As an example of this phenomenon, $(\mathcal{E} \oplus \mathcal{F}, \pi_{\mathcal{E},\mathcal{F}})$ is isomorphic to $(\mathcal{E} \oplus \mathcal{F}, -\pi_{\mathcal{E},\mathcal{F}})$ via $e+f \mapsto e-f$ or via $e+f \mapsto -e+f$. In the port-Hamiltonian literature, one of these two isomorphisms is used to get around the need to "change the direction" for a bond space; see [14, Definition 2.3] for an example. The presentation in terms of bond spaces obviates the need to choose one of these isomorphisms, as we can simply work with $(\mathcal{E} \oplus \mathcal{F}, -\pi_{\mathcal{E},\mathcal{F}})$ without needing to make it isomorphic to $(\mathcal{E} \oplus \mathcal{F}, \pi_{\mathcal{E},\mathcal{F}})$.

Proposition 4.13. *There is a symmetric monoidal structure on the category* Power *given by* $(V,\pi)\oplus(V',\pi')=(V\oplus V',\pi+\pi')$*, and with unit* $(\mathbb{R}^0,0)$*, which we denote* $\mathbb{B}^0$*; the 0-dimensional standard bond space.*

**Proof.** For this proof, we follow the following strategy. We first construct a category of vector spaces with quadratic forms on them, via the Grothendieck construction. We then use the work of Moeller and Vasilakopolou [58, Theorem 3.10] to put a symmetric monoidal structure on this category. Finally, we show that the subcategory of vector spaces with split quadratic forms (i.e. Power) is closed under the monoidal operations, and so is a symmetric monoidal category itself.

Consider the functor $Q\colon \mathsf{Lin}^{\mathrm{op}}\to\mathsf{Set}$ defined by

$$Q(V)=\{\text{quadratic forms on } V\}$$

$$Q(L\colon V\to W)(q\colon W\to\mathbb{R})=q\circ L\colon V\to\mathbb{R}$$

The Grothendieck construction of $Q$, $\int Q$, is the category where the objects are vector spaces with quadratic forms, and the morphisms are linear maps that preserve those quadratic forms. Power is the full subcategory consisting of vector spaces with split quadratic forms. We put a symmetric monoidal structure on $\int Q$ by giving $Q$ a lax symmetric monoidal structure and applying Moeller and Vasilakopolou [58, Theorem 3.10].

We now prove that the following maps give a lax symmetric monoidal structure to $Q$:

$$\mu_{X,Y}\colon Q(X)\times Q(Y)\to Q(X\oplus Y)$$

$$\mu_{X,Y}(q_X,q_Y)=q_X+q_Y$$

$$\varepsilon\colon 1\to Q(\mathbb{R}^0)$$

$$\varepsilon(*)=0$$

We can prove this by diagram-chase. The first condition for these maps giving a symmetric monoidal structure is the commutativity of

$$\begin{array}{ccc}
(Q(X)\times Q(Y))\times Q(Z) & \xrightarrow{a_{Q(X),Q(Y),Q(Z)}} & Q(X)\times(Q(Y)\times Q(Z)) \\
\downarrow{\scriptstyle \mu_{X,Y}\times 1_{F(z)}} & & \downarrow{\scriptstyle 1_{Q(X)}\times\mu_{Y,Z}} \\
Q(X\oplus Y)\times Q(Z) & & Q(X)\times Q(Y\oplus Z) \\
\downarrow{\scriptstyle \mu_{X\oplus Y,Z}} & & \downarrow{\scriptstyle \mu_{X,Y\oplus Z}} \\
Q((X\oplus Y)\oplus Z) & \xrightarrow[Q(a_{X,Y,Z})]{} & Q(X\oplus(Y\oplus Z))
\end{array}$$

Starting with $((q_X,q_Y),q_Z)$, it is not hard to see that either way we chase this, we end up with $q_X+q_Y+q_Z$. The second condition is commutativity of the following diagram (and the diagram for the right unitor, not pictured)

$$\begin{array}{ccc}
1\times Q(X) & \xrightarrow{\epsilon\times 1_{Q(X)}} & Q(0)\times Q(X) \\
\downarrow{\scriptstyle l_{Q(X)}} & & \downarrow{\scriptstyle \mu_{0,X}} \\
Q(X) & \xleftarrow[Q(l_X)]{} & Q(0)\times X
\end{array}$$

After chasing these, we find that commutativity is equivalent to $0+q=q=q+0$, and thus holds. Lastly, the diagram for showing that this is compatible with the symmetry is

$$\begin{array}{ccc}
F(X)\otimes_{\mathsf{D}}F(Y) & \xrightarrow{B_{F(X),F(Y)}} & F(Y)\otimes_{\mathsf{D}}F(X) \\
\downarrow{\scriptstyle \mu_{X,Y}} & & \downarrow{\scriptstyle \mu_{Y,X}} \\
F(X\otimes_{\mathsf{C}}Y) & \xrightarrow[F(B_{X,Y})]{} & F(Y\otimes_{\mathsf{C}}X)
\end{array}$$

This can be chased because $q_X+q_Y=q_Y+q_X$.

We have now shown that $\int Q$ is a symmetric monoidal category. To show that Power is also a symmetric monoidal category, note that if $\pi: V \to \mathbb{R}$ is of signature $(n, n)$ and $\pi': V' \to \mathbb{R}$ is of signature $(m, m)$, then $\pi + \pi': V \oplus V' \to \mathbb{R}$ is of signature $(n+m, n+m)$, so Power is closed under the monoidal product and thus a symmetric monoidal category in its own right. □

Note that this is neither a cartesian or cocartesian monoidal structure, because the injection and projection maps are not power-preserving. As mentioned earlier, we do not as of yet know of much use for Power on its own; its main use in this thesis is to provide a natural category in which to discuss isomorphisms of bond spaces and also to provide a symmetric monoidal structure for PowerRel, which we discuss in the next section.

## 4.3. Power-preserving relations

Just as we generalized linear maps to linear relations, we also generalize power-preserving maps to power-preserving relations.

DEFINITION 4.14. *A* ***power-preserving relation*** *between bond spaces* $(V, \pi)$ *and* $(V', \pi')$ *is a relation* $R: V \nrightarrow V'$ *such that if* $(v, v') \in R$, $\pi(v) = \pi(v')$.

PROPOSITION 4.15. *The composite of power-preserving relations is again a power-preserving relation.*

DEFINITION 4.16. *The category* PowerRel *has bond spaces as objects and power-preserving relations as morphisms.*

We now want to make PowerRel into a symmetric monoidal category analogous to LinRel. To do this, we define a functor $\oplus$: PowerRel × PowerRel → PowerRel. On objects of PowerRel, $\oplus$ acts as it does for Power. That is, $(V, \pi_V) \oplus (W, \pi_W) = (V \oplus W, \pi_V + \pi_W)$. Then on morphisms of PowerRel, $\oplus$ acts as it does for LinRel. That is, if $R_1: V_1 \nrightarrow W_1$ and $R_2: V_2 \nrightarrow W_2$ are power-preserving relations between $(V_i, \pi_{V_i})$ and $(W_i, \pi_{W_i})$, then $R_1 \oplus R_2: V_1 \oplus W_1 \nrightarrow V_2 \oplus W_2$ is defined to be the monoidal composition of relations in LinRel; it is not hard to see that this is likewise power-preserving as a relation between $(V_1 \oplus V_2, \pi_{V_1} + \pi_{V_2})$ and $(W_1 \oplus W_2, \pi_{W_1} + \pi_{W_2})$.

Now, just as Lin and $\mathsf{Lin}^{\mathrm{op}}$ are both subcategories of LinRel, so too are Power and $\mathsf{Power}^{\mathrm{op}}$ subcategories of PowerRel, via wide (surjective on objects) embeddings graph: Power → PowerRel and cograph: Power → PowerRel. All of the requisite isomorphisms needed to make $\oplus$ into a symmetric monoidal structure on PowerRel come from this graph embedding, and thus $\oplus$ is a valid symmetric monoidal structure on PowerRel, with monoidal unit $\mathbb{B}^0$ as before.

We now reveal one purpose of Dirac diagrams: Dirac diagrams can be used to represent morphisms in PowerRel.

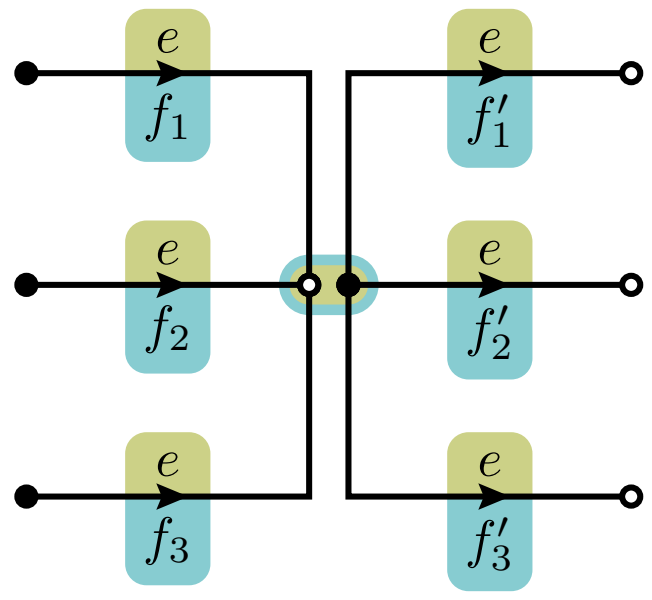

**Figure 4.3.** An effort-matching, flow-summing junction (0-junction). The outside color of the junction (blue) represents the quantity that is summed, and the inside color (green) represents the quantity that is matched. The junction enforces that $f_1 + f_2 + f_3 = f'_1 + f'_2 + f'_3$.

**Example 4.17.** Fix a bond space of the form $(\mathcal{E} \oplus \mathcal{F}, \pi_{\mathcal{E},\mathcal{F}})$. Then the **effort-matching, flow-summing junction** is the relation $J: (\mathcal{E} \oplus \mathcal{F})^{\oplus m} \nrightarrow (\mathcal{E} \oplus \mathcal{F})^{\oplus n}$ defined by

$$((e_1 + f_1, \ldots, e_m + f_m), (e'_1 + f'_1, \ldots, e'_n + f'_n)) \in J$$

if and only if

$$((e_1,\ldots,e_m),(e'_1,\ldots,e'_n))\in\circ \text{ and } ((f_1,\ldots,f_m),(f'_1,\ldots,f'_n))\in\bullet$$

That is, $e_i=e'_j=e$ for some $e$, and $\sum_{i=1}^m f_i=\sum_{j=1}^n f_j$. This is pictured in Figure 4.3. In bond graphs, this is typically called a **0-junction**, and sometimes we use this terminology as well.

This is power-preserving, because for $((e_1+f_1,\ldots,e_m+f_m),(e'_1+f'_1,\ldots,e'_n+f'_n))\in R$, we have

$$\begin{aligned}\sum_{i=1}^m \langle e_i,f_i\rangle &= \left\langle e,\sum_{i=1}^m f_i\right\rangle\\ &= \left\langle e,\sum_{j=1}^n f'_j\right\rangle\\ &= \sum_{j=1}^n \langle e'_j,f'_j\rangle\end{aligned}$$

where $e_i=e=e'_j$ for all $i,j$.

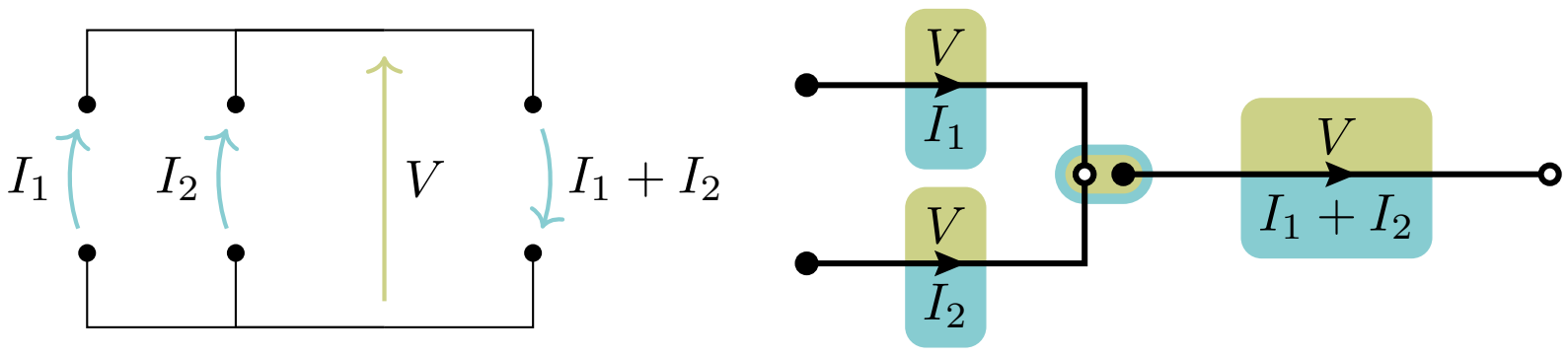

**Figure 4.4.** 0-junction in electronic circuit and Dirac diagram notation. In the gaps marked with blue current arrows, we assume that some circuit component is connected, with the noted current passing through it.

**Example 4.18.** The 0-junction in the electrical context represents composition in parallel. In Figure 4.4, we see the traditional circuit diagram version of a zero junction, along with its Dirac diagram counterpart. Note that the sum of the currents in is equal to the currents out.

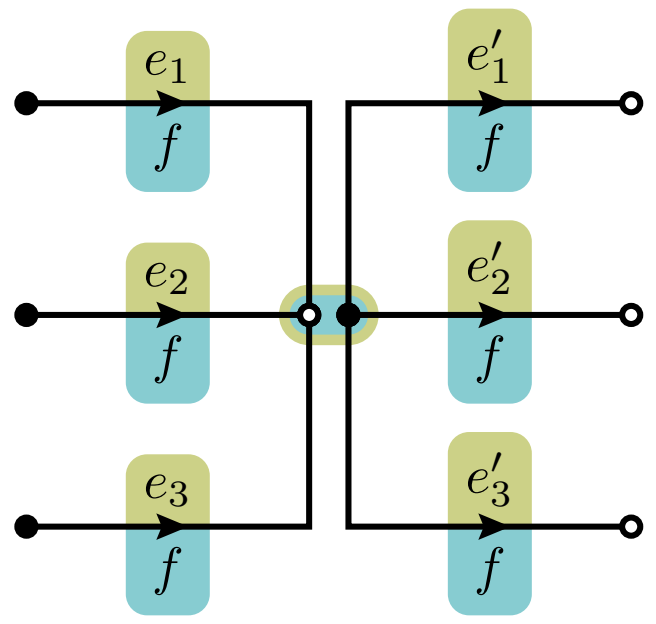

**Figure 4.5.** Flow-matching, effort-summing junction (1-junction). The outside color of the junction (green) represents the quantity that is summed, and the inside color (blue) represents the color that is matched. The junction enforces that $e_1+e_2+e_3=e'_1+e'_2+e'_3$.

**Example 4.19.** Fix a bond space of the form $(\mathcal{E}\oplus\mathcal{F},\pi_{\mathcal{E},\mathcal{F}})$. Then the **flow-matching, effort-summing junction** is the relation $J\colon(\mathcal{E}\oplus\mathcal{F})^{\oplus m}\nrightarrow(\mathcal{E}\oplus\mathcal{F})^{\oplus n}$ defined by

$$((e_1+f_1,\ldots,e_m+f_m),(e'_1+f'_1,\ldots,e'_n+f'_n))\in J$$

if and only if

$$((e_1,\ldots,e_m),(e'_1,\ldots,e'_n))\in\bullet \text{ and } ((f_1,\ldots,f_m),(f'_1,\ldots,f'_n))\in\circ$$

This is the dual of the 0-junction, and is also called the **1-junction**. It is power-preserving by the same logic as the 0-junction.

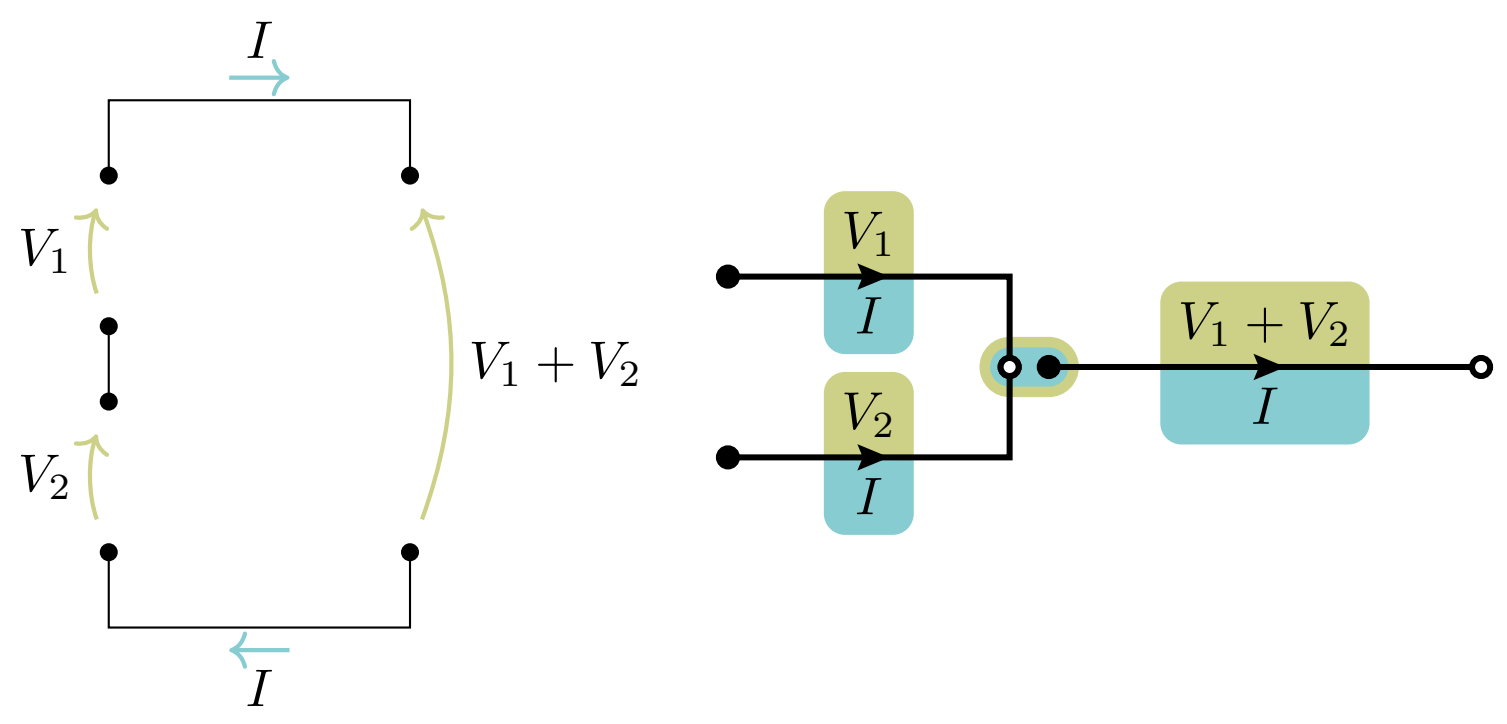

**Figure 4.6.** 1-junction in electronic circuit and Dirac diagram notation.

**Example 4.20.** The 1-junction in the electrical context represents composition in series. In Figure 4.6, we see the traditional circuit diagram version of a 1-junction, along with its Dirac diagram counterpart. Note that the sum of the voltages in is equal to the sum of the voltages out.

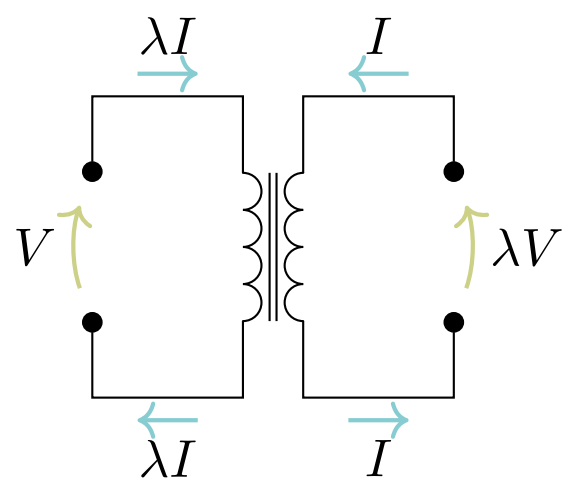

**Figure 4.7.** A circuit for a transformer

**Example 4.21.** Another common example of a power-preserving relation is the graph of $e+f\mapsto \lambda e+\frac{1}{\lambda}f$. This represents mechanical advantage.

## 4.4. Poisson structures and pre-symplectic structures

An important special case of a power preserving relation is one between $\mathbb{B}^0$ and $(V,\pi)$. These consist of subspaces $U\subset V$ with $\pi|_U=0$. Consequently, a power-preserving relation between $(V,\pi_V)$ and $(W,\pi_W)$ can be represented as a power-preserving relation between $\mathbb{B}^0$ and $(V\oplus W,\pi_V-\pi_W)$. We now state this all more formally.

Definition 4.22. *An* ***isotropic subspace*** *of a bond space* $(V,\pi)$ *is a subspace* $U\subset V$ *such that* $\pi|_U=0$.

Proposition 4.23. *A power-preserving relation between V and W is an isotropic subspace of* $\bar{V}\oplus W$ *(where* $\bar{V}$ *is the dual of V, as defined in Definition 4.6. This is a consequence of the compact closed structure on* PowerRel.

The categorical consequences of the previous fact are stated in the following proposition.

Proposition 4.24. PowerRel *is a compact closed category, where the dual of a bond space as defined in Definition 4.6 is the categorical dual.*

**Proof.** Recall from Definition 2.20 that a compact closed category is a symmetric monoidal category in which every object has a dual object.

We show that the the dual of $V$ in PowerRel is $\bar{V}$ (recall Definition 4.6). We must show that this satisfies the properties of a dual object. Recall that in LinRel, every object is its own dual, via the morphisms $\cap\colon 0\nrightarrow V\oplus V$ and $\cup\colon V\oplus V\nrightarrow 0$ both defined by the diagonal $\Delta V=\{(v,v)\mid v\in V\}\subset V\oplus V$. When viewed as maps $\cap\colon 0\nrightarrow V\oplus\bar{V}$ and $\cup\colon\bar{V}\oplus V\nrightarrow 0$, in PowerRel, these relations are power-preserving, because the power of $(v,v)\in V\oplus\bar{V}$ is $\pi(v)-\pi(v)=0$.

It remains to show that the zigzag relations hold, but this is a simple consequence of the forgetful functor from PowerRel to LinRel being faithful. □

Thus, we can study power-preserving relations by studying isotropic subspaces. We now consider some important examples of isotropic subspaces.

Definition 4.25. *Suppose that $X$ is a vector space, and $J\colon X^* \to X$ is a linear map such that $\langle \varphi, J\varphi \rangle = 0$ for all $\varphi \in X^*$. Then we call $J$ a* ***Poisson structure*** *on $X$, and we call*

$$\{(\varphi, J\varphi) \mid \varphi \in X^*\} \subset X^* \oplus X$$

*the* ***induced isotropic subspace*** *of the Poisson structure.*

Using the machinery of Section 3.5, a Poisson structure $J\colon X^* \to X$ is equivalent to an alternating map $X^* \times X^* \to \mathbb{R}$. This is then equivalent to a linear map $\Lambda^2 X^* \to \mathbb{R}$, which corresponds to an element of $\Lambda^2 X$. Thus, a Poisson structure on $X$ is nothing more than a choice of an element of $\Lambda^2 X$.

We can just as easily go the other way as well.

Definition 4.26. *A Poisson structure on $X^*$ is called a* ***pre-symplectic structure*** *on $X$. This is simply a map $J\colon X \to X^*$ such that $\langle Jx, x \rangle = 0$ for all $x \in X$. We call*

$$\{(Jx, x) \mid x \in X\} \subset X^* \oplus X$$

*the* ***induced isotropic subspace*** *of the pre-symplectic structure.*

By similar reasoning as before, a pre-symplectic structure is equivalent to an element of $\Lambda^2 X^*$.

Definition 4.27. *A linear bijection $J\colon X \to X^*$ such that $\langle x, Jx \rangle = 0$ for all $x \in X$ is called a* ***symplectic structure*** *on $X$. Note that this is also a pre-symplectic struction* and *also induces a Poisson structure on $X$. A vector space with a symplectic structure on it is called a* ***symplectic vector space****.*

**Example 4.28.** There is a natural symplectic structure on $X = V \oplus V^*$ for any vector space $V$, given by the map

$$J\colon V \oplus V^* \to (V \oplus V^*)^* \cong V^* \oplus V^{**} \cong V^* \oplus V$$

defined by

$$J(v, \varphi) \mapsto (\varphi, -v)$$

This is a symplectic structure because

$$\langle (v, \varphi), (\varphi, -v) \rangle = \langle \varphi, v \rangle - \langle \varphi, v \rangle = 0$$

Symplectic structures are the foundation of Hamiltonian mechanics, of which port-Hamiltonian mechanics is a generalization, and we further discuss this thread more in the next two chapters.

## 4.5. Dirac relations

In every physically relevant power preserving relation $R\colon V \nrightarrow W$ that we have seen it has been the case

$$\dim R = \frac{1}{2}(\dim V + \dim W)$$

The precise reason for this is not clear to us. However, relations of this form can in fact be characterized in many different ways, which lends a sort of naturality to this condition. For instance, an isotropic subspace is maximal if and only if it is of this dimension. In any event, the literature on port-Hamiltonian systems focuses on relations with this property. Thus, in this chapter we restrict PowerRel to only these relations, which we call Dirac relations.

In several ways, Dirac relations are highly analogous to Lagrangian relations. As a result, for the proofs in this section we were able to follow Fong [7, Chapter 6], who studies Lagrangian relations. It seems like Dirac and Lagrangian relations are somehow dual, though we have not yet figured out how to make this precise.

The main theorem in this section is a proof that the composite of Dirac relations is again a Dirac relation. However, it takes a while to get there, and we must first build some general technology for working with bond spaces and their subspaces. We start with the following definition.

DEFINITION 4.29. *Given a bond space $(V, \pi)$, let the **power bilinear form** $\langle -, - \rangle_+$ be the associated symmetric bilinear form to $\pi$, defined as in Definition 3.28 by*

$$\langle v, w \rangle_+ = \frac{1}{2}(\pi(v+w) - \pi(v) - \pi(w))$$

Note that $\langle v, v \rangle_+ = \pi(v)$. Moreover, for vector spaces $\mathcal{E}$, $\mathcal{F} = \mathcal{E}^*$, if $e \in \mathcal{E}$ and $f \in \mathcal{F}$, and $\langle -, - \rangle_+$ is the power bilinear form corresponding to $\pi_{\mathcal{E},\mathcal{F}}$ then

$$\langle e_1 + f_1, e_2 + f_2 \rangle_+ = \frac{1}{2}\langle e_1, f_2 \rangle + \frac{1}{2}\langle e_2, f_1 \rangle$$

and thus

$$\langle e, f \rangle_+ = \frac{1}{2}(\pi(e+f) - \pi(e) - \pi(f)) = \frac{1}{2}\pi(e+f) = \frac{1}{2}\langle e, f \rangle$$

DEFINITION 4.30. *Given a subspace $U \subset V$ of a bond space $(V, \pi)$, define its **polar** $U^\circ$ by*

$$U^\circ = \{v \in V \mid \forall u \in U, \langle v, u \rangle_+ = 0\}$$

PROPOSITION 4.31. *$U$ is isotropic if and only if $U \subset U^\circ$.*

**Proof.** First of all, if $U \subset U^\circ$, then for all $v \in U$, $\langle v, v \rangle_+ = 0$, so $\pi|_U = 0$ as required. Now, assume that $\pi|_U = 0$. Then for $v, v' \in U$, $v + v'$ and $v - v'$ are both in $U$, so

$$\langle v, v' \rangle = \frac{1}{2}(\pi(v+v') - \pi(v) - \pi(v')) = 0$$

and thus $U \subset U^\circ$. □

PROPOSITION 4.32. *For $U \subset V$,*

$$\dim U + \dim U^\circ = \dim V$$

**Proof.** This follows from the rank-nullity theorem, because if $u_1, \ldots, u_n$ is a basis for $U$, then $U^\circ$ is the kernel of $L: V \to \mathbb{R}^n$ defined by

$$Lv = (\langle u_1, v \rangle_+, \ldots, \langle u_n, v \rangle_+)$$

and $L$ is surjective because $\langle -, - \rangle_+$ is non-degenerate. □

We also have the following identities for subspaces $U, W \subset V$

$$\begin{aligned} (U^\circ)^\circ &= U \\ (U + W)^\circ &= U^\circ \cap W^\circ \\ (U \cap W)^\circ &= U^\circ + W^\circ \end{aligned}$$

DEFINITION 4.33. *$U$ is **coisotropic** if $U^\circ$ is isotropic, or equivalently if $U^\circ \subset U$. $D$ is a **Dirac structure** if $D$ is isotropic and coisotropic.*

PROPOSITION 4.34. *The following are equivalent for a subspace of a bond space $D \subset V$.*

*i. $D$ is a Dirac structure*

*ii. $D$ is a maximal isotropic subspace*

*iii. $D$ is a minimal coisotropic subspace*

*iv.* $D^\circ = D$

*v.* $D$ *is isotropic and* $\dim D = \frac{1}{2}\dim V$

*vi.* $D$ *is coisotropic and* $\dim D = \frac{1}{2}\dim V$

**Example 4.35.** Definitions 4.25, 4.26, and 4.27 all give Dirac structures on $X \oplus X^*$.

Definition 4.36. *A* ***Dirac relation*** *between V and W is a Dirac structure on* $\bar{V} \oplus W$.

Every Dirac relation is also a power-preserving relation, and Examples 4.21 and 4.17 are Dirac relations.

Now, we would like to make the following definition.

Definition 4.37. $\mathsf{DiracRel}$ *is the wide subcategory of* $\mathsf{PowerRel}$ *consisting of Dirac relations. That is,* $\mathsf{DiracRel}$ *is the subcategory of* $\mathsf{PowerRel}$ *containing all of the objects in* $\mathsf{PowerRel}$*, but only the morphisms that are Dirac relations.*

For this to be well-defined, we must show that $\mathsf{DiracRel}$ is closed under composition.

Proposition 4.38. *If* $R$ *is a Dirac relation from* $(U, \pi_U)$ *to* $(V, \pi_V)$*, and* $R'$ *is a Dirac relation from* $(V, \pi_V)$ *to* $(W, \pi_W)$*, then* $R' \circ R$ *is a Dirac relation from* $(U, \pi_U)$ *to* $(W, \pi_W)$*.*

We prove this via a sequence of lemmas.

Lemma 4.39. *Let* $D \subset V$ *be a Dirac structure on a bond space V, and let* $U \subset V$ *be an isotropic subspace of V. Then* $(D \cap U^\circ) + U \subset V$ *is a Dirac structure.*

**Proof.** From Proposition 4.34, a subspace is a Dirac structure if and only if it is equal to its polar. Then we compute using the way that the polar interacts with sums and intersections.

$$\begin{aligned}((D \cap U^\circ) + U)^\circ &= (D \cap U^\circ)^\circ \cap U^\circ \\ &= (D^\circ + (U^\circ)^\circ) \cap U^\circ \\ &= (D + U) \cap U^\circ \\ &= (D \cap U^\circ) + (U \cap U^\circ) \\ &= (D \cap U^\circ) + U\end{aligned}$$

We are done. □

This next lemma states that every Dirac structure has a complementary Dirac structure. This generalizes the following example.

**Example 4.40.** For a vector space $\mathcal{E}$, and for $\mathcal{F} = \mathcal{E}^*$, both $\mathcal{E}$ and $\mathcal{F}$ are Dirac structures on $\mathcal{E} \oplus \mathcal{F}$.

Lemma 4.41. *If* $D \subset V$ *is a Dirac structure, then there exists another Dirac structure* $D' \subset V$ *such that* $D \cap D' = \{0\}$*, and* $D + D' = V$*. Moreover, the map* $v \mapsto 2\langle v, - \rangle_+$ *from* $D'$ *to* $D^*$ *is an isomorphism, and if we extend this map to a map* $D + D' \to D \oplus D^*$*, we get an isomorphism of bond spaces* $V \cong D \oplus D^*$*, where we put the standard bond space structure on* $D \oplus D^*$*.*

**Proof.** We construct $D'$ in the following way. Start with an isotropic subspace $D'$ complementary with $D$. If this has dimension $\frac{1}{2}\dim V$, then we are done. Otherwise, we claim that $D + (D')^\circ = V$. For, suppose that $\langle v, w\rangle = 0$ for all $w \in D + (D')^\circ$. Then $v \in D$ because $D^\circ = D$, and $v \in D'$ because $((D')^\circ)^\circ = D'$. Thus, $v = 0$, so $(D + (D')^\circ)^\circ = \{0\}$, whence $D + (D')^\circ = V$. Consequently, we can write $V$ as

$$V = D \oplus D' \oplus D''$$

where $D'' \subset (D')^\circ$. We can then add any $v \in D''$ to $D'$ to make a larger isotropic subspace disjoint with $D$.

Repeating this process, we get $D'$ such that $D'$ is isotropic and of dimension $\frac{1}{2}\dim V$, so $D'$ is a Dirac structure by Proposition 4.34. Moreover, $D + D' = V$ by disjointness of $D$ and $D'$.

Finally, for $d_1, d_2 \in D$ and $d_1', d_2' \in D'$, we have

$$\langle d_1 + d_1', d_2 + d_2' \rangle_+ = \langle d_1, d_2' \rangle_+ + \langle d_1', d_2 \rangle_+$$

and this gives us our isomorphism of bond spaces $D \oplus D' \cong D \oplus D^*$. □

The above lemma shows that each Dirac structure gives a natural way of breaking a bond space into analogues of effort and flow.

LEMMA 4.42. *Let $U$ be an isotropic subspace of $V$. Then $U^\circ/U$ is a bond space with power form $\pi'(v + U) = \pi(v)$.*

**Proof.** $\pi'$ is well-defined, as for any $u \in U$, $v \in U^\circ$

$$\pi(u + v) = \langle u + v, u + v \rangle = \langle u, u \rangle + 2\langle u, v \rangle + \langle v, v \rangle = \langle v, v \rangle = \pi(v)$$

$\pi'$ is clearly also still a quadratic form. Additionally, the power bilinear form

$$\langle u + U, v + U \rangle_+ = \frac{1}{2}(\pi'(u + v + U) - \pi'(u + U) - \pi'(v + U))$$

is also well-defined.

To show that $\pi$ has the right signature, let $U \subset D \subset U^\circ$ be a Dirac structure, and then pick another Dirac structure $D'$ such that $D \oplus D' = V$, which we can do by Lemma 4.41. Let $D_U = D/U \subset U^\circ/U$ and $D_U' = (D' \cap U^\circ)/U \subset U^\circ/U$; note that $U^\circ/U = D_U \oplus D_U'$.

Now, let $e_1, \ldots, e_n$ be a basis for $D$ such that $U = \text{span}\{e_1, \ldots, e_k\}$, and let $f_1, \ldots, f_n$ is the corresponding dual basis for $D'$, given by the isomorphism $D' \cong D^*$ of Lemma 4.41. Then $e_{k+1} + U, \ldots, e_n + U$ is a basis for $D_U$, and $f_{k+1} + U, \ldots, f_n + U$ is a basis for $D_U'$. Moreover,

$$e_{k+1} + f_{k+1} + U, \ldots, e_n + f_n + U, e_{k+1} - f_{k+1} + U, \ldots, e_n - f_n + U$$

is a diagonalizing basis for $\pi' : U^\circ/U$; it is a basis because $U^\circ/U = D_U \oplus D_U'$, and it is diagonalizing by direct computation, as

$$\langle e_i + U, f_j + U \rangle = \delta_{ij}$$

Thus, $\pi'$ has signature $(0, n-k, n-k)$, as required. □

COROLLARY 4.43. *If $U$ is an isotropic subspace of a bond space $V$, and $D \supset U$ is a Dirac structure, then $D_U = D/U \subset U^\circ/U$ is a Dirac structure on $U^\circ/U$.*

We can now prove Proposition 4.38.

**Proof.** Let $R$ be a Dirac relation from $(U, \pi_U)$ to $(V, \pi_V)$, and $R'$ be a Dirac relation from $(V, \pi_V)$ to $(W, \pi_W)$. We must show that $R' \circ R$ is a Dirac relation from $(U, \pi_U)$ to $(W, \pi_W)$.

Recall that this means that $R$ is a Dirac structure on $\bar{U} \oplus V$, and $R'$ is a Dirac structure on $\bar{V} \oplus W$. Now, let $\Delta$ be the diagonal subspace

$$\Delta = \{(0, v, v, 0) | v \in V\} \subset \bar{U} \oplus V \oplus \bar{V} \oplus W$$

$\Delta$ is isotropic because of the way that the power form is defined on $V \oplus \bar{V}$. The polar of $\Delta$ is

$$\Delta^\circ = \{(u, v, v, w) | u \in U, v \in V, w \in W\} \subset \bar{U} \oplus V \oplus \bar{V} \oplus W$$

Consider the projection map $\Delta^\circ \to \bar{U} \oplus W$. The kernel of this map is $\Delta$, and this map is surjective. Moreover, this map is power-preserving, and thus $\bar{U} \oplus W \cong \Delta^\circ/\Delta$.

We now show that $R' \circ R \subset \bar{U} \oplus W$ is a Dirac structure by applying Lemma 4.39 and Corollary 4.43. Recall that

$$R' \circ R = \{(u, w) \mid \text{there exists } v \in V \text{ such that } (u, v) \in R \text{ and } (v, w) \in R'\}$$

We want to show that this is Dirac. To do this, we start with the Dirac structure $R \oplus R'$ on $\bar{U} \oplus V \oplus \bar{V} \oplus W$:

$$R \oplus R' = \{(u, v, v', w) \mid (u, v) \in R, (v', w) \in R'\}$$

Now, take the intersection of this with $\Delta^\circ$

$$(R \oplus R') \cap \Delta^\circ = \{(u, v, v, w) \mid (u, v) \in R, (v, w) \in R'\}$$

The projection of the above onto $\bar{U} \oplus W$ is clearly $R' \circ R$. But note that this projection is equivalent to quotienting by $\Delta$, i.e.

$$R' \circ R = (((R \oplus R') \cap \Delta^\circ) + \Delta) / \Delta$$

Now, Lemma 4.39 implies that $((R \oplus R') \cap \Delta^\circ) + \Delta$ is a Dirac structure on $\bar{U} \oplus V \oplus \bar{V} \oplus W$, and then Corollary 4.43 implies that $R' \circ R$ is a Dirac structure on $\Delta^\circ / \Delta = \bar{U} \oplus W$, which was what we wanted. □

We now show that not only is DiracRel a subcategory of PowerRel, but in fact it also has all of the nice structure of PowerRel.

Proposition 4.44. DiracRel *inherits the compact closed structure from* PowerRel*.*

**Proof.** To show this, it suffices to show that all of the relevant operations in PowerRel restrict to DiracRel. First of all, as noted before, the direct product preserves Dirac relations. Moreover, any isomorphism in Power is a Dirac relation, so all of the natural isomorphisms for the symmetric monoidal structure that come from Power are still present in DiracRel. Finally, the cap and the cup that characterize the duals in PowerRel are Dirac relations, and so the compact closed structure carries over as well. □

We finish this chapter by summarizing what we have developed. We wanted to characterize connections between physical systems that exchanged energy. To do this, we started with formalizing a specification of a interface through which power can travel using *bond spaces*, which are vector spaces that represent abstract effort and flow. We then formalized connections between the interfaces of systems using *power-preserving relations*, which are linear relations that represent interactions where power is conserved. Finally, we identified an important subclass of power-preserving relations which we called Dirac relations. The category of Dirac relations has a lot of nice structure: specifically it is a compact closed category. We represent Dirac relations with Dirac diagrams.

The next chapter takes what we have done and lifts it to a nonlinear context, while additionally providing the foundation for port-Hamiltonian systems.

# 5. Categorical Differential Geometry

## 5.1. Vector bundles

Up to now we have been primarily concerned with vector spaces, linear maps and linear relations. However, many systems are not linear. In order to deal with these non-linearities, we use differential geometry. In this chapter, we review some relevant definitions, as always from a category-theoretic viewpoint.

We assume that the reader is familiar with basic definitions in differential geometry, such as manifolds, tangent and cotangent vectors, and smooth maps, and we seek to place these definitions into a categorical framework. References for this material can be found in Frankel [59] or Baez [60].

Definition 5.1. *The category* Man *has as objects manifolds and as morphisms smooth maps. For each $n$,* Man *has a full subcategory* $\mathsf{Man}_n$ *consisting of $n$-dimensional manifolds.*

Definition 5.2. *For $X$ a manifold, we define* $\mathsf{Bund}(X) = \mathsf{Man}/X$ *to be the category of* ***bundles*** *over $X$. An object of this category is a manifold $E$ with a map $p: E \to X$, and a morphism from $p: E \to X$ to $p': E' \to X$ is a map $f: E \to E'$ such that the following diagram commutes.*

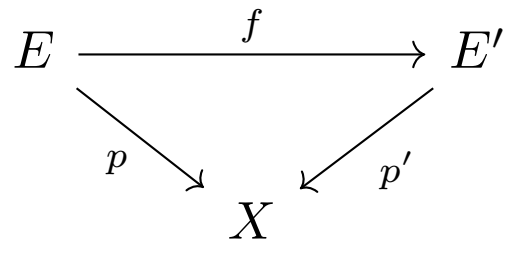

We think of a bundle as a manifold "smoothly parameterized" by $X$. That is, for each point $x \in X$, there is a space $E_x = p^{-1}(x)$ called the **fiber** that "smoothly depends" on $x$, or is "modulated" by $x$. The simplest version of this is given by the so-called trivial bundle, where $E_x$ is constant.

Definition 5.3. *There is a functor* $\mathrm{Triv}_X\colon \mathsf{Man} \to \mathsf{Bund}(X)$ *defined by* $\mathrm{Triv}_X(A) = X \times A \xrightarrow{\pi_1} X$. $\mathrm{Triv}_X(A)$ *is called the* ***trivial bundle with fiber A****.*

A generalization of this is given by the case where the bundle is only *locally* trivial.

Definition 5.4. *For $X$ a manifold, we may consider the subcategory $\mathsf{LTBund}(X) \subset \mathsf{Bund}(X)$ consisting of* ***locally trivial bundles*** *on $X$. A locally trivial bundle on $X$ is a bundle $p\colon E \to X$ along with a manifold $F$ called the* ***standard fiber*** *such that for all $x \in X$ there exists an open neighborhood $U \ni x$ and an isomorphism $U \times F \cong p^{-1}(U)$ in $\mathsf{Bund}(U)$. That is, an isomorphism $f$: $U \times F \to p^{-1}(U)$ such that*

$$\begin{array}{ccc} U \times F & \xrightarrow{\quad f \quad} & p^{-1}(U) \\ & {}_{\pi_1}\searrow \quad \swarrow_{p} & \\ & U & \end{array}$$

*is a commutative triangle.*

Note that $F \cong E_x$ for any $x$. Now, what's the use of a manifold that "varies smoothly" with points in $X$, if the manifold is the same for every point? Classically, one reason for this is that such bundles can have interesting topological structure. For instance, a Möbius strip is an interesting non-trivial vector bundle over $S^1$. However, even with trivial bundles, there still is a reason to think about them, because as always in category theory, we care about the *morphisms*, not the objects. A morphism betweeen bundles is a "smoothly parameterized" map $f_x\colon E_x \to E'_x$. If $E_x \cong \mathbb{R}^n$ for every $x$, $E'_x \cong \mathbb{R}^m$, and $f_x$ is linear for every $x$, then this allows us to use the tools of linear algebra, even when $f_x$ depends non-linearly on $x$. Such bundles are called *vector bundles*, and the category of vector bundles on $X$, $\mathsf{Lin}(X)$, has properties similar to $\mathsf{Lin}$. We now prepare to define what a vector bundle is.

Proposition 5.5. *The category $\mathsf{LTBund}(X)$ has all finite products.*

**Proof.** If $\mathsf{Man}$ had all (categorical) pullbacks, then products in $\mathsf{Bund}(X) = \mathsf{Man}/X$ would be computed by pullback. However this is not the case. That being said, given two locally trivial fiber bundles $p_1\colon E_1 \to X$, $p_2\colon E_2 \to X$ with canonical fibers $F$ and $F'$ we can construct their product by taking the pullback $p\colon E_1 \times_X E_2 \to X$ in *topological* spaces, and then putting the smooth structure on it given by isomorphisms $F \times F' \times U \cong p^{-1}(U)$. □

Definition 5.6. *If $X$ is a manifold, then a* ***rank-n vector bundle*** *over $X$ is a locally trivial fiber bundle $E$ with canonical fiber $\mathbb{R}^n$, along with maps in $\mathsf{LTBund}(X)$ (remember, these are fiberwise maps)*

- $(- + -)\colon E \times E \to E$ *(addition)*
- $(- \cdot -)\colon \mathrm{Triv}(\mathbb{R}) \times E \to E$ *(scalar multiplication)*

*such that the vector space axioms hold, i.e.*

- $(E, +)$ *is an abelian group*
- $\lambda \cdot (e_1 + e_2) = \lambda \cdot e_1 + \lambda \cdot e_2$
- $\lambda \cdot (\gamma \cdot e) = (\lambda \gamma) \cdot e$

*We also require that for every point $x \in X$, there exists a neighborhood $U \ni x$ and a local trivialization isomorphism $U \times \mathbb{R}^n \cong p^{-1}(U)$ that preserves the linear structure. A morphism $L\colon E \to E'$ of vector bundles is then a bundle morphism that respects the vector space structure. The category of vector bundles and vector bundle morphisms on $X$ we call $\mathsf{Lin}(X)$.*

We think of $\mathsf{Lin}(X)$ as "the category of vector spaces parameterized by $X$". $\mathsf{Lin}(X)$ has a distinguished object, which plays a crucial role.

**Example 5.7.** For a manifold $X$, the **tangent bundle** $p: TX \to X$ is defined by

$$TX = \{(x, v) \mid x \in X, v \in T_x X\}, \quad p(x, v) = x$$

where $T_x X$ is the **tangent space** to $X$ at a point $x$. The reader can refer to Frankel or Baez for the definition of tangent space, and for how to put a manifold structure on $TX$. If $X$ is $n$-dimensional, then $TX$ is a rank-$n$ vector bundle.

We now show how all of our work putting interesting structures on $\mathsf{Lin}$ can be neatly extended to $\mathsf{Lin}(X)$. This comes from the theory of *smooth functors*, as described in Szilasi [61]. Without going into too much detail, we give an overview of the idea. Our presentation differs from Szilasi, so as to harmonize better with the other chapters of this thesis.

Definition 5.8. *Let $\mathsf{Lin}_{\cong}$ be the category of vector spaces and vector space isomorphisms.*

Definition 5.9. *A **smooth functor** is a functor $F$ from $\mathsf{Lin}_{\cong}^m$ to $\mathsf{Lin}$ such that for all vector spaces $V_1, \ldots, V_m, W_1, \ldots, W_m$, the map*

$$F: \mathsf{Lin}_{\cong}(V_1, W_1) \times \cdots \times \mathsf{Lin}_{\cong}(V_m, W_m) \to \mathsf{Lin}(F(V_1, \ldots, V_m), F(W_1, \ldots, W_m))$$

*is* smooth *with respect to the vector space structures on the domain and codomain.*

Smooth functors describe constructions on vector spaces that can be lifted to the category of vector bundles on a manifold. The fact that the map on morphisms is smooth roughly means that we can apply it to the transition maps between charts in a manifold.

**Example 5.10.** The tensor product $\otimes: \mathsf{Lin}_{\cong} \times \mathsf{Lin}_{\cong} \to \mathsf{Lin}$ is a smooth functor.

We use $\mathsf{Lin}_{\cong}$ because it allows us to handle a greater variety of operations. For instance, contravariant functors are no problem because we can take inverses.

**Example 5.11.** Define a smooth functor $(-^{-1})^*: \mathsf{Lin}_{\cong} \to \mathsf{Lin}$ by sending $V$ to its dual space $V^*$, and $L: V \to W$ to $(L^{-1})^*: V^* \to W^*$.

The same trick can be applied to functors of multiple arguments that are covariant in one and contravariant in the other.

**Example 5.12.** The closed structure $[-,-]: \mathsf{Lin}_{\cong} \times \mathsf{Lin}_{\cong} \to \mathsf{Lin}$ is a smooth functor.

**Example 5.13.** The functor $V \mapsto [V, V]$ is a smooth functor $\mathsf{Lin}_{\cong} \to \mathsf{Lin}$. Note that as $[-,-]$ is contravariant in its first argument and covariant in its second, it would not be possible to define a contravariant *or* covariant functor from $\mathsf{Lin}$ to $\mathsf{Lin}$ with this value on objects.

**Example 5.14.** The direct sum $\oplus: \mathsf{Lin}_{\cong} \times \mathsf{Lin}_{\cong} \to \mathsf{Lin}$ is a smooth functor.

**Example 5.15.** For any vector space $V$, there is a functor $1 \to \mathsf{Lin}$ that picks out $V$ and is trivially smooth.

**Example 5.16.** The $n$th exterior power $\Lambda^n: \mathsf{Lin}_{\cong} \to \mathsf{Lin}$ is a smooth functor.

Proposition 5.17. *If $X$ is a manifold, then any smooth functor $F: \mathsf{Lin}_{\cong}^m \to \mathsf{Lin}$ induces a functor $\mathsf{Lin}(X)(F): \mathsf{Lin}(X)_{\cong}^m \to \mathsf{Lin}(X)$.*

**Proof.** Conceptually, this is done "fiberwise". That is, if $E_1, \ldots, E_n$ are vector bundles over $X$ with projection maps $p_i: E_i \to X$ and $F$ is a smooth functor, then we construct a vector bundle $\mathsf{Lin}(X)(F)(E_1, \ldots, E_n)$ by letting the fiber over $x \in X$ be the vector space

$$F(p_1^{-1}(x), \ldots, p_n^{-1}(x))$$

The fact that $F$ is a smooth functor then allows these fibers to be stitched together into a vector bundle, see Szilasi [61] for details. □

This proposition allows us to "lift" almost all of the constructions that we performed in Chapter 3 into the world of vector bundles. We perform this lifting implicitly, for instance writing $E \oplus E'$ instead of $\mathsf{Lin}(X)(\oplus)(E, E')$. For instance, this can be applied to all of Examples 5.11, 5.10, 5.12, 5.14, 5.15, 5.16; we freely use the notation $E^*$, $E \otimes E'$, $[E, E']$, $E \oplus E'$, and $\Lambda^n E$ for vector bundles $E$.

The following proposition extends Proposition 5.17 to also apply to the construction of maps between vector spaces.

Proposition 5.18. *If $F, G\colon \mathsf{Lin}^m_{\cong} \to \mathsf{Lin}$ are smooth functors, and $\alpha\colon F \Rightarrow G$ is a natural transformation, then for any smooth manifold $X$ we have a natural transformation*

$$\mathsf{Lin}(X)(\alpha)\colon \mathsf{Lin}(X)(F) \Rightarrow \mathsf{Lin}(X)(G)$$

*defined fiberwise by $\alpha$. That is, if $E_1, \ldots, E_n$ are vector bundles over $X$ with projection maps $p_i$: $E_i \to X$, and $x \in X$, then the fiber of*

$$\mathsf{Lin}(X)(\alpha)_{E_1, \ldots, E_n}\colon \mathsf{Lin}(X)(F)(E_1, \ldots, E_n) \to \mathsf{Lin}(X)(G)(E_1, \ldots, E_n)$$

*over $x$ is*

$$\alpha_{p_1^{-1}(x), \ldots, p_n^{-1}(x)}\colon F(p_1^{-1}(x), \ldots, p_n^{-1}(x)) \to G(p_1^{-1}(x), \ldots, p_n^{-1}(x))$$

**Example 5.19.** Any linear map $L$ from a vector space $V$ to a vector space $W$ can be lifted by Proposition 5.18 to a bundle map $V \times X \to W \times X$.

**Example 5.20.** The natural map $\mathbb{R} \to V \otimes V^*$ given by sending $1 \in \mathbb{R}$ to the identity in $[V, V] \cong V \otimes V^*$ can be lifted to a natural map $\mathbb{R} \times X \to E \otimes E^*$ for any vector bundle $E$.

**Example 5.21.** The natural map $[V, W] \to [W^*, V^*]$ given by taking the dual of a linear map can be lifted to a map of bundles $[E_1, E_2] \to [E_2^*, E_1^*]$.

As this lifting is defined fiberwise, it clearly respects composition. Thus, all the operations, natural isomorphisms, and coherence conditions that give $\mathsf{Lin}$ its two monoidal structures carry directly into $\mathsf{Lin}(X)$.

Proposition 5.22. *For any manifold $X$, $(\mathsf{Lin}(X), \oplus, \mathbb{R}^0 \times X)$ is a symmetric monoidal category (see Proposition 3.2).*

Proposition 5.23. *For any manifold $X$, $(\mathsf{Lin}(X), \otimes, \mathbb{R} \times X)$ is a symmetric monoidal category (see Proposition 3.17).*

Proposition 5.24. *For any manifold $X$, $(\mathsf{Lin}(X), \otimes, \mathbb{R} \times X)$ is a compact closed category (see Corollary 3.25).*

As mentioned before, the original use of string diagrams was for tensor algebra. In fact, it was for tensor algebra that varied over a manifold, i.e. this precise category. However, so the reader does not become confused, we should point out that the "cap and cup" in the setting of Penrose graphical notation typically referred to the metric tensor (discussed in the next section), not the dual maps.

Now that we have discussed a wide variety of ways of constructing vector bundles, the reader might be wondering what one does with them after they have been constructed. The next section gives one answer to this question.

## 5.2. Sections of vector bundles

Just as the elements $v$ of a vector space $V$ are represented by linear maps out of $\mathbb{R}$ which send 1 to $v$, so too we might think of the "elements" of a vector *bundle* $E$ over $X$ as being maps from $\mathbb{R} \times X$. However, these are not simply the elements of $E$ considered as a manifold, rather this is a choice of an element in $E_x$ for each $x$. We call this a *section*.

Definition 5.25. *A* ***section*** *of a vector bundle* $p: E \to X$ *consists of a smooth function* $\varphi: X \to E$ *such that* $p \circ \varphi = 1_X$*. The space of sections we call* $\Gamma(E)$*.*

Proposition 5.26. $\Gamma(E)$ *has a natural vector space structure on it, given by pointwise operations.*

Proposition 5.27. $\Gamma(E)$ *is naturally isomorphic to* $\mathrm{Hom}_{\mathsf{Lin}(X)}(\mathbb{R} \times X, E)$*.*

**Proof.** A vector bundle morphism $f: \mathbb{R} \times X \to E$ induces a section $x \mapsto f(1, x)$, and moreover any section $\varphi$ induces a vector bundle morphism via $f(\lambda, x) = \lambda\, \varphi(x)$. □

As the analogous object to "elements of a vector space", "sections of a vector bundle" play an immensely important role in differential geometry, and by extension, any field that relies on it.

Definition 5.28. *A* ***dynamical system*** *consists of a manifold* $X$ *along with a section* $v$ *of the tangent bundle* $TX$*. Such a section is typically called a* ***vector field****. A* ***trajectory*** *of such a dynamical system consists of a smooth map* $\gamma: [a, b] \to X$ *such that* $\gamma'(t) = v(\gamma(t))$ *for all* $t \in [a, b]$*.*

Dynamical systems are a particularly simple type of ordinary differential equation, and are often written in a more equational form as something like

$$\dot{x} = v(x)$$

A solution of this equation is precisely a trajectory of the dynamical system.

**Example 5.29.** Let SBF: $\mathsf{Lin}_{\cong} \to \mathsf{Lin}$ be the functor that sends a vector space to the vector space of symmetric bilinear forms on that vector space. A **Riemannian manifold** is a manifold $X$ equipped with a section $g$ of the bundle $\mathrm{SBF}(TX)$, such that $g(x): T_x X \times T_x X \to \mathbb{R}$ is positive-definite for each $x$. The section $g$ is known as the **metric tensor**. The study of Riemannian manifolds is called *Riemannian geometry*.

There are two particularly important classes of sections in differential geometry.

Definition 5.30. *A* ***rank*** $(n, m)$ ***tensor*** *over a manifold* $X$ *is a section of the vector bundle* $(TX)^{\otimes n} \otimes (T^* X)^{\otimes m}$*.*

**Example 5.31.** The metric tensor on a Riemannian manifold is a rank $(0, 2)$ tensor.

Definition 5.32. *A* ***differential form of order*** $n$ *over a manifold* $X$ *is a section of the vector bundle* $\Lambda^n T^* X$*. The space of differential forms of order* $n$ *we call* $\Omega^n X = \Gamma(\Lambda^n T^* X)$*.*

## 5.3. Nonlinear Dirac relations

Just as we can discuss vector bundle maps, which we think of as linear maps that are nonlinearly parametrized by a manifold $X$, we can also discuss relations between vector bundles, which we think of as linear relations that are nonlinearly parametrized by a manifold $X$. A relation between vector spaces $V$ and $W$ is a linear subspace of $V \oplus W$; the exact same definition works for vector bundles.

Definition 5.33. *A* ***linear sub-bundle*** *of a vector bundle* $p: E \to X$ *is a subobject of* $E$ *in* $\mathsf{Lin}(X)$*, that is a bundle* $E'$ *along with a monomorphism* $E' \to E$*.*

Definition 5.34. *A* ***linear bundle relation*** *between a vector bundle* $E_1$ *and a vector bundle* $E_2$ *over the same manifold* $X$ *is a linear sub-bundle* $R \subset E_1 \oplus E_2$*.*

Now, we would like to simply show that $\mathsf{Lin}(X)$ is a regular category and then follow the standard construction to make a category of linear bundle relations. However, this does not work; $\mathsf{Lin}(X)$ is not a regular category. This is because maps in $\mathsf{Lin}(X)$ are not necessarily of constant rank, so the image/kernel of a map in $\mathsf{Lin}(X)$ may not have constant dimension, and thus might not be a vector bundle. Thus, the naive composite of two linear bundle relations (i.e. via composing the linear relations fiberwise) might not even be a vector bundle, as it might not have constant dimension across its fibers.

On the other hand, the composite of two Dirac relations of given dimensions is of a fixed dimension, because it is also a Dirac relation. Thus, we might hope that if we restrict our attention only to Dirac relations, this constant-rank property might lead to composition being well-defined. Unfortunately, it is unclear that the composite is a submanifold, and we have not yet been able to prove it.

It could be that we are missing assumptions; Dirac structures in the setting of differential geometry are quite complex to define and have additional properties. To treat them properly, we would have to talk about some geometric properties that are absent in the linear treatment (i.e. integrability), and this would lead us farther into differential geometry than we would like to travel in this thesis. The interested reader can reference Merker [62] for more information on this subject.

Fortunately, for what we do in this thesis, we do not need any of these geometric properties. We thus go forwards with a not completely satisfactory treatment of the subject, which is nevertheless sufficient for what we want to prove later.

Definition 5.35. *Let $Q\colon \mathsf{Lin}_{\cong} \to \mathsf{Lin}$ send a vector space $V$ to the vector space of quadratic forms on $V$.*

Definition 5.36. *A **bond bundle** over a manifold $X$ is a vector bundle $E$ over $X$ along with a section $\pi$ of $Q(E)$ called the **power form** such that $\pi_x\colon E_x \to \mathbb{R}$ is a split quadratic form for all $x \in X$.*

**Example 5.37.** For any vector bundle $E$ over a manifold $X$, $E \oplus E^*$ is a bond bundle with the power form

$$\pi_x(e + f) = \langle e, f \rangle$$

where $e \in E_x$, $f \in E^*_x$.

**Example 5.38.** For any manifold $X$, the bundle $TX \oplus T^*X$ is a bond bundle, equipped with the power form from Example 5.37.

**Example 5.39.** If $X$ is a manifold, and $(V, \pi)$ is a bond space, then $(V \times X, \pi)$ is a bond bundle over $X$.

Definition 5.40. *If $B$ is a bond bundle over a manifold $X$ with power form $\pi$, then define $\bar{B}$ to be the bond bundle with underlying vector bundle the same as $B$ and power form $-\pi$.*

Definition 5.41. *The category $\mathsf{Power}(X)$ has as objects bond bundles over $X$ and as morphisms power-preserving vector bundle morphisms.*

We now give our definition of a Dirac relation between two bond bundles.

Definition 5.42. *If $A$ and $B$ are two bond bundles over a manifold $X$, then a **Dirac relation** $D\colon A \nrightarrow B$ between them consists of a Dirac relation $D_x \subset \bar{A}_x \oplus B_x$ for every $x \in X$.*

Note that we do not require that $D$ is itself a linear subbundle of $\bar{A} \oplus B$; this is why the definition is not quite adequate. This makes the following proposition trivial to prove.

Proposition 5.43. *If $A$, $B$, and $C$ are all bond bundles over a manifold $X$, and $D\colon A \nrightarrow B$ and $D'\colon B \nrightarrow C$ are both Dirac relations, then $D' \circ D\colon A \nrightarrow C$ is also a Dirac relation, where $(D' \circ D)_x = D'_x \circ D_x$.*

Definition 5.44. *The category $\mathsf{DiracRel}(X)$ has as objects bond bundles over $X$ and as morphisms Dirac relations between them.*

Proposition 5.45. *There is a functor $\mathrm{Triv}_X\colon \mathsf{DiracRel} \to \mathsf{DiracRel}(X)$ that sends a bond space $(V, \pi)$ to the bond bundle $(V \times X, \pi)$, and a Dirac relation $R \subset \bar{V} \oplus W$ to the Dirac bundle relation $R \times X \subset (\bar{V} \oplus W) \times X$.*

# 6. Port-Hamiltonian Systems

## 6.1. Hamiltonian mechanics

Before discussing port-Hamiltonian systems, we give a brief review of Hamiltonian mechanics. Port-Hamiltonian systems generalize Hamiltonian mechanics in several ways, but Hamiltonian mechanics remain an important special case. The reader can refer to Sussman and Wisdom [63], Taylor [64] or Arnol'd [65] for an account of the material of this section in much greater depth.

The object of Hamiltonian mechanics is to derive *system dynamics* for a *closed*, *energy-conserving system* from the *energy* of the system. To put this more mathematically, we are given a state space $X$ (which we assume to be a manifold) and a smooth function $H\colon X\to\mathbb{R}$ called the *Hamiltonian*, which is supposed to model the energy of a point in the state space. We want to derive a vector field $\xi_H\in\Gamma(TX)$ such that the system dynamics are described by the differential equation

$$\dot{x}=\xi_H(x) \tag{6.1}$$

Without any additional information, this is not in general possible. However, it turns out that there is a certain structure that we can put on our state space $X$, natural in many cases, that gives us a way of taking any smooth function $H\colon X\to\mathbb{R}$ and producing a vector field $\xi_H\in\Gamma(TX)$. This structure is called a *Poisson structure*, but before we get to defining it, we give examples of this setup.

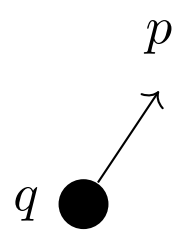

**Figure 6.1.** A free particle at position $q$ with momentum $p$.

**Example 6.1.** A particle of mass $m$ in 3-dimensional space has state space $X=\mathbb{R}^6$, with coordinates $(q,p)$ describing position and momentum respectively. The more traditional (and isomorphic) state space is $T^*\mathbb{R}^3$. Suppose that this particle is acted on by a conservative force field $F$. A typical Hamiltonian for such a particle is given by

$$H(q,p)=\frac{1}{2\,m}\|p\|^2+U(q)$$

where $\frac{1}{2m}\|p\|^2$ is called the **kinetic energy** and $U(q)$ is called the **potential energy**, defined by

$$U(q)=-\int_a^b F(x)\cdot\gamma'(t)\,\mathrm{d}t$$

where $\gamma\colon[a,b]\to\mathbb{R}^3$ starts at a fixed point $q_0$ and ends at $q$; conservativity of $F$ implies that it does not matter which $\gamma$ we choose, though choosing a different $q_0$ would change $U$ by a constant. The reason there is a negative sign is that $F$ is the force applied to the particle by the field, but when want to calculate the potential energy (i.e., the energy "stored in" the field), we say that the potential energy is the work done on the field by the particle, and thus there is a negative sign.

Note that as a result of the fundamental theorem of calculus, $F=-\nabla U$, and thus $F=-\frac{\partial H}{\partial q}$. That is, we can deduce the force from the Hamiltonian. This observation is the key to Hamiltonian mechanics.

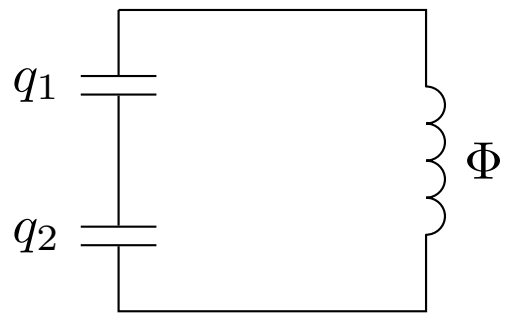

**Figure 6.2.** A circuit with charges $q_1$ and $q_2$ on two capacitors, and flux linkage $\Phi$ through an inductor

**Example 6.2.** Consider the circuit diagram in Figure 6.2. The state space for this is $X = \mathbb{R}^3$, as the state of the circuit can be described by the charges $q_1$ and $q_2$ on each of the capacitors and the flux linkage $\Phi$ on the inductor. Then the Hamiltonian is

$$H(q_1, q_2, \Phi) = \frac{1}{2C_1} q_1^2 + \frac{1}{2\,C_2} q_2^2 + \frac{1}{2L} \Phi^2$$

where $C_1$ and $C_2$ are the respective capacitances of the capacitors and $L$ is the inductance of the inductor.

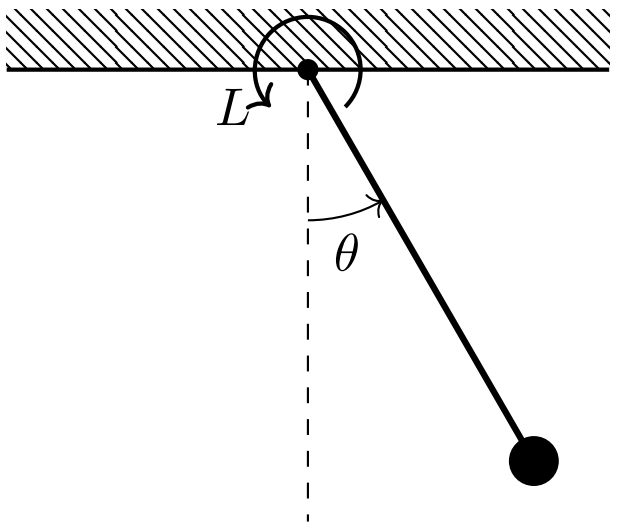

**Figure 6.3.** Pendulum at angle $\theta$ with angular momentum $L$

**Example 6.3.** A pendulum with a rod of length $l$ and mass $m$ concentrated at the tip is described by its angle $\theta$ and angular momentum $L$, and has Hamiltonian

$$H(\theta, L) = m\,g\,l \sin\theta + \frac{1}{2\,m\,l^2} L^2$$

**Example 6.4.** Consider $n$ point masses attached to each other with rigid rods, spinning freely in space around their common center of mass. This rigid body is described by vectors measuring the displacement of each particle from the center of mass $\{r_k \in \mathbb{R}^3\}$, and the mass of each particle $\{m_k \in \mathbb{R}_{>0}\}$. The state of this rigid body is given by a rotation $\Theta \in \mathrm{SO}(3)$ and an angular momentum $L \in T^*_\Theta(\mathrm{SO}(3))$, so the state space is $X = T^*(\mathrm{SO}(3))$. To compute the kinetic energy, we recall that the angular momentum $L$ is related to the angular velocity $\omega$ via the equation

$$L = \tilde{I}\,\omega$$

where $\tilde{I}$ is the **moment of inertia tensor**, which can be computed from the $r_k$ and the $m_k$. Formally speaking, $\tilde{I}$ is a section of $[T(\mathrm{SO}(3)), T^*(\mathrm{SO}(3))]$, giving an invertible linear map $T_\Theta(\mathrm{SO}(3)) \to T^*_\Theta(\mathrm{SO}(3))$ for each $\Theta \in \mathrm{SO}(3)$. Now, recalling that kinetic energy is $\frac{1}{2}\langle w, L\rangle$, we can now write

$$H(\Theta, L) = \frac{1}{2}\langle \tilde{I}_\Theta^{-1}(L), L\rangle$$

Thus, once we have used $r_k$ and $m_k$ to compute $\tilde{I}_\Theta$, we have a problem posed only in the 6 dimensions of $T^*(\mathrm{SO}(3))$, rather than a problem in the $6n$ dimensions it would take to represent the positions of all the point masses directly.

Now, before we talk about how we determine $\xi_H$ from $H$, let us discuss what properties such a determination should satisfy, speaking from physical intuition. By no means is the following argument a *deduction* of Hamiltonian mechanics. However, it should show at least that Hamiltonian mechanics is a reasonable idea for modeling closed mechanical systems.

First of all, as we are dealing with a *closed* system, we need for $H$ to be *conserved* along the motion of $\xi_H$. That is, if $\gamma: \mathbb{R} \to X$ is a solution to (6.1), then

$$\frac{\mathrm{d}}{\mathrm{d}t} H(\gamma(t)) = 0 \tag{6.2}$$

Using the chain rule, we get

$$0 = \frac{\mathrm{d}}{\mathrm{d}t} H(\gamma(t)) = \langle \mathrm{d}H(\gamma(t)), \gamma'(t)\rangle = \langle \mathrm{d}H(\gamma(t)), \xi_H(\gamma(t))\rangle$$

Or, put in a simpler manner,

$$\langle \mathrm{d}H, \xi_H\rangle = 0 \tag{6.3}$$

To derive our second assumption, note that through physical experiments we cannot determine the total energy of a system, only changes in that total energy. Thus if $H'(x) = H(x) + c$ for some constant $c$, we should have $\xi_H = \xi_{H'}$. We write this as

$$\xi_H = \xi_{H+c} \tag{6.4}$$

Finally, we assume that there is no "action at a distance". That is, if $H$ and $H'$ agree in some small neighborhood around $x$, then $\xi_H(x)$ should be equal to $\xi_{H'}(x)$; the evolution of a system should only depend on how the energy changes *locally*, and not on changes of energy that might be far away. Thus, for $U$ open and $x \in U$

$$H|_U = H'|_U \text{ implies that } \xi_H(x) = \xi_{H'}(x) \tag{6.5}$$

A simple determination of $\xi_H$ from $H$ that satisfies Equations 6.4 and 6.5 is for $\xi_H(x)$ to just depend linearly on $\mathrm{d}H(x)$, i.e.

$$\xi_H(x) = J(x)\,\mathrm{d}H(x) \tag{6.6}$$

Then, in order to satisfy Equation 6.3, we require $J(x): T_x^* X \to T_x X$ to be a *skew-symmetric* linear map, i.e.

$$\langle \varphi, J(x)\,\varphi \rangle = 0$$

for all $\varphi \in T_x^* X$. Finally, we also require that $J(x)$ depends smoothly on $x$. Now, a linear map $J(x): T_x^* X \to T_x X$ is the same thing as an element of $(T_x^* X)^* \otimes T_x X \cong T_x X \otimes T_x X$, and then skew-symmetry is expressed by $J(x)$ residing in the second exterior power of the tensor product, $\Lambda^2 T_x X$. Thus, $J$ can be expressed as a section of $\Lambda^2 TX$.

It turns out that many mechanical systems (including Examples 6.1, 6.2, 6.3, 6.4) can indeed be described by such a structure. Thus, we make the following definition.

DEFINITION 6.5. *An **almost Poisson manifold** is a manifold $X$ along with a smooth section $J$ of the bundle $\Lambda^2 TX$. This gives a Poisson structure (see Definition 4.25) to $T_x X$ for each $x$. We thus also say that $J$ is a **almost Poisson structure** on $X$.*

The name of the previous definition suggests that there is an another condition needed to make an almost Poisson manifold into a Poisson manifold. And this is indeed the case; Poisson manifolds also have have an integrability condition on $J$ which is necessary for the flow of $\xi_H$ to preserve $J$. As we will not be talking about integrability, we will not go into this in depth, though we will not discuss any examples that *aren't* integrable, so all of our examples will be of proper Poisson manifolds. For a more in-depth look at integrability and Poisson geometry, see Weinstein [66].

**Example 6.6.** In this example, we derive an appropriate Poisson structure on $T^*\mathbb{R}^3$ by using our physics intuition applied to Example 6.1. Recall that in Example 6.1, we described the state of a free particle with coordinates $(q, p)$ representing position and momentum respectively, and we had a Hamiltonian

$$H(q, p) = \frac{1}{2\,m}\|p\|^2 + U(q)$$

where

$$U(q) = -\int_a^b F(x) \cdot \gamma'(t)\,\mathrm{d}t$$

for some conservative force $F$, and any $\gamma: [a, b] \to \mathbb{R}^3$ such that $\gamma(a) = q_0$, $\gamma(b) = q$. Then we can derive $F$ from $U$ by

$$F = -\frac{\partial U}{\partial q}$$

In fact, we can derive $F$ from $H$ via

$$F = -\frac{\partial H}{\partial q}$$

because the kintetic energy does not depend on $q$.

Now, Newton's second law is commonly written as

$$F = m\,a = m\,\ddot{x}$$

However, using $p = m\,\dot{x}$, we can equivalently write it as

$$F = \dot{p}$$

This is the form that we want to use for the present derivation, because we can now write

$$\dot{p} = -\frac{\partial H}{\partial q}$$

We are halfway to a differential equation for the evolution of $(q, p)$; we now must derive $\dot{q}$. This is more straightforward, as $\dot{q}$ is velocity, and momentum is mass times velocity,

$$\dot{q} = \frac{1}{m}\,p$$

If we now take the partial derivative of $H$ with respect to $p$, we get precisely this.

$$\frac{\partial H}{\partial p} = \frac{1}{m} p$$

This seems perhaps a bit convoluted, but it turns out that the equation

$$\dot{q} = \frac{\partial H}{\partial p}$$

is in fact more fundamental than $\dot{q} = \frac{1}{m} p$; when relativity is brought into the picture $\dot{q} = \frac{1}{m} p$ is no longer true, but $\dot{q} = \frac{\partial H}{\partial q}$ remains true. In any case, this now allows us to write a set of equations for the evolution of $(q, p)$ in terms of partial derivatives of $H$

$$\begin{aligned} \dot{q} &= \frac{\partial H}{\partial p} \\ \dot{p} &= -\frac{\partial H}{\partial q} \end{aligned}$$

If we write

$$\mathrm{d}H = \begin{bmatrix} \frac{\partial H}{\partial q} \\ \frac{\partial H}{\partial q} \end{bmatrix}$$

and then let

$$J = \begin{bmatrix} 0 & I \\ -I & 0 \end{bmatrix}$$

where $I$ is the identity matrix, we can then write this equation as

$$\begin{bmatrix} \dot{q} \\ \dot{p} \end{bmatrix} = \begin{bmatrix} 0 & I \\ -I & 0 \end{bmatrix} \begin{bmatrix} \frac{\partial H}{\partial q} \\ \frac{\partial H}{\partial q} \end{bmatrix} = J\,\mathrm{d}H$$

This $J$ is then a Poisson structure on $X = T^*\mathbb{R}^3$.

When we write $J$ out as a matrix, the fact that it is an alternating map can be expressed by saying that it is **skew-symmetric**, that is $J_{ij} = -J_{ji}$.

Note that we have been rather cavalier here with regard to differential geometry, in that we have identified both $T_{(q,p)}X$ and $T^*_{(q,p)}X$ with $\mathbb{R}^6$ in order to write down $J$ as a matrix.

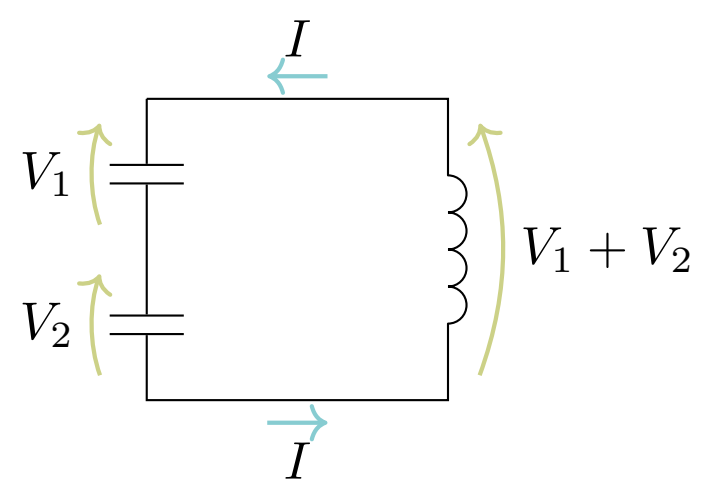

**Figure 6.4.** Currents and voltages for the circuit from Figure 6.2

**Example 6.7.** We now put a Poisson structure on the state space $\mathbb{R}^3$ of the circuit in Example 6.2. We want to make a Poisson structure that gives the right dynamics when applied to the Hamiltonian

$$H(q_1, q_2, \Phi) = \frac{1}{2C_1} q_1^2 + \frac{1}{2\,C_2} q_2^2 + \frac{1}{2L} \Phi^2$$

As all of the circuit components are connected in series, there is a single current $I$ that runs through all of them, and this current is related to the flux linkage in the inductor by

$$L\,I = \Phi_B$$

where $L$ is the inductance of the inductor. We can get this current out of the Hamiltonian via

$$I = \frac{\partial H}{\partial \Phi}$$

This current feeds into both of the capacitors, so

$$\dot{q}_1 = \dot{q}_2 = I = \frac{\partial H}{\partial \Phi}$$

To derive the dynamics of $\Phi$, recall that the voltage across a capacitor satisfies

$$C_i\, V_i = q_i$$

This can be written in terms of the Hamiltonian as

$$V_i = \frac{\partial H}{\partial q_i}$$

Then the change in the flux linkage is given by

$$\dot{\Phi} = -(V_1 + V_2) = -\left(\frac{\partial H}{\partial q_1} + \frac{\partial H}{\partial q_2}\right)$$

There is a negative sign there because of how we chose the directionality of voltage. In total, we have

$$\begin{aligned}
\dot{q}_1 &= \frac{\partial H}{\partial \Phi} \\
\dot{q}_2 &= \frac{\partial H}{\partial \Phi} \\
\dot{\Phi} &= -\left(\frac{\partial H}{\partial q_1} + \frac{\partial H}{\partial q_2}\right)
\end{aligned}$$

or equivalently,

$$\begin{bmatrix} \dot{q}_1 \\ \dot{q}_2 \\ \dot{\Phi} \end{bmatrix} = \begin{bmatrix} 0 & 0 & 1 \\ 0 & 0 & 1 \\ -1 & -1 & 0 \end{bmatrix} \begin{bmatrix} \frac{\partial H}{\partial q_1} \\ \frac{\partial H}{\partial q_2} \\ \frac{\partial H}{\partial \Phi} \end{bmatrix}$$

If we let $J$ be the matrix above, then we have a Poisson structure on the manifold $\mathbb{R}^3$. In later sections, we show how this Poisson structure can be deduced directly from the interconnection pattern (which is a 1-junction).

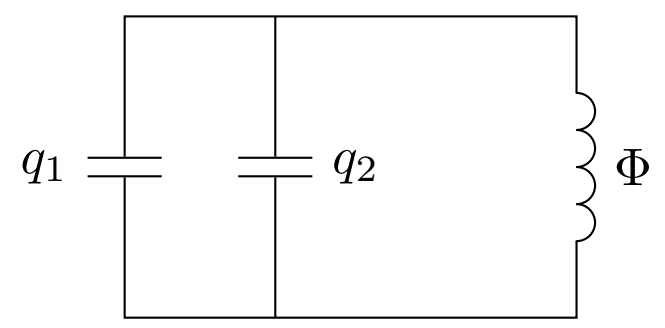

**Figure 6.5.** Two capacitors in parallel

In general, however, a circuit does *not* represent a Poisson manifold. For one, circuits often have resistors, which do not conserve energy. However, even a circuit without resistors are not in general represented by a Poisson manifold. The reason why can be illustrated by Figure 6.5. In that circuit, the *state* of the circuit is constrained by the fact that the voltages across each of those capacitors must be equal, so the charges on the capacitors must always satisfy

$$\frac{1}{C_1} q_1 = \frac{1}{C_2} q_2$$

Such a constraint of systems cannot be expressed within the formalism of Poisson manifolds.

But most importantly, Poisson manifolds cannot express open systems, i.e. systems where energy comes in and out. The generalization to port-Hamiltonian systems provides the flexibility to tackle these problems.

## 6.2. Port-Hamiltonian systems

Port-Hamiltonian systems add two innovations to Poisson manifolds. The first is that a port-Hamiltonian system has an *interface* by which it can transmit and receive energy. This takes the form of a bond space. Recall from section 4.2 that a bond space is a vector space $V$ along with a split quadratic form $\pi$ called the power form, and the architypical example of bond spaces is $\mathcal{E} \oplus \mathcal{F}$, where $\mathcal{E}$ is a vector space and $\mathcal{F} = \mathcal{E}^*$, along with the quadratic form $\pi(e + f) = \langle e, f \rangle$.

For instance, a capacitor is a port-Hamiltonian system with interface $\mathbb{B} = \mathbb{R} \oplus \mathbb{R}^*$. An element $(I, V) \in \mathbb{R} \oplus \mathbb{R}^*$ represents the current running through and voltage across the capacitor, respectively, and their product $I V$ represents the power running into the capacitor. Typically port-Hamiltonian systems also model dissiapative elements, such as resistors, however we have chosen to leave out dissiapative elements for simplicity.

The second innovation is that we drop the requirement that we can uniquely determine $\dot{x}$ from $\mathrm{d}H(x)$, as is true in Poisson systems. We replace this with a more general *relation* between $\dot{x}$, $\mathrm{d}H(x)$, and the interface variables.

DEFINITION 6.8. *If $B$ is a bond space, then* ***port-Hamiltonian system*** *on $B$ consists of*

- *A manifold $X$. This represents the* ***state*** *of the system, and is called the* ***state manifold****.*
- *A smooth function $H: X \to \mathbb{R}$. This represents the* ***energy*** *of the system, and is called the* ***Hamiltonian****.*
- *An element $D$ of $\mathsf{DiracRel}(X)(TX \oplus T^*X, \mathrm{Triv}_X(B))$. That is, for every $X$, $D_x$ is a Dirac relation between $T_xX \oplus T_x^*X$ and $B$. This represents the* ***connection*** *between the system and its interface, and is called the* ***connecting Dirac relation****.*

In order to understand port-Hamiltonian systems, we really must talk about their semantics in terms of *trajectories*, but for concreteness, we give an example before we get into that.

**Example 6.9.** A capacitor is a port-Hamiltonian system on interface $\mathbb{B} = \mathbb{R} \oplus \mathbb{R}^*$, with

- $X = \mathbb{R}$. The state of a capacitor is the amount of charge stored, so we use the letter $q$ to represent a point in $X$.
- $H(q) = \frac{1}{2C} q^2$. The energy of a capacitor with $q$ charge stored is $\frac{1}{2C} q^2$, where $C$ is the capacitance.
- $D_q \subset (TX \oplus T^*X) \oplus (\mathbb{R} \oplus \mathbb{R}^*)$ defined by

  $$D_q = \{(\xi, \varphi, V, I) \mid \xi = I, V = \varphi\}$$

  This connecting Dirac relation makes more sense when we talk about trajectories, but roughly speaking, this mean that the derivative of the charge stored is the current, and the voltage is the derivative of $H$ with respect to the charge stored.

We now discuss the “behavior” of port-Hamiltonian systems. Unlike with Poisson manifolds, a port-Hamiltonian system does not necessarily have a unique forwards evolution. This is because the evolution of a port-Hamiltonian system depends crucially on inputs and outputs: the values over time of the interface variables.

Definition 6.10. *Suppose that $(X, H, D)$ is a port-Hamiltonian system on interface $B$. Let $[a,b] \subset \mathbb{R}$ be an interval, and let $b\colon [a,b] \to B$ be a smooth path in $B$. A **trajectory** of $(X,H,D)$ **compatible with** $b$ consists of a map $\gamma\colon [a,b] \to X$ such that for all $t \in [a,b]$,*

$$(\gamma'(t), \mathrm{d}H(\gamma(t)), b(t)) \in D_{\gamma(t)}$$

In the equation for a trajectory, we relate $\gamma'(t)$, $\mathrm{d}H(\gamma(t))$ and $b(t)$. Often it is convenient to simply use $\gamma'(t)$ and $\mathrm{d}H(\gamma(t))$ in place of generic “dummy” variables when specifying the relation between $TX \oplus T^*X$ and $B$. For instance, we might describe a capacitor by the equations

$$\dot{q} = I$$

$$\frac{\partial H}{\partial q} = C q = V$$

instead of using variables $\xi$ and $\varphi$ that later get identified with $\dot{q}$ and $\frac{\partial H}{\partial q}$ as in Example 6.9. We call these the **kinetic equations**, and this should be understood as a notational convenience for conserving variable names; the formal structure remains unchanged.

The compatibility condition from Definition 6.10 is precisely what we need to prove the following theorem.

Theorem 6.11. *If $(X,H,D)$ is a port-Hamiltonian system on interface $B$, $b\colon [a,b] \to B$ is a smooth path through the interface, and $\gamma\colon [a,b] \to X$ is a compatible trajectory, then for all $t \in [a,b]$*

$$H(\gamma(t)) - H(\gamma(a)) = \int_a^t \pi(b(t))\, \mathrm{d}t$$

**Proof.** As a Dirac relation is power-preserving,

$$(\gamma'(t), \mathrm{d}H(\gamma(t)), b(t)) \in D_{\gamma(t)}$$

implies that

$$\pi(b(t)) = \langle \gamma'(t), \mathrm{d}H(\gamma(t)) \rangle$$

The result then follows from the fact that

$$\frac{\mathrm{d}}{\mathrm{d}t} H(\gamma(t)) = \langle \gamma'(t), \mathrm{d}H(\gamma(t)) \rangle$$

and the fundamental theorem of calculus. □

Thus, we can calculate the energy difference along a trajectory purely by looking at the power put into the interface.

**Example 6.12.** In this example, we construct a trajectory of the capacitor, as defined in Example 6.9. Let $[a,b]$ be an interval, let $I\colon [a,b] \to \mathbb{R}$ be any smooth function, and fix $q_0 \in \mathbb{R}$. Then define $q\colon [a,b] \to X = \mathbb{R}$ by

$$q(t) = q_0 + \int_a^t I(t)\, \mathrm{d}t$$

and $V\colon [a,b] \to \mathbb{R}$ by

$$V(t) = \frac{1}{C} q(t)$$

Then $V \oplus I\colon [a,b] \to \mathbb{R} \oplus \mathbb{R}^* = \mathbb{B}$ is a path through the bond space $\mathbb{B}$, and $q\colon [a,b] \to X$ is a compatible trajectory, as

$$\frac{\partial}{\partial q} H(q) = \frac{\partial}{\partial q} \frac{1}{2C} q^2 = \frac{1}{C} q$$

We can see that for an arbitrary path $b\colon [a,b]\to\mathbb{B}$, there is probably no compatible trajectory $q\colon [a,b]\to X$, and in general there exists at most one path. Thus, we can think of the port-Hamiltonian system as putting a *constraint* on its interface variables.

In this scenario, Theorem 6.11 states that

$$\frac{1}{2C}\,q(t)^2-\frac{1}{2\,C}\,q(a)^2=\int_a^t V(\tau)\,I(\tau)\,\mathrm{d}\tau$$

which is a well-known fact for a capacitor.

**Example 6.13.** An immovable wall is a very simple port-Hamiltonian system on $\mathbb{B}$. It has

- State space $X=\mathbb{R}^0$. The wall only has one state.
- Energy function $H=0$. The energy of that one state is 0.
- Connecting Dirac relation $\{(F,0)|F\in\mathbb{R}\}\subset\mathbb{B}\cong T\mathbb{R}^0\oplus T^*\mathbb{R}^0\oplus(\mathbb{B}\times\mathbb{R}^0)$. When the state space is trivial, the connecting Dirac relation is simply a Dirac structure on the interface.

A trajectory exists for a path $(F,v)\colon [a,b]\to\mathbb{R}$ if and only if $v(t)=0$ for all $t\in[a,b]$. This enforces the property that one can push and pull on a immovable wall as much as one likes, but it won't go anywhere. In this scenario, Theorem 6.11 states that pushing on an immovable wall does not impart any energy.

**Example 6.14.** If $(X,J)$ is a Poisson manifold, and $H\colon X\to\mathbb{R}$ is any smooth function, then define a Dirac relation between $TX\oplus T^*X$ and the empty interface $\mathbb{B}^0$ in the following way. First note that such a Dirac relation is simply a Dirac structure on $TX\oplus T^*X$. Then define a Dirac structure $D$ by

$$D_x=\mathrm{graph}(J(x))=\{(J(x)\,\varphi,\varphi)\mid\varphi\in T_x^*X\}\subset T_xX\oplus T_x^*X$$

A trajectory of $(X,H,D)$ compatible with the trivial path $b\colon [a,b]\to\mathbb{B}^0$ is then a map $\gamma\colon [a,b]\to X$ such that

$$\gamma'(t)=J(\gamma(t))\,\mathrm{d}H(\gamma(t))$$

This is a solution to the differential equation

$$\dot{x}=J(x)\,\mathrm{d}H(x)$$

which is Hamilton's equation. Thus we see that Hamiltonian mechanics is subsumed by port-Hamiltonian systems.

**Example 6.15.** In this example, we consider the example of a *forced pendulum*. This is a port-Hamiltonian system on interface $\mathbb{B}$ with

- State space $X=T^*S^1$, representing rotation $\theta$ and rotational momentum $L$
- Hamiltonian

  $$H(\theta,L)=m\,g\,l\sin\theta+\frac{1}{2\,I}L^2$$

  where $m$ is the weight at the end of the pendulum, $l$ is the length of the pendulum, and $I=m\,l^2$ is the moment of inertia of the pendulum.
- Dirac relation

  $$D_{(\theta,L)}=\{((\omega,\tau_{\mathrm{in}}+l\cos\theta\,F_{\mathrm{ap}}),(-\tau_{\mathrm{in}},\omega),(F_{\mathrm{ap}},l\cos\theta\,\omega))\mid\omega,F_{\mathrm{ap}},\tau_{\mathrm{in}}\in\mathbb{R}\}$$

  In this relation, $\omega$ is the angular velocity, $\tau_{\mathrm{in}}$ is the "internal torque" produced from the Hamiltonian (in this case produced by gravity), and $F_{\mathrm{ap}}$ is the linear applied force to the pendulum. Then the applied torque is

  $$\tau_{\mathrm{ap}}=l\cos\theta\,F_{\mathrm{ap}}$$

and the velocity at which it is applied is

$$v = l\cos\theta\,\omega$$

If we fix $\tau_{\text{ap}}\colon [a,b]\to\mathbb{R}$, the differential equations for a trajectory become

$$\begin{aligned}\dot{\theta} &= \frac{\partial}{\partial L}H = \frac{1}{I}L\\ \dot{L} &= -\frac{\partial}{\partial\theta}H + F_{\text{ap}}(t)\,l\cos\theta = -m\,g\,l\cos\theta + \tau_{\text{ap}}(t)\end{aligned}$$

If we let $v(t) = l\cos\theta\,\dot{\theta}(t)$, then a solution to these differential equations is a trajectory compatible with $(\tau_{\text{ap}},\omega)\colon [a,b]\to\mathbb{B}$.

One natural question to ask about port-Hamiltonian systems is what is the correct notion of morphism between them. To answer this, we must think about what one might want to accomplish with a morphism of port-Hamiltonian systems. One such goal could be "simplifying" a port-Hamiltonian system. That is, suppose that a port-Hamiltonian system consisted of two parts that did not interact with each other, and that only one of the parts interacted with the interface. Then from the perspective of the interface, we could drop the non-interacting part.

There likely could be other definitions of morphisms between port-Hamiltonian systems, but the following notion is at least useful in the above regard, and we do not consider other types of morphism. To our knowledge, the following definition and the propositions depending on it are original to this thesis, however it is similar in spirit to work from Barbero-Linam [67].

Definition 6.16. *If $(X,H,D)$ and $(X',H',D')$ are port-Hamiltonian systems on the same interface $B$, then a **forwards morphism** between them is a smooth map $f\colon X\to X'$ such that*

1. *$H = H'\circ f$*
2. *For every $x\in X$, $\xi\in T_x X$, $\varphi\in T^*_{f(x)}X'$, $b\in B$,*

$$(\xi, f^*\varphi, b)\in D_x \Rightarrow (f_*\xi, \varphi, b)\in D'_{f(x)} \tag{6.7}$$

The following proposition shows that we can "push forwards" trajectories of port-Hamiltonian systems along theses morphisms.

Proposition 6.17. *Suppose that $(X,H,D)$ and $(X',H',D')$ are port-Hamiltonian systems on the same interface $B$, and $f\colon X\to X'$ is a forwards morphism between them. Then fix a map $b$: $[a,b]\to B$, and let $\gamma\colon[a,b]\to X$ be a trajectory of $(X,H,D)$ compatible with $b$. Then $f\circ\gamma\colon[a,b]\to X'$ is a trajectory of $(X',H',D')$ compatible with $b$.*

**Proof.** We must show that for every $t\in[a,b]$,

$$((f\circ\gamma)'(t), \mathrm{d}H'(f(\gamma(t))), b(t))\in D'_{f(\gamma(t))}$$

By the assumption that $\gamma$ is a trajectory of $(X,H,D)$ compatible with $b$, we have

$$(\gamma'(t), \mathrm{d}H(\gamma(t)), b(t))\in D_{\gamma(t)}$$

Now, as $H = H'\circ f$,

$$\mathrm{d}H(\gamma(t)) = \mathrm{d}(H'\circ f)(\gamma(t)) = f^*(\mathrm{d}H'(f(\gamma(t))))$$

Thus, we can rewrite (6.17) as

$$(\gamma'(t), f^*(\mathrm{d}H'(f(\gamma(t)))), b(t))\in D'_{f(\gamma(t))}$$

If we let $\varphi = \mathrm{d}H'(f(\gamma(t)))$, $\xi = \gamma'(t)$, and recall that $(f\circ\gamma)'(t) = f_*(\gamma'(t))$, then (6.7) implies that

$$((f\circ\gamma)'(t), \mathrm{d}H'(f(\gamma(t))), b(t))\in D'_{f(\gamma(t))}$$

as required. □

We would like to now form a category of port-Hamiltonian systems. To do this, we must show that forwards morphisms of port-Hamiltonian systems compose.

PROPOSITION 6.18. *Suppose that $(X, H, D)$, $(X', H', D')$ and $(X'', H'', D'')$ are all port-Hamiltonian systems on the same interface $B$, and that $f: X \to X'$ and $g: X' \to X''$ are forwards morphisms between them. Then $(g \circ f): X \to X''$ is also a forwards morphism.*

**Proof.** Clearly property 1 holds. Now, suppose that

$$(v, (g \circ f)^* \varphi, b) \in D_x$$

Then by property 2 for $f$,

$$(f_* v, g^* \varphi, b) \in D'_{f(x)}$$

and then by property 2 for $g$,

$$(g_* f_* v, \varphi, b) \in D''_{g(f(x))}$$

Note that we have implicitly used $(g \circ f)_* = g_* \circ f_*$ and $(g \circ f)^* = f^* \circ g^*$. □

DEFINITION 6.19. *For a bond space $B$, let $\mathrm{pH}(B)$ be the category where the objects are port-Hamiltonian systems on $B$ and the morphisms are forwards morphisms.*

We can now express the content of Proposition 6.17 as defining a functor.

DEFINITION 6.20. *Let $B$ be a bond space, and fix a smooth function $b: [a, b] \to B$. Then define a functor $\mathrm{Traj}_b: \mathrm{pH}(V) \to \mathsf{Set}$ by*

$$\begin{aligned} \mathrm{Traj}_b(X, H, D) &= \{\gamma: [a, b] \to X \mid \gamma \text{ is compatible with } b\} \\ \mathrm{Traj}_b(f: (X, H, D) \to (X', H', D')) &= \gamma \mapsto f \circ \gamma \end{aligned}$$

# 7. Operads

## 7.1. Operads

Composing functions is quite straightforward. Functions have a well-defined “input” and “output”, and if the output of one function matches up with the input of another, they can be composed. We can represent the composition of functions using a category.

However, composing physical systems is not so simple. A physical system does not necessarily have an “input” and an “output”. Moreover, we don’t necessarily want to only think about composing two physical systems at a time: we might want to compose many! So then what is the mathematical structure by which we can represent composition of physical systems?

One answer lies in *choosing* parts of the system to consider as the input and output, and then composing systems as morphisms in a category. Then various structures in the category allow one to “twist around” inputs to outputs and vice-versa. This idea has been formalized with decorated cospans and structured cospans; see Fong [68], Baez and Courser [43], and Baez, Courser, and Vasilakopolou [69].

In this thesis, we take a different approach to this problem, which makes the following tradeoff. We do not have to choose inputs and outputs for systems. However, we do have to choose how we compose systems; there is no longer any sort of canonical composition. This approach comes from the theory of *operads* and *operad algebras*. To understand how this works, we make the following analogy to group actions.

DEFINITION 7.1. *If $G$ is a group, and $X$ is a set, then a group action of $G$ on $X$ is a group homomorphism from $G$ into the group of bijections from $X$ to $X$.*

One can think of a group as a description of abstract symmetries. For instance, the group $\mathbb{Z}/7$ describes sevenfold symmetry. And then a group action provides the thing that actually has that symmetry, along with how that symmetry works.

Analogously, operads provide a description of "abstract composition operations." Then an operad algebra provides something for those composition operations to act on, and a description of how those composition operations work. Without further vagueries, we plunge into the formal definition.[7.1]

Definition 7.2. *An **operad** $\mathcal{O}$ consists of*

- *A collection of **types** $\mathcal{O}_0$*
- *For every $X_1, \ldots, X_n, Y \in \mathcal{O}_0$, a collection of **operations** $\mathcal{O}(X_1, \ldots, X_n; Y)$*

*along with*

- *For every type $X$, an operation $1_X \in \mathcal{O}(X; X)$ called **the identity operation***
- *For every tuple of types $X_1, \ldots, X_n$, $Y_1, \ldots, Y_m$, $Z$, $i \in \{1, \ldots, m\}$, a **composition map***

$$\circ_i : \mathcal{O}(X_1, \ldots, X_n; Y_i) \times \mathcal{O}(Y_1, \ldots, Y_m; Z) \to \mathcal{O}(Y_1, \ldots, Y_{i-1}, X_1, \ldots, X_n, Y_{i+1}, \ldots, Y_m; Z)$$

- *For every tuple of types $X_1, \ldots, X_n$, $Y$, and every permutation $\sigma \in S_n$, a **symmetry map***

$$\sigma^* : \mathcal{O}(X_1, \ldots, X_n; Y) \to \mathcal{O}(X_{\sigma(1)}, \ldots, X_{\sigma(n)}; Y)$$

*such that several laws hold. It is quite an ordeal to write down these laws, so we instead explain how they are derived using some graphical principles, and refer the reader to Leinster [70] for the details.*

There are multiple different ways of drawing operations in an operad, but we start with just one for now. We call this "tree notation", and in tree notation an operation $f \in \mathcal{O}(X_1, X_2, X_3; Y)$ would look something like this:

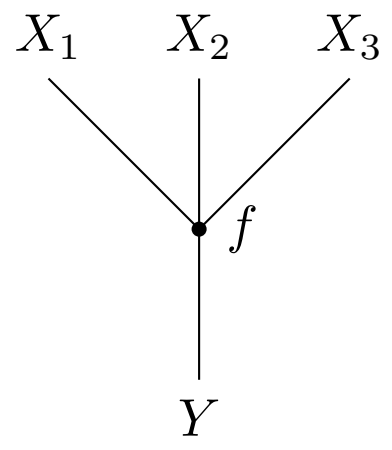

The reader might be thinking that this looks somewhat similar to string diagram notation, and there is in fact a connection here that we get into soon. Composition of multiple operations is performed simply by stacking multiple trees on top of each other, so that if $f \in \mathcal{O}(X_1, X_2; Y_2)$ and $g \in \mathcal{O}(Y_1, Y_2, Y_3; Z)$, their composition $f \circ_2 g$ looks like the following:

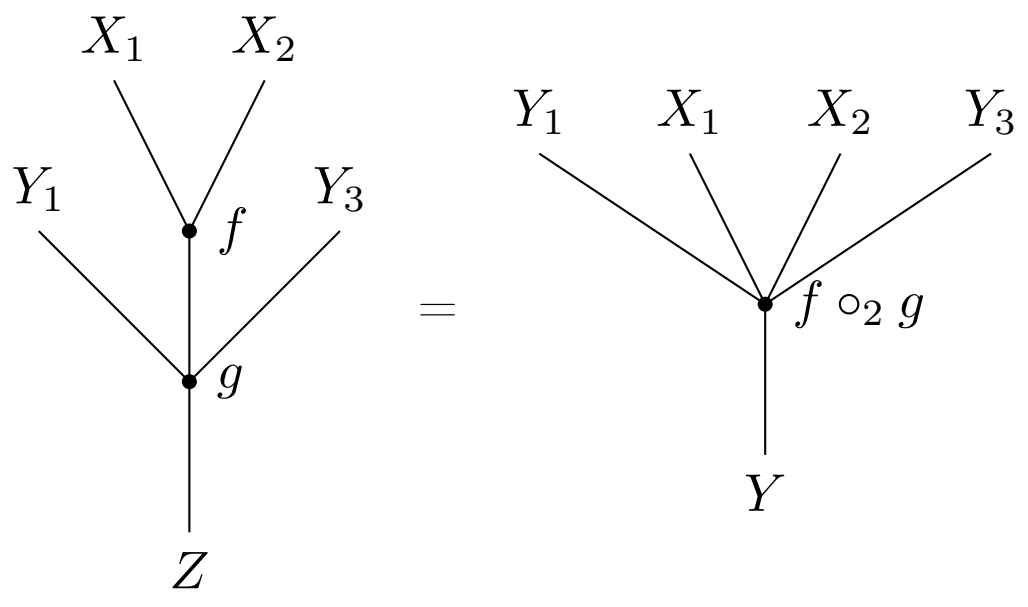

Associativity of composition means that it does not matter in what order we stack trees. This is implied naturally by the tree syntax. The following tree could either represent $(h \circ_2 f) \circ_2 g$, or $h \circ_3 (f \circ_2 g)$; it is a law that both of these must be equal.

7.1. Note that when we say "operad" we always mean the multi-colored version of operads, also called symmetric multicategories; if the reader has not encountered operads before they can safely disregard this footnote.

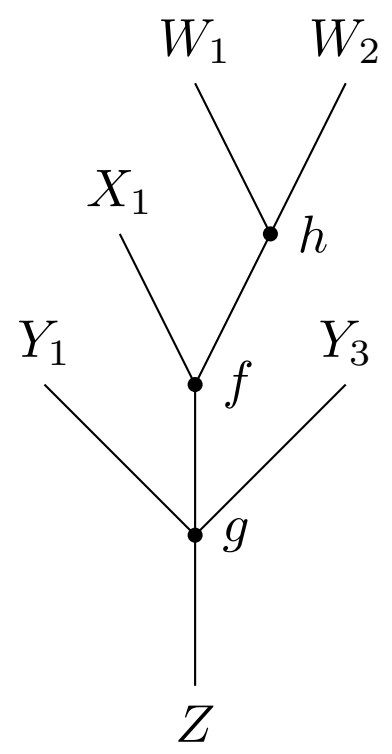

The composite is also independent of whether we stack trees on top of each other, or compose "in parallel", i.e. it does not matter which of $f$ or $f'$ we stack first in the following tree

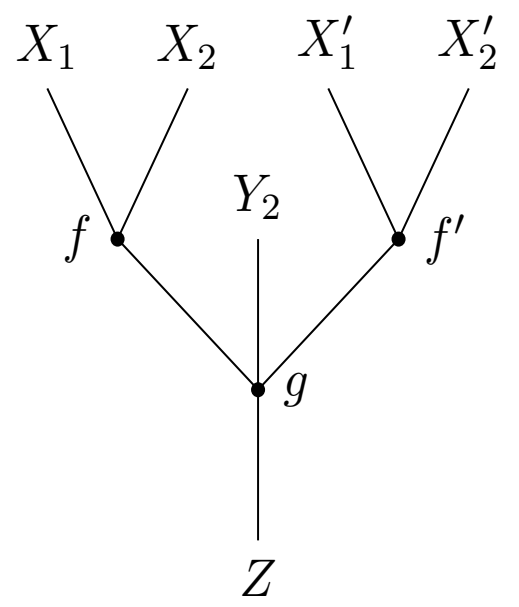

Just as in string diagrams, we draw the identity as a plain wire, that is

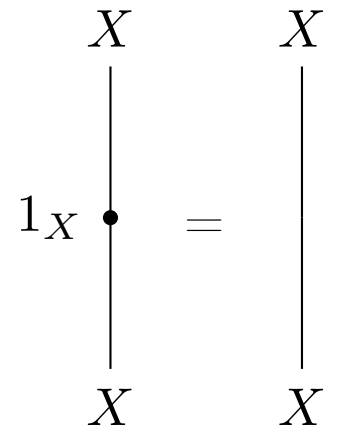

Composing with the identity leaves whatever you had unchanged.

There are in fact two different ways of defining composition in an operad, which are equivalent. The first is the way that we have already stated, where we compose two operations at a time. This looks like a function

$$\circ_i\colon \mathcal{O}(X_1,\ldots,X_n;Y_i)\times\mathcal{O}(Y_1,\ldots,Y_m;Z)\to\mathcal{O}(Y_1,\ldots,Y_{i-1},X_1,\ldots,X_n,Y_{i+1},\ldots,Y_m;Z)$$

The second way is that we compose an operation $g\in\mathcal{O}(Y_1,\ldots,Y_n;Z)$ with $n$ operations simultaniously, i.e. with $f_1\in\mathcal{O}(X_{1,1},\ldots,X_{1,k_1};Y_1),\ldots,f_n\in\mathcal{O}(X_{n,1},\ldots,X_{n,k_n};Y_n)$. This looks like a function

$$\circ\colon \mathcal{O}(Y_1,\ldots,Y_n;Z)\times\mathcal{O}(X_{1,1},\ldots,X_{1,k_1};Y_1)\times\cdots\times\mathcal{O}(X_{n,1},\ldots,X_{n,k_n};Y_n)\to\mathcal{O}(X_{1,1},\ldots,X_{n,k_n};Z)$$

We write the composition of $f_1,\ldots,f_n$ with $g$ as $g\circ(f_1,\ldots,f_n)$.

These two ways of defining composition are equivalent, because we can define the second composition using the first simply by applying the first composition multiple times, and we can define the first composition using the second by putting in identities in all of the slots except for the $i$th slot.

The first type of composition is more flexible in practice, because it is often the case that we want to only compose two operations. However, the second type of composition can often be easier to work with mathematically. We use both types of composition.

Finally, also have some laws for the symmetry map. Specifically, for $\sigma,\tau\in S_n$, $(\sigma\circ\tau)^*=\sigma^*\circ\tau^*$ (consequently $\sigma^*$ is a bijection). Moreover for $g\in\mathcal{O}(Y_1,\ldots,Y_n;Z)$, and $f_i\in\mathcal{O}(X_{i,1},\ldots,X_{i,k_i};Y_i)$

$$(\sigma^*g)\circ(f_1,\ldots,f_n)=\tau^*(g\circ(f_1,\ldots,f_n))$$

where $\tau \in S_{\sum_{i=1}^n k_i}$ is the permutation that swaps the $k_i$-sized blocks in $\{1, \ldots, \sum_{i=1}^n k_i\}$ according to $\sigma$. Thus, the symmetry map respects composition.

Definition 7.3. *If an operad $\mathcal{O}$ has only one type, we call it a* ***single-colored operad****. In this case, we refer to the set of n-ary operations with $\mathcal{O}_n$. That is, if the single type is $X$, then*

$$\mathcal{O}_n = \mathcal{O}(\underbrace{X, \ldots, X}_{n}; X)$$

**Example 7.4.** The **operad of convex combinations**, which we write as $\mathcal{CC}$, is a single-colored operad with

$$\mathcal{CC}_n = \left\{ (\lambda_1, \ldots, \lambda_n) \mid \lambda_i \in [0, 1], \sum_{i=1}^n \lambda_i = 1 \right\}$$

Composition is defined in the following way. Suppose that $(\lambda_1, \ldots, \lambda_n) \in \mathcal{CC}_n$, and that $(\gamma_{1,1}, \ldots, \gamma_{1,k_1}) \in \mathcal{CC}_{k_1}, \ldots, (\gamma_{n,1}, \ldots, \gamma_{n,k_n}) \in \mathcal{CC}_{k_n}$. Then let

$$(\lambda_1, \ldots, \lambda_n) \circ ((\gamma_{1,1}, \ldots, \gamma_{1,k_1}), \ldots, (\gamma_{n,1}, \ldots, \gamma_{n,k_n})) = (\lambda_1 \gamma_{1,1}, \ldots, \lambda_n \gamma_{n,k_n}) \in \mathcal{CC}_{\sum_{i=1}^n k_i}$$

We think of an operation in this operad as being a recipe for how to mix $n$ ingredients. That is, we take $\lambda_1$ of the first ingredient, $\lambda_2$ of the second ingredient, and so on. Then we can think of composition saying that we are going to take $\lambda_1$ of the result of the recipe $(\gamma_{1,1}, \ldots, \gamma_{1,k_1})$, $\lambda_2$ of the result of the recipe $(\gamma_{2,1}, \ldots, \gamma_{1,k_2})$, and this is the same as taking $\lambda_1 \gamma_{1,1}$ of the first ingredient of the first recipe, $\lambda_2 \gamma_{2,1}$ of the first ingredient of the second recipe, and so on.

To prove directly that something is an operad can be a somewhat tedious affair. Fortunately, we have a construction that does most of the heavy lifting for us.

Proposition 7.5. *If $(\mathsf{C}, \otimes, I)$ is a symmetric monoidal category, then we can build an operad which we call $\mathrm{Op}(\mathsf{C})$ in the following way.*

- *Let $\mathrm{Op}(\mathsf{C})_0 = \mathsf{C}_0$, so that a type in $\mathrm{Op}(\mathsf{C})$ is an object in $\mathsf{C}$*
- *For $X_1, \ldots, X_n, Y \in \mathrm{Op}(\mathsf{C})_0$, let*

$$\mathrm{Op}(\mathsf{C})(X_1, \ldots, X_n; Y) = \mathrm{Hom}_{\mathsf{C}}(X_1 \otimes \cdots \otimes X_n, Y)$$

*If $g \in \mathrm{Op}(\mathsf{C})(Y_1, \ldots, Y_n; Z)$ and $f_i \in \mathrm{Op}(\mathsf{C})(X_{i,1}, \ldots, X_{i,k_i}; Y_i)$, then we define*

$$g \circ (f_1, \ldots, f_n) = g \circ (f_1 \otimes \cdots \otimes f_n)$$

*The identity operations are given by the identity morphisms. Finally, the action of $\sigma \in S_n$ on operations is given by composing with the braiding.*

**Proof.** See Leinster [70, Chapter 2]. □

## 7.2. Finite cospans and wiring diagrams

In this section, we give an extended discussion on the operad of undirected wiring diagrams. Undirected wiring diagrams are widely used to describe the interconnection of systems, and are a very well-studied operad. We do not use them directly in this paper, but they have indirectly influenced our approach to a great extent. The reader is encouraged to peruse Spivak [10], Yau [71], Vagner [72], Fong [12], and Libkind [9] for more information on the subject.

Consider the following symmetric monoidal category, which we call $\mathsf{FinCospan}$. The construction of this category is dual to the construction of span categories, as discussed in Definition 2.34. The objects of this category are finite sets, and a morphism between finite sets $A$ and $B$ is an equivalence class of **cospans**, which are diagrams of the form

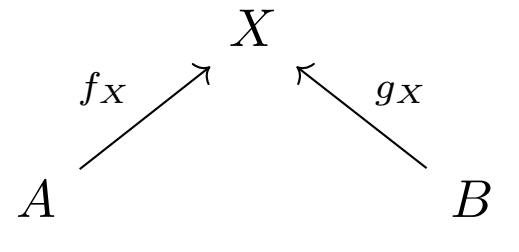

where $X$ is a finite set. Two cospans are equivalent if there exists an isomorphism $\phi$ as in the following commutative diagram.

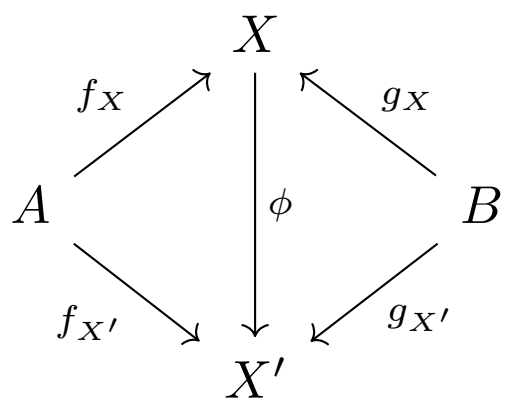

We visualize a cospan of finite sets in Figure 7.1.

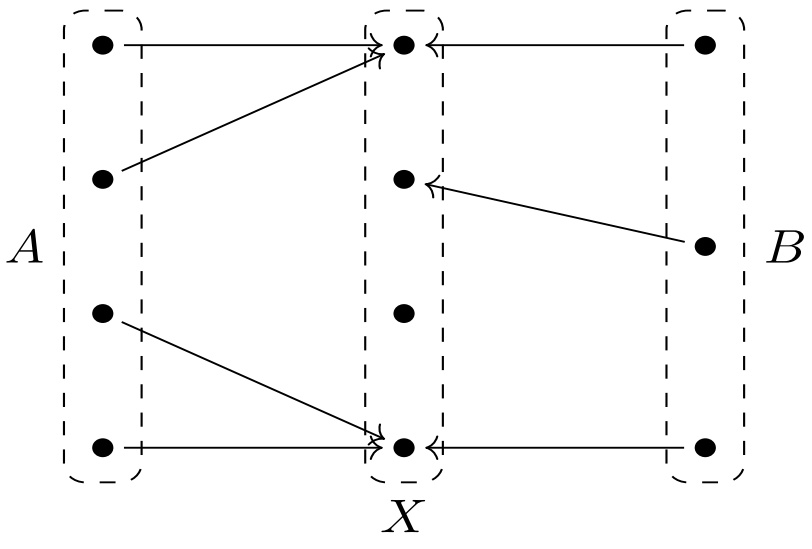

**Figure 7.1.** A morphism of FinCospan from $A$ to $B$

We think of this as providing a description of connections. That is, nodes on the outside that map to the same inner node are "connected." Composition of cospans is done via pushout.

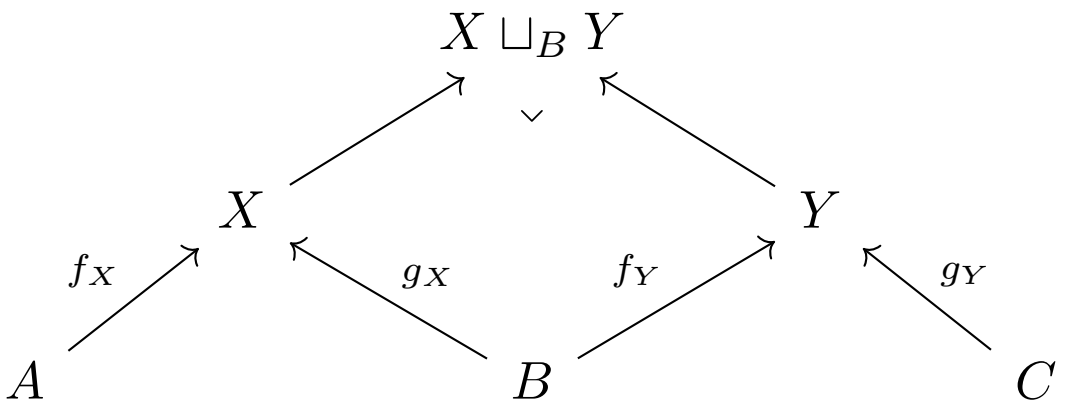

This definition for composition captures the notion that if $a \in A$ is connected to $b \in B$, and $b$ is conneted to $c \in C$, then $a$ should be connected to $c$ when the two connection schemes are merged. We draw this as in Figure 7.2.

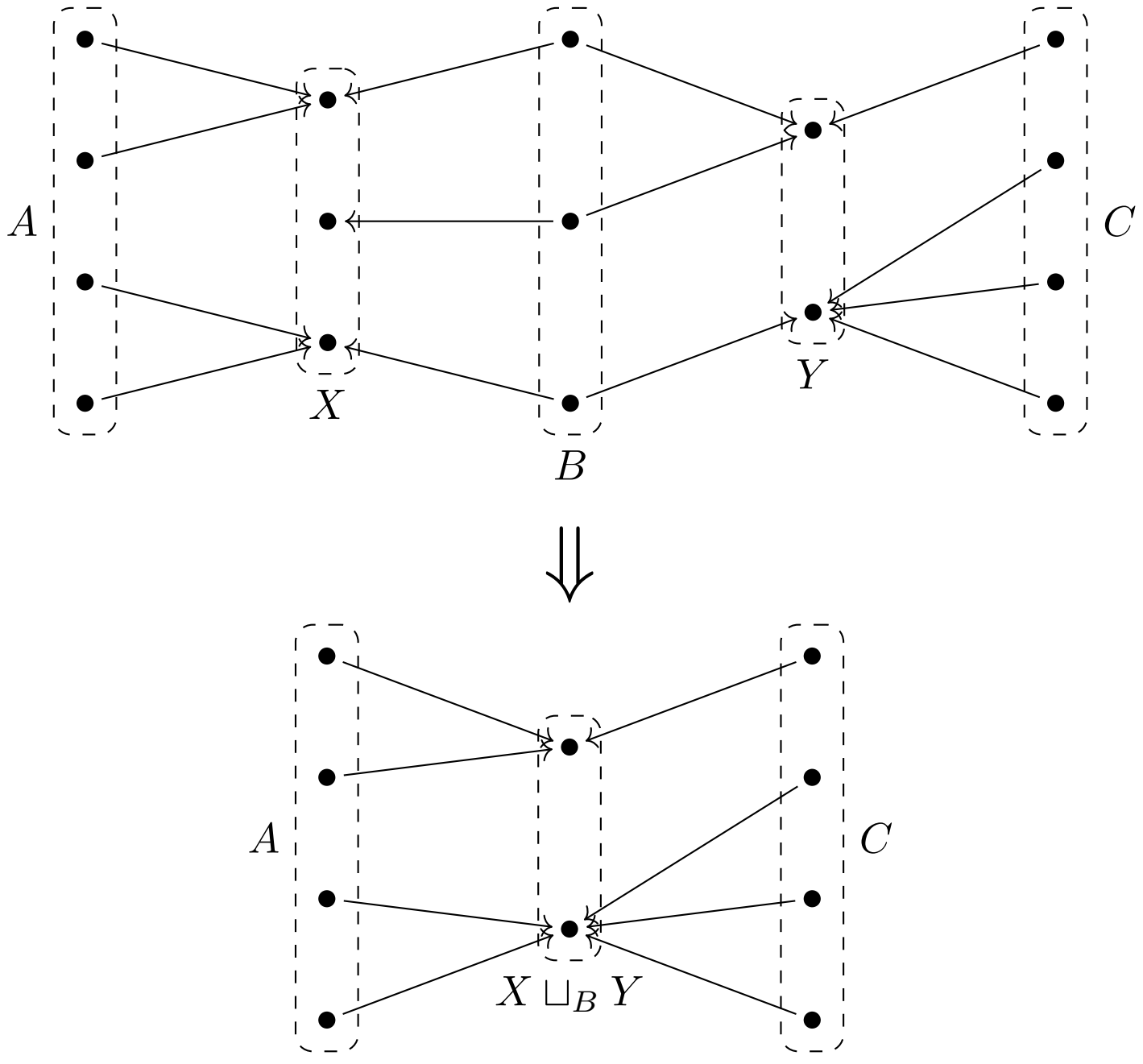

**Figure 7.2.** Composition of finite cospans

In analogy to Definition 2.34, we need to consider equivalence classes of cospans in order for composition to be well-defined, as pushout is defined only up to natural isomorphism.

We then put a monoidal structure on FinCospan via disjoint union. Given two cospans

$$A_1 \xrightarrow{f_{X_1}} X_1 \xleftarrow{g_{X_1}} B_1, \quad A_2 \xrightarrow{f_{X_2}} X_2 \xleftarrow{g_{X_2}} B_2$$

their monoidal product is

$$A_1 + A_2 \xrightarrow{f_{X_1}+f_{X_2}} X_1 + X_2 \xleftarrow{g_{X_1}+g_{X_2}} B_1 + B_2$$

which we visualize as in Figure 7.3.

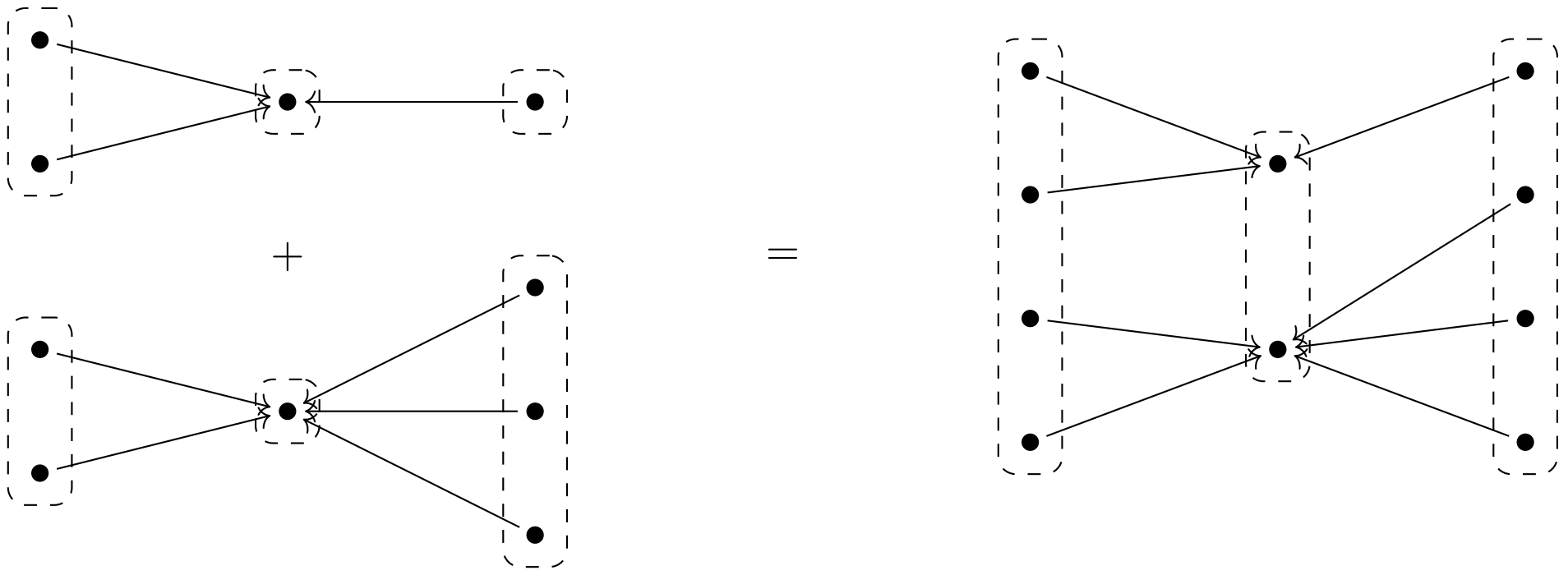

**Figure 7.3.** Monoidal composition of finite cospans

All of the required natural isomorphisms for this to be a symmetric monoidal product come from the fact that the category FinSet of finite sets and functions between them is a wide subcategory of FinCospan via the functor that sends a map $f: A \to B$ to the cospan

$$A \xrightarrow{f} B \xleftarrow{1_B} B$$

and this monoidal product is the cocartesian monoidal product for FinSet.

Now, we have a symmetric monoidal category $(\mathsf{FinCospan}, +, \emptyset)$, and we can build an operad out of it, Op(FinCospan). We call this operad $\mathcal{UWD}$, which stands for **undirected wiring diagrams**. We think of an operation in $\mathcal{UWD}$ as providing a description of how several boxes can be wired together, in an undirected fashion, and also exposing some "external ports." This corresponds to the visualization in Figure 7.4

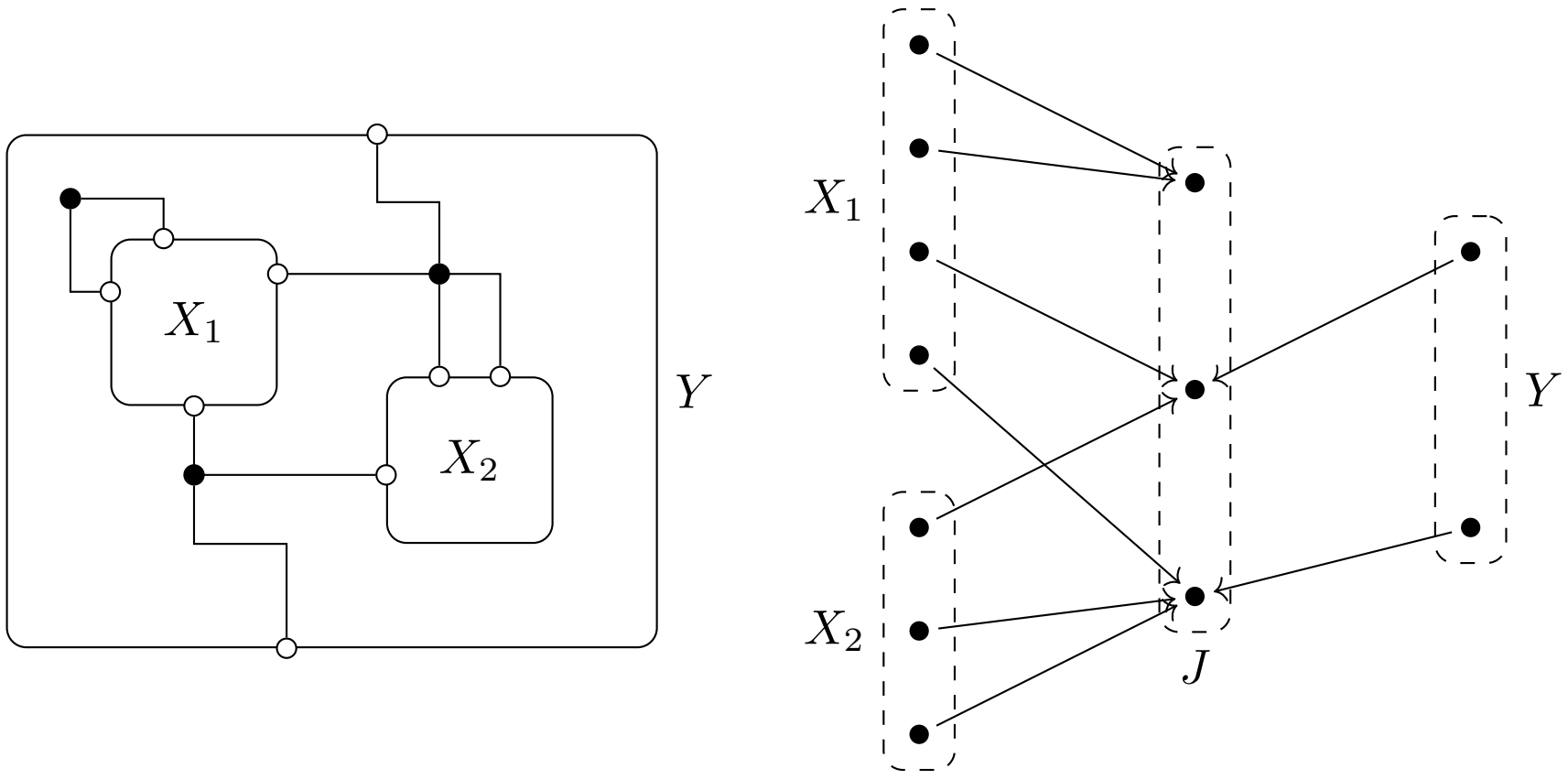

**Figure 7.4.** An undirected wiring diagram, and its underlying cospan $X_1 + X_2 \to J \leftarrow Y$ of finite sets

In the diagrams for undirected wiring diagrams, composition is drawn by nesting, as pictured in Figure 7.5.

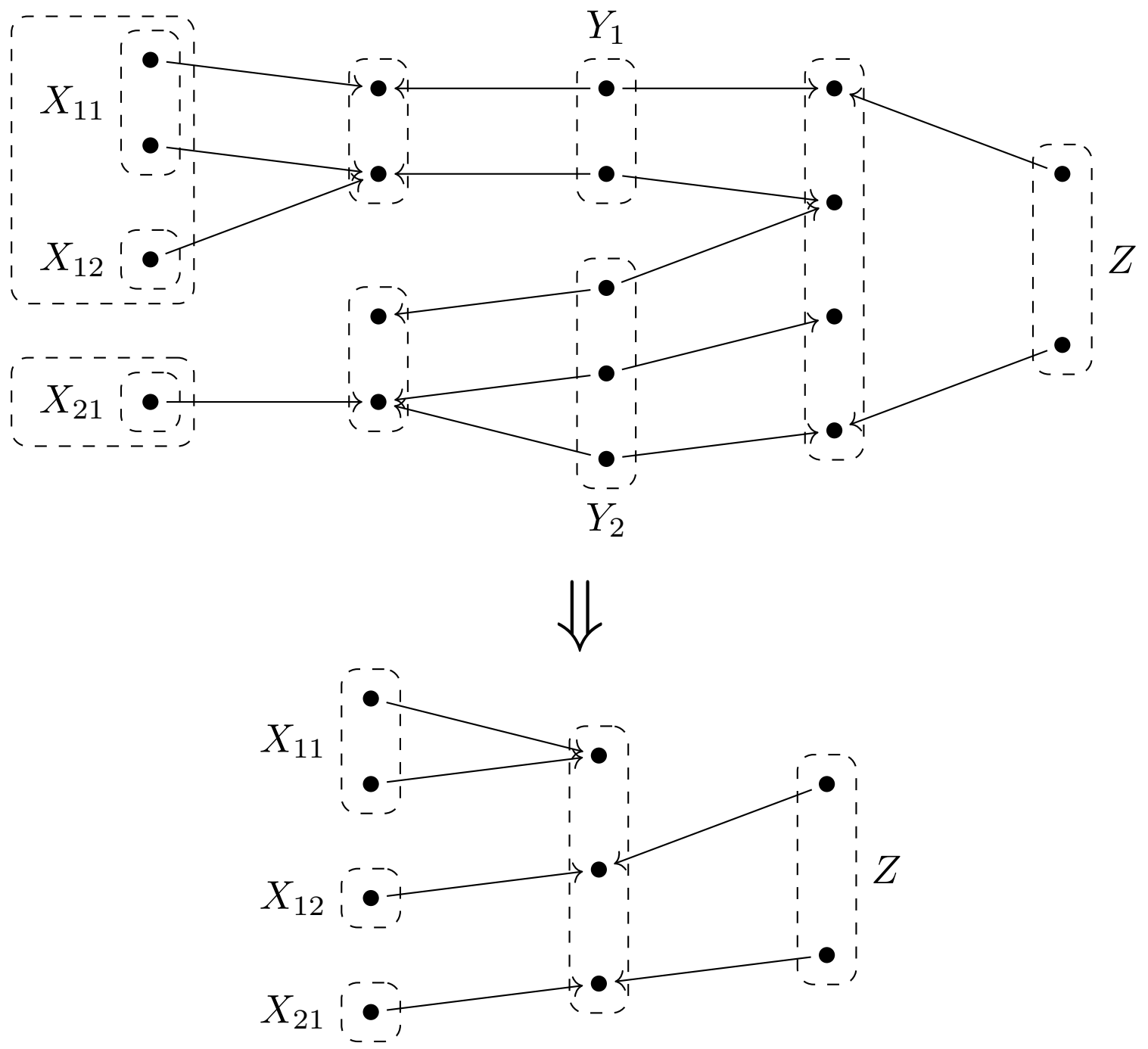

(a) Literal cospan diagram

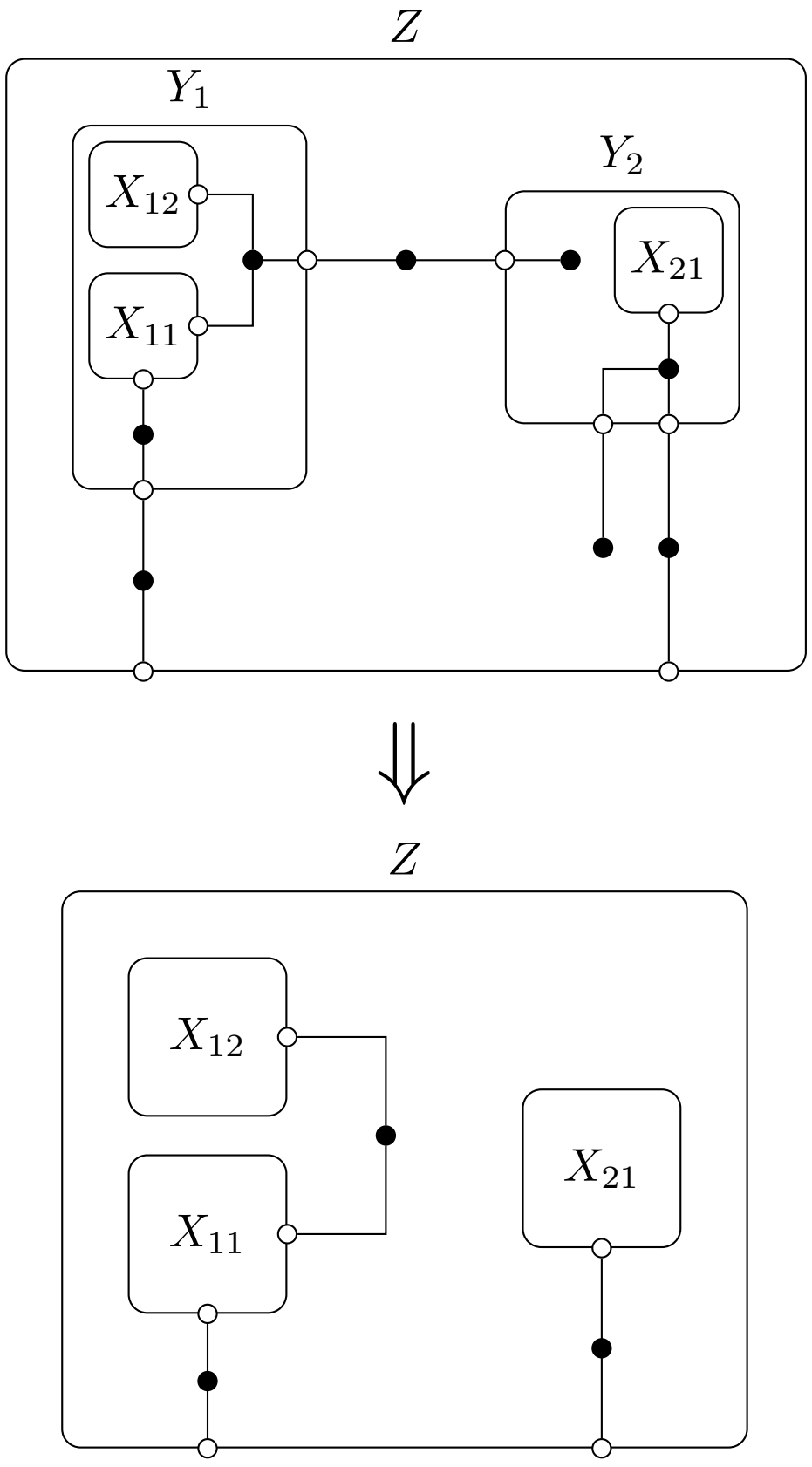

(b) Undirected wiring diagram

**Figure 7.5.** Two ways of drawing how to compose cospans.

## 7.3. Dirac relations as an operad

The main operad we use in this paper is built out of the symmetric monoidal category DiracRel (see section 4.5 for exposition of DiracRel).

Definition 7.6. *We define $\mathcal{DR}$ to be the operad* $\mathrm{Op}(\mathsf{DiracRel})$.

We picture operations in this operad using Dirac diagrams. We explain how this works using a sequence of diagrams. Note that there is not a one-to-one correspondence between Dirac diagrams and Dirac relations; the same Dirac relation could be pictured in several different ways.

**Example 7.7.** The simplest Dirac relation is the identity on the bond space $\mathbb{B} = \mathbb{R} \oplus \mathbb{R}^*$. This is pictured in Figure 7.6 (a). This is an operation that "takes in" something of type $\mathbb{B}$ in the smaller empty box, and produces the exact same thing. Figure 7.6 (b) pictures the identity on the bond space $\bar{\mathbb{B}}$, which has the same underlying space as $\mathbb{B}$ but has the opposite power structure. Unless otherwise specified, bonds in a Dirac diagram represent the bond space $\mathbb{B}$ or $\bar{\mathbb{B}}$, with the choice of which given by the directionality.

Note that we use a junction on the border of the outer box to represent the "external" ports. This allows us, when we compose relations by nesting as in Figure 7.5, to simply erase the box and be left with a valid Dirac diagram. We give an example of this later.

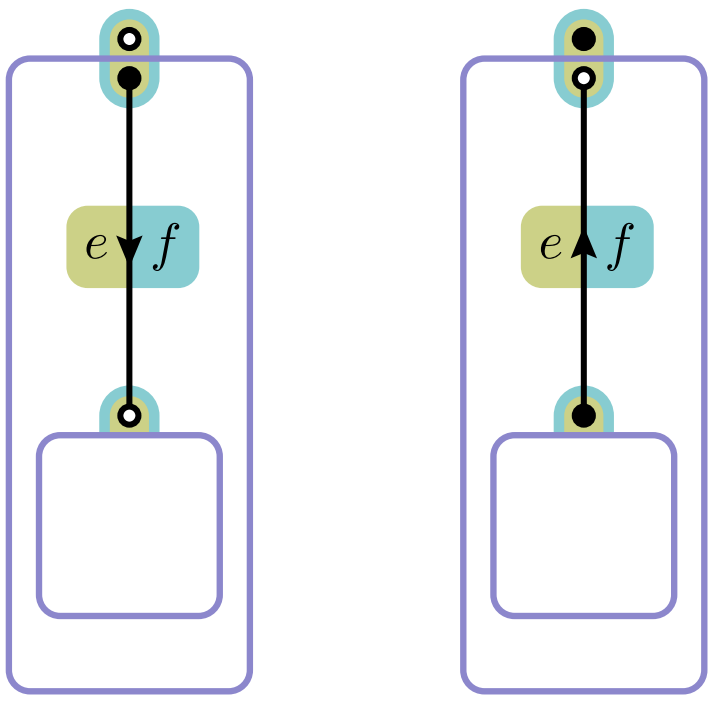

(a) $\mathrm{id}_{\mathbb{B}}$ (b) $\mathrm{id}_{\bar{\mathbb{B}}}$

**Figure 7.6.** The identity Dirac relation as an operation in $\mathrm{Op}(\mathsf{DiracRel})$, pictured in Dirac diagram style

**Example 7.8.** Another simple Dirac relation is $\cap \colon \bar{\mathbb{B}} \oplus \mathbb{B} \twoheadrightarrow 0$ (from Proposition 4.24) pictured in Figure 7.7. As an operation in $\mathcal{DR}$, this takes two things of type $\bar{\mathbb{B}}$ and $\mathbb{B}$ respectively, and produces a closed system

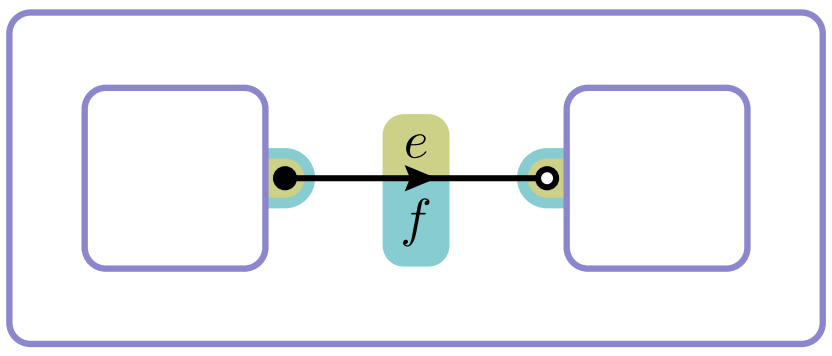

**Figure 7.7.** Creating a closed system by connecting two boxes.

**Example 7.9.** There are two ways of connecting two things of type $\mathbb{B}$. One is by inverting efforts,

and one is by inverting flows. Both are pictured in Figure 7.8. Recall from Chapter 4 that a junction works by *matching* the inner color, and *summing* the outer color. So in Figure 7.8 (a), the efforts are matched, and the flows are summed. Conversely, in Figure 7.8 (b) the flows are matched and the efforts are summed.

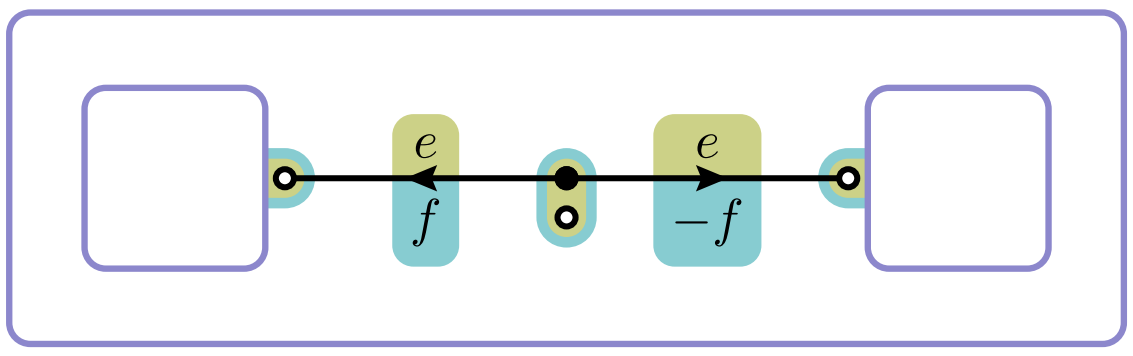

(a)

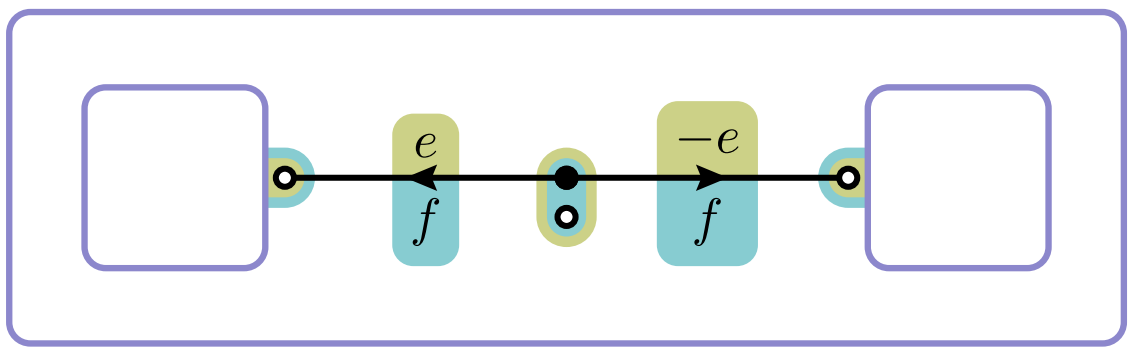

(b)

**Figure 7.8.** Two ways of connecting two things of type $\mathbb{B}$

**Example 7.10.** In Figure 7.9, we see how the composition of Dirac relations is pictured with the Dirac diagram syntax for Op(DiracRel). We visualize the setup for composition via nesting Dirac diagrams. Then to compose, we simply erase the intervening boxes. Optionally, we can then simplify the diagram, which simply means writing down the same relation in a different way.

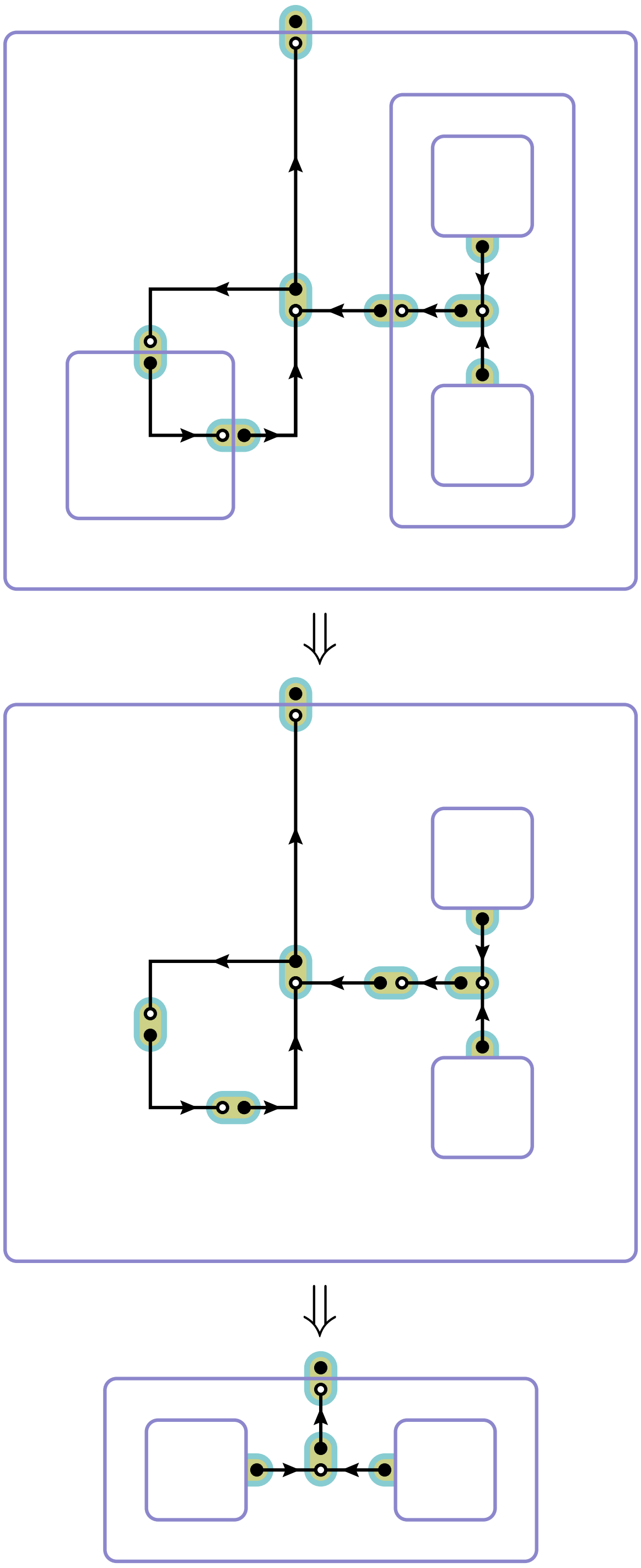

**Figure 7.9.** Visualizing composition in Op(DiracRel)

## 7.4. Operad algebras

If an operad describes ways of composing things, an operad algebra describes a choice of what those things actually are, and how the ways work to compose them. Just as a group action can be seen as a homomorphism from a group $G$ to the endomorphism group of a set $X$, we can also see an operad algebra as a certain homomorphism. To do this, we must define what a morphism of operads is.

Definition 7.11. *If $\mathcal{O}$ and $\mathcal{O}'$ are operads, a morphism $F$ between them consists of*

- *A type $F(X) \in \mathcal{O}'_0$ for each type $X \in \mathcal{O}_0$*

- *An operation* $F(f) \in \mathcal{O}'(F(X_1), \ldots, F(X_n); F(Y))$ *for every* $f \in \mathcal{O}(X_1, \ldots, X_n; Y)$

*such that*

- *if* $f, g_1, \ldots, g_n$ *are composable, then*

$$F(f \circ (g_1, \ldots, g_n)) = F(f) \circ (F(g_1), \ldots, F(g_n))$$

- *for* $\mathrm{id}_X \in \mathcal{O}(X; X)$*,*

$$F(\mathrm{id}_X) = \mathrm{id}_{F(X)}$$

- *for* $f \in \mathcal{O}(X_1, \ldots, X_n; Y)$ *and* $\sigma \in S_n$*,*

$$F(\sigma^*(f)) = \sigma^*(F(f)) \in \mathcal{O}(F(X_{\sigma(1)}), \ldots, F(X_{\sigma(n)}); Y)$$

DEFINITION 7.12. *Let* $\mathcal{S}$ *be the operad* $\mathrm{Op}(\mathsf{Set})$*, where we give* $\mathsf{Set}$ *the cartesian monoidal structure. So a type in* $\mathcal{S}$ *is a set, and an operation* $f \in \mathcal{S}(X_1, \ldots, X_n; Y)$ *is a function*

$$f\colon X_1 \times \cdots \times X_n \to Y$$

DEFINITION 7.13. *Let* $\mathcal{O}$ *be an operad. Then an operad algebra* $F$ *of* $\mathcal{O}$ *is an operad morphism from* $\mathcal{O}$ *to* $\mathcal{S}$*. That is, an operad algebra consists of*

- *A set* $F(X)$ *for every* $X \in \mathcal{O}$
- *A function* $F(f)\colon F(X_1) \times \cdots \times F(X_n) \to Y$ *for every operation* $f \in \mathcal{O}(X_1, \ldots, X_n; Y)$

*such that the conditions in Definition 7.11 hold.*

**Example 7.14.** Suppose that $A \subset \mathbb{R}^n$ is convex, that is $\lambda\, a + (1 - \lambda)\, a' \in A$ whenever $a, a' \in A$, $\lambda \in [0, 1]$. Then we can make an operad algebra of $\mathcal{CC}$, where the single type of $\mathcal{CC}$ is sent to $A$, and $(\lambda_1, \ldots, \lambda_n) \in \mathcal{CC}_n$ is sent to the function $A^n \to A$ defined by

$$(a_1, \ldots, a_n) \mapsto \lambda_1\, a_1 + \cdots + \lambda_n\, a_n$$

Now, in the special case that $\mathcal{O}$ and $\mathcal{O}'$ are operads derived from symmetric monoidal categories, a morphisms between them corresponds to a certain type of functor between the categories. Recall that in section 2.4, we defined three types of functors between monoidal categories: strict monoidal functors, monoidal functors, and lax monoidal functors. It turns out that lax monoidal functors correspond to a morphisms of operads, in a way that we make precise in the next proposition.

Recall that a lax symmetric monoidal functor between symmetric monoidal categories $(\mathsf{C}, \otimes_{\mathsf{C}}, I_{\mathsf{C}})$ and $(\mathsf{D}, \otimes_{\mathsf{D}}, I_{\mathsf{D}})$ consists of a functor $F\colon \mathsf{C} \to \mathsf{D}$, along with natural transformations

$$\varepsilon\colon I_{\mathsf{D}} \to F(I_{\mathsf{C}})$$

$$\mu_{X,Y}\colon F(X) \otimes_{\mathsf{D}} F(Y) \to F(X \otimes_{\mathsf{C}} Y)$$

and coherence conditions that imply that we can form a unique morphism

$$\mu_{X_1, \ldots, X_n}\colon F(X_1) \otimes_{\mathsf{D}} \cdots \otimes_{\mathsf{D}} F(X_n) \to F(X_1 \otimes_{\mathsf{C}} \cdots \otimes_{\mathsf{C}} X_n)$$

for every $X_1, \ldots, X_n \in \mathsf{C}$.

DEFINITION 7.15. *Suppose that* $(\mathsf{C}, \otimes_{\mathsf{C}}, I_{\mathsf{C}})$ *and* $(\mathsf{D}, \otimes_{\mathsf{D}}, I_{\mathsf{D}})$ *are symmetric monoidal categories, and* $(F, \varepsilon, \mu)$ *is a lax monoidal functor between them. Then define* $\mathrm{Op}(F)$*, an operad morphism between* $\mathrm{Op}(\mathsf{C})$ *and* $\mathrm{Op}(\mathsf{D})$*, by*

- *On types* $F \in \mathrm{Op}(\mathsf{C})_0 = \mathsf{C}_0$*, define*

$$\mathrm{Op}(F)(X) = F(X)$$

- *On operations* $f \in \mathrm{Op}(\mathsf{C})(X_1 \ldots, X_n; Y) = \mathrm{Hom}_{\mathsf{C}}(X_1 \otimes_{\mathsf{C}} \cdots \otimes_{\mathsf{C}} X_n, Y)$*, define*

$$\mathrm{Op}(F)(f) = F(X_1) \otimes_{\mathsf{D}} \cdots \otimes_{\mathsf{D}} F(X_n) \xrightarrow{\mu_{X_1, \ldots, X_n}} F(X_1 \otimes_{\mathsf{C}} \cdots \otimes_{\mathsf{C}} X_n) \xrightarrow{F(f)} F(Y)$$

PROPOSITION 7.16. *The previous definition makes* $\mathrm{Op}$ *into a functor from* $\mathsf{SMC}$*, the category of symmetric monoidal categories and lax monoidal functors between them, to* $\mathsf{Op}$*, the category of (symmetric) operads and operad morphisms between them.*

**Proof.** See Leinster [70] Chapter 2. ☐

Corollary 7.17. *If $F$ is a lax symmetric monoidal functor from a SMC $(\mathsf{C}, \otimes, I)$ to $(\mathsf{Set}, \times, 1)$, we can construct an operad algebra* $\mathrm{Op}(F)$ *of* $\mathrm{Op}(\mathsf{C})$.

## 7.5. Open graphs as an operad algebra of undirected wiring diagrams

In general, algebras of the operad of undirected wiring diagrams are all about "gluing together" things. Such algebras are also known as cospan-algebras, and have been studied in depth in Fong and Spivak [73]. One of the simplest things that we can glue together is graphs. In this example, we work out precisely what this means.

Definition 7.18. *A* ***graph*** *$G$ is a functor from the category*

$$\mathsf{Gr} \;:=\; E \underset{\mathrm{tgt}}{\overset{\mathrm{src}}{\rightrightarrows}} V$$

*into* $\mathsf{FinSet}$*. That is, a graph consists of*

- *A finite set of* ***vertices*** $G(V)$
- *A finite set of* ***edges*** $G(E)$
- *A* ***source map*** $G(\mathrm{src})\colon G(E) \to G(V)$
- *A* ***target map*** $G(\mathrm{tgt})\colon G(E) \to G(V)$

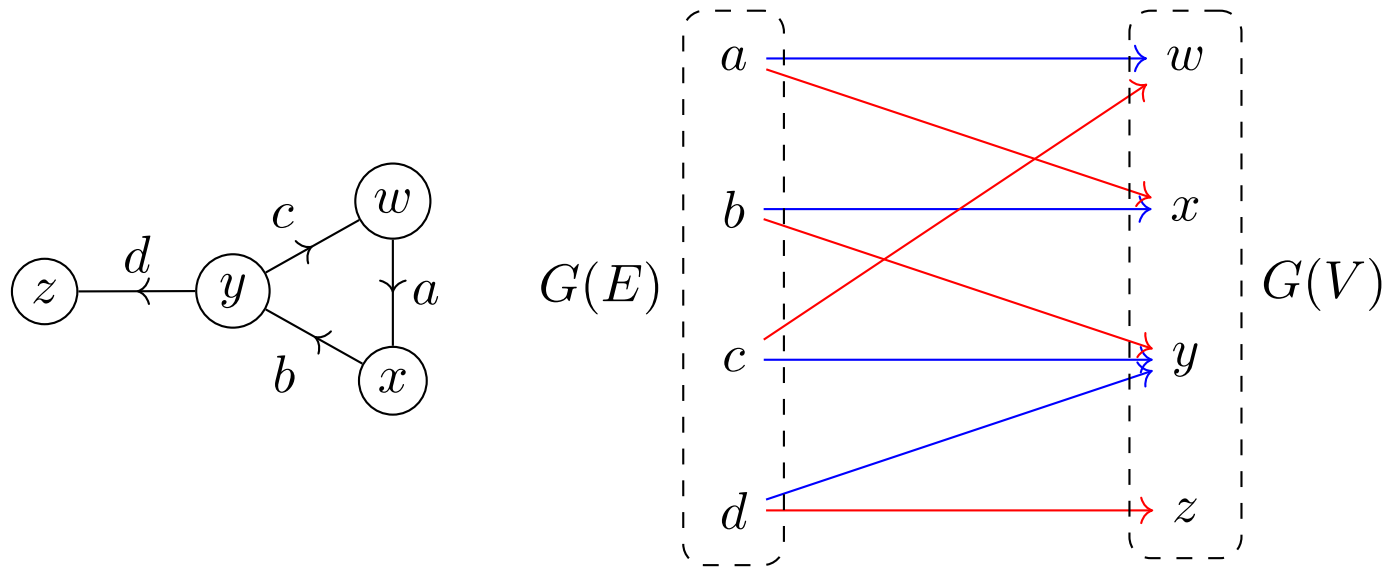

**Figure 7.10.** A graph along with its representation as a functor. The blue arrows represent the source mapping and the red arrows represent the target mapping.

See Figure 7.10 for an example graph, along with a more abstract representation in as a functor into $\mathsf{FinSet}$.

Definition 7.19. *An* ***open graph*** *with interface $X$ is a graph $G$ along with a map $\iota\colon X \to G(V)$. We call $\iota$ the* ***interface map****.*

We think of $\iota$ as picking out the "external" vertices of $G$. An example open graph is pictured in Figure 7.11.

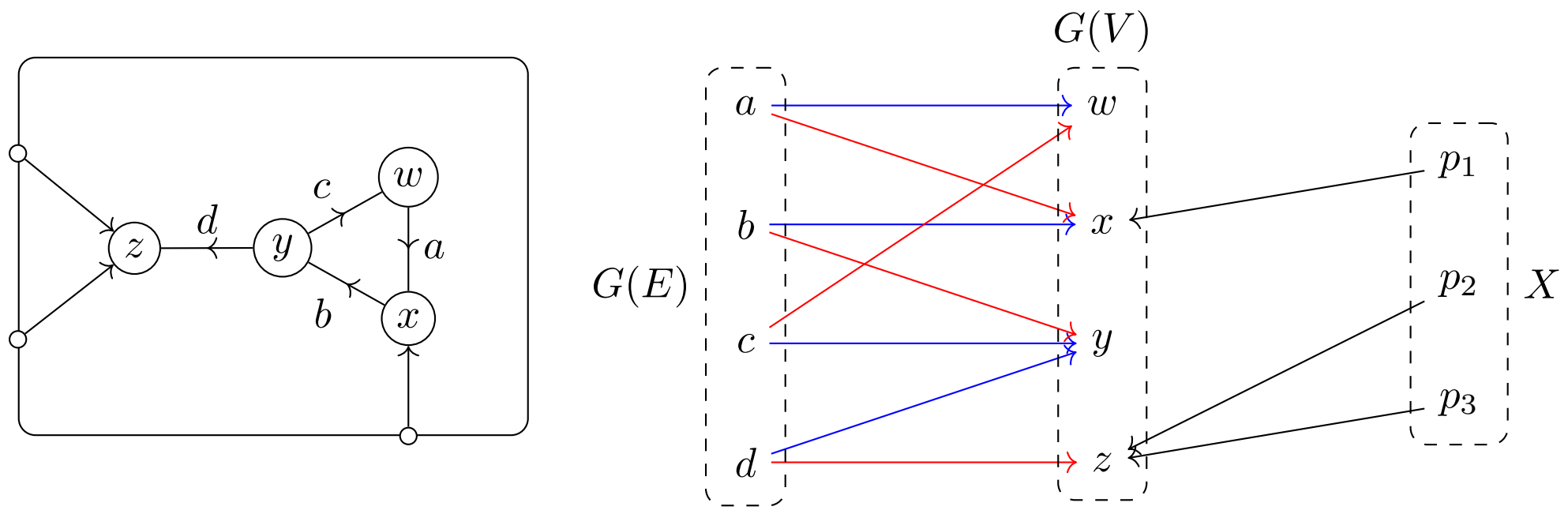

**Figure 7.11.** An open graph, pictured as a graph and as a mapping of sets.

We make open graphs into an operad algebra of operad of wiring diagrams by constructing a lax symmetric monoidal functor OpenGr from FinCospan to Set in the following way. For a finite set $X$, we let OpenGr$(X)$ be the set of equivalence classes of open graphs with interface $X$, where we declare that two open graphs isomorphic in a way that respects the interface maps are equivalent. Discussion of this process of taking equivalence classes can be found in Baez and Pollard [74, Theorem 6]. Then, if $G$ is an open graph on $X$, and $f: X \to Z \leftarrow Y: g$ is a cospan, we define OpenGr$(f, g)(G)$ using the following diagram.

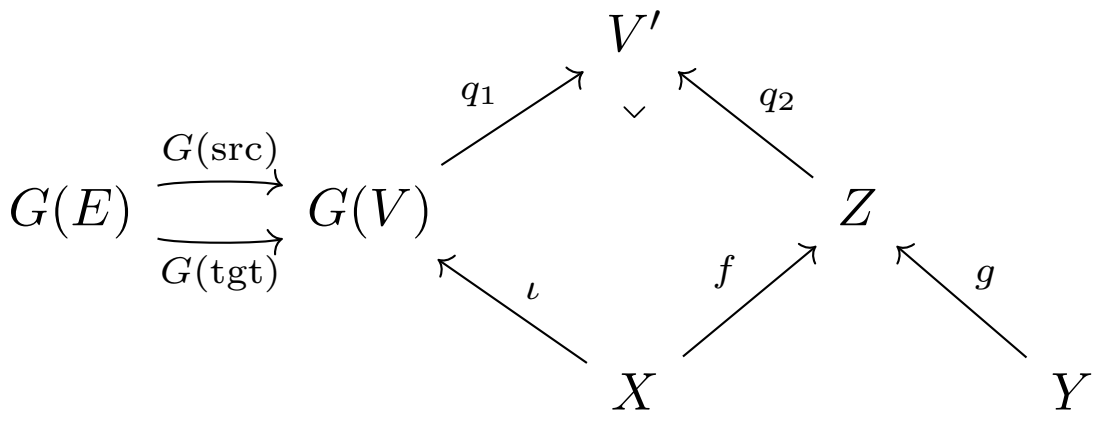

We then let

$$\begin{aligned}\mathrm{OpenGr}(f,g)(G)(V) &= V'\\ \mathrm{OpenGr}(f,g)(G)(E) &= G(E)\\ \mathrm{OpenGr}(f,g)(G)(\mathrm{src}) &= q_1 \circ G(\mathrm{src})\\ \mathrm{OpenGr}(f,g)(G)(\mathrm{tgt}) &= q_1 \circ G(\mathrm{tgt})\end{aligned}$$

Finally, the interface map from $Y$ to OpenGr$(f, g)(G)(V)$ is simply $q_2 \circ g$. The idea here is that vertices in the graphs are "glued together" by the cospan.

The laxator $\mu$ for OpenGr simply takes the disjoint union of graphs. That is, if $G$ is an open graph on $X$ and $G'$ is an open graph on $X'$, then $G \sqcup G'$ is naturally an open graph on $X \sqcup X'$, and from this we get

$$\mu: \mathrm{OpenGr}(X) \times \mathrm{OpenGr}(X') \to \mathrm{OpenGr}(X \sqcup X')$$

When we put this all together, we get to compose graphs using undirected wiring diagrams.

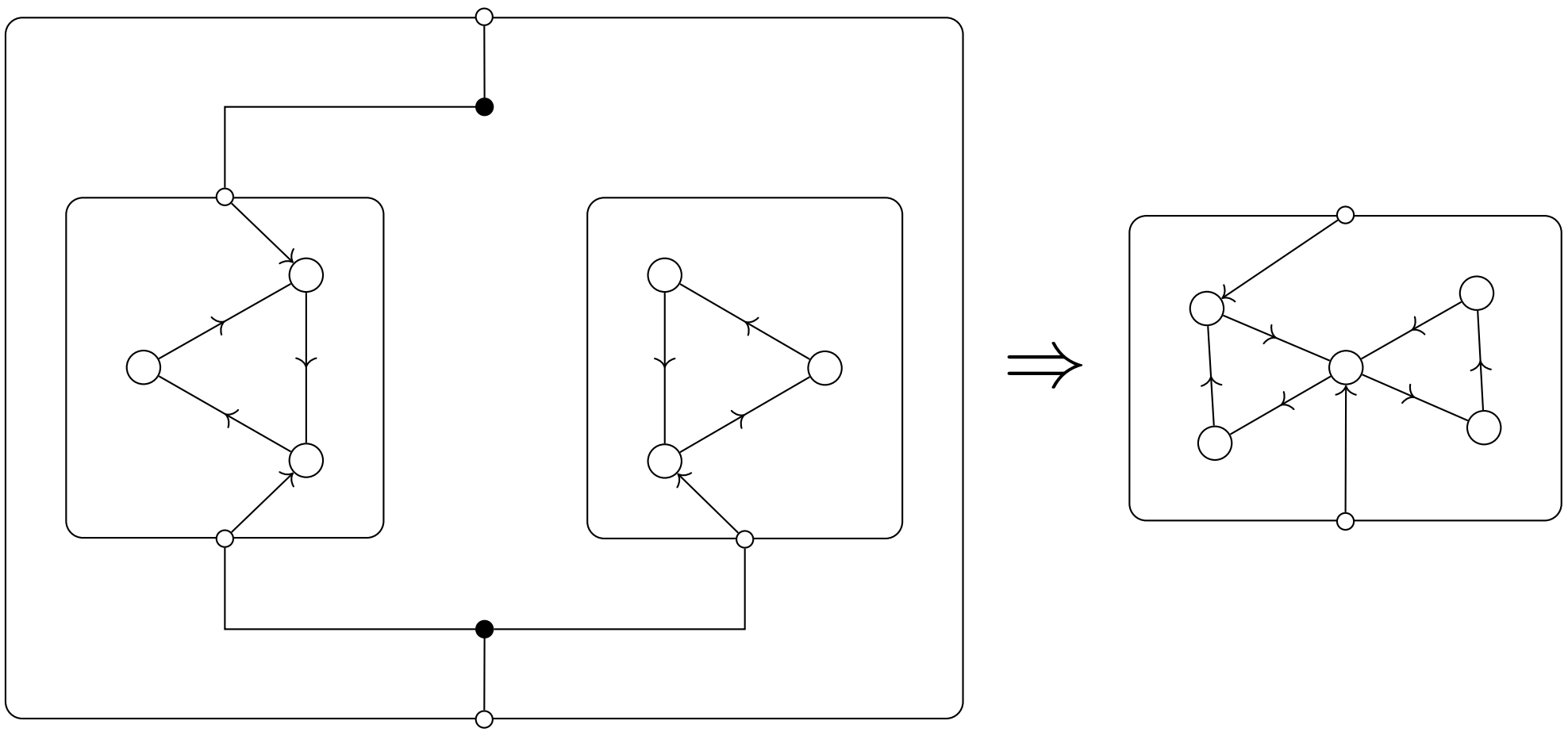

**Figure 7.12.** Composition of open graphs using undirected wiring diagrams.

## 7.6. Higher operad algebras

The mantra of categorification says to replace sets with categories. We can do this in the definition of an operad algebra by replacing the operad morphism into $\mathcal{S} = \mathrm{Op}(\mathsf{Set})$ with an operad morphism into $\mathcal{C} = \mathrm{Op}(\mathsf{Cat})$, where Cat has the symmetric monoidal structure given by cartesian product. When this happens, we must weaken the laws to only hold up to natural isomorphism in the categories.

We as of yet do not know precisely the correct definition of "higher operad" and "higher operad algebra" that would allow for this weakened notion of operad algebra. However, in the special case where the operad comes from a symmetric monoidal category, and the algebra comes from a lax symmetric monoidal functor, the correct type of "weakening" has been studied, and is called a lax symmetric monoidal pseudofunctor. This generalization to lax symmetric pseudofunctors allows the coherence conditions for the lax symmetric monoidal structure on the functor only to have to hold up to natural isomorphism. It is not worth going over the details of this before we see it in use, because it makes the most sense in context. The definition of lax symmetric monoidal pseudofunctors was originally given in Street and Day [75], and was recently used in Baez, Courser and Vasilakopoulou [69] in order to formalize a "higher" notion of decorated cospans.

# 8. Composition of Port-Hamiltonian Systems

## 8.1. Definitions

Finally, we have laid all of the ground-work to compose port-Hamiltonian systems. As said before, we compose port-Hamiltonian systems by means of the theory of operads and operad algebras. However, as noted before, our approach is only "operad-inspired", as we do not yet have a notion of higher operad for which a lax symmetric monoidal pseudofunctor produces a higher operad algebra. Thus, the main goal of this chapter is simply to lay out a sketch of how to construct a lax symmetric monoidal pseudofunctor from DiracRel to Cat. This will only be a sketch because we have not as of yet proved the higher coherence conditions required for a lax symmetric monoidal pseudofunctor. We then explain how this functor models composition of port-Hamiltonian systems by giving several examples.

The fundamental idea for how to compose port-Hamiltonian systems has been understood on a non-categorical level for a while, for instance in van der Schaft and Jeltsema [14, Chapter 6], or van der Schaft and Maschke [76]. Our contribution is to systematize what it means to compose port-Hamiltonian systems, in the same way that the idea of a group systematizes what it means to compose symmetries.

Now, to make the functor from DiracRel to Cat, recall Definition 6.19 of the category $\mathrm{pH}(B)$ for a bond space $B$. We now extend pH into a lax symmetric monoidal functor from DiracRel to Cat. First, we define the action of pH on morphisms of DiracRel, i.e. Dirac relations.

Construction 8.1. *If $B$ and $B'$ are bond spaces, and $R\colon B \nrightarrow B'$ is a Dirac relation between them, define $\mathrm{pH}(R)\colon \mathrm{pH}(B) \to \mathrm{pH}(B')$ in the following way.*

- *On objects $(X, H, D) \in \mathrm{pH}(B)$ (which are port-Hamiltonian systems on the interface $B$), define*

  $$\mathrm{pH}(R)(X, H, D) = (X, H, \mathrm{Triv}(R) \circ D)$$

  *where $\mathrm{Triv}_X\colon \mathsf{DiracRel} \to \mathsf{DiracRel}(X)$ is as defined in 5.45.*

  *This is a port-Hamiltonian system on the interface $B'$ because*

  $$D \in \mathsf{DiracRel}(X)(TX \oplus T^*X, \mathrm{Triv}(B))$$

  *and*

  $$\mathrm{Triv}_X(R) \in \mathsf{DiracRel}(X)(\mathrm{Triv}(B), \mathrm{Triv}(B'))$$

  *so*

  $$\mathrm{Triv}_X(R) \circ D \in \mathsf{DiracRel}(X)(TX \oplus T^*X, \mathrm{Triv}(B'))$$

  *as required for $(X, H, \mathrm{Triv}_X(R) \circ D)$ to be a port-Hamiltonian system on the interface $B'$.*

- *On forwards morphisms $f\colon (X, H, D) \to (X', H', D')$, $\mathrm{pH}(R)$ is just the identity. To see why this works, recall from Definition 6.16 that a forwards morphism $f\colon (X, H, D) \to (X', H', D')$ is a smooth map $f\colon X \to X'$ such that*

  1. *$H = H' \circ f$*

2. *For every* $x \in X$*,* $\xi \in T_x X$*,* $\varphi \in T_x^* X$*,* $b \in B$*,*

$$(\xi, f^*\varphi, b) \in D_x \Rightarrow (f_* \xi, \varphi, b) \in D'_{f(x)}$$

*Well, upon inspection, such a map* $f: X \to X'$ *is* also *a forwards morphism*

$$f: (X, H, \mathrm{Triv}(R) \circ D) \to (X', H', \mathrm{Triv}(R) \circ D')$$

*To see why, first of all property 1 clearly still holds. Moreover, for any* $x \in X$*,* $\xi \in T_x X$*,* $\varphi \in T_x^*$*,* $b' \in B'$*, if*

$$(\xi, f^*\varphi, b') \in (\mathrm{Triv}(R) \circ D)_x = R \circ D_x$$

*then there exists* $b \in B$ *such that* $(b, b') \in R$*, and*

$$(\xi, f^*\varphi, b) \in D_x$$

*Thus,*

$$(f_* \xi, \varphi, b) \in D'_{f(x)}$$

*But then, as* $(b, b') \in R$*,*

$$(f_* \xi, \varphi, b') \in R \circ D'_{f(x)}$$

*as required. As* pH$(R)$ *is the identity on morphisms, it clearly respects composition, and so* pH$(R)$ *is a functor.*

PROPOSITION 8.2. *As defined above,* pH *is a functor from* DiracRel *to* Cat.

**Proof.** We must show that for $R: B_1 \nrightarrow B_2$, $R': B_2 \to B_3$, $\mathrm{pH}(R') \circ \mathrm{pH}(R) = \mathrm{pH}(R' \circ R)$. To see this, on objects of $\mathrm{pH}(B_1)$ we have

$$\begin{aligned}(\mathrm{pH}(R') \circ \mathrm{pH}(R))(X, H, D) &= \mathrm{pH}(R')(X, H, \mathrm{Triv}(R) \circ D)\\ &= (X, H, \mathrm{Triv}(R') \circ \mathrm{Triv}(R') \circ D)\\ &= (X, H, \mathrm{Triv}(R' \circ R) \circ D)\\ &= \mathrm{pH}(R' \circ R)(X, H, D)\end{aligned}$$

On morphisms of $\mathrm{pH}(B_1)$, $\mathrm{pH}(R') \circ \mathrm{pH}(R)$ and $\mathrm{pH}(R' \circ R)$ are both the identity, so they are clearly equal. We are done. □

It turns out that pH is in fact a strict functor. The reason why we treat it as a pseudofunctor, however, is that the lax symmetric monoidal structure that we put on it is be weakened in an appropriate way.

We now put a lax symmetric monoidal structure on pH. To do this, we must define morphisms $\varepsilon: 1 \to \mathrm{pH}(\mathbb{B}^0)$ (where 1 is the terminal category), $\mu_{B,B'}: \mathrm{pH}(B) \times \mathrm{pH}(B') \to \mathrm{pH}(B \oplus B')$, show that they are natural, and then finally show that they obey the coherence conditions for making $(\mathrm{pH}, \varepsilon, \mu)$ into a lax symmetric monoidal pseudofunctor. The proof that these obey the coherence conditions will only be a sketch.

CONSTRUCTION 8.3. *Define* $\varepsilon: 1 \to \mathrm{pH}(\mathbb{B}^0)$ *by sending the single element of* 1 *to the port-Hamiltonian system* $(\mathbb{R}^0, 0, D)$*, where* $D$ *is the trivial relation between* $T\mathbb{R}^0 \oplus T^*\mathbb{R}^0$ *and* $\mathrm{Triv}_X(\mathbb{B}^0)$*.*

*For bond spaces* $A, B$*, define*

$$\mu_{A,B}: \mathrm{pH}(A) \times \mathrm{pH}(B) \to \mathrm{pH}(A \oplus B)$$

*by sending a pair* $((X_A, H_A, D_A), (X_B, H_B, D_B))$ *of port-Hamiltonian systems on* $A$ *and* $B$ *respectively to the port-Hamiltonian system*

$$(X_A \times X_B, H_A + H_B, D_A \oplus D_B)$$

*which is a port-Hamiltonian system on the interface* $A \oplus B$*. Here* $H_A + H_B: X_A \times X_B \to \mathbb{R}$ *is defined by*

$$(H_A + H_B)(x_A, x_B) = H_A(x_A) + H_B(x_B)$$

*and $D_A \oplus D_B$ is defined by taking the isomorphisms*

$$T(X_A \times X_B) \oplus T^*(X_A \times X_B) \cong (TX_A \oplus T^*X_B) \times (TX_A \oplus T^*X_B)$$

*and*

$$\mathrm{Triv}_{X_A}(A) \times \mathrm{Triv}_{X_B}(B) \cong \mathrm{Triv}_{X_A \times X_B}(A \oplus B)$$

*and composing them with*

$$D_A \times D_B \subset (TX_A \oplus T^*X_A \oplus \mathrm{Triv}_{X_A}(A)) \times (TX_B \oplus T^*X_B \oplus \mathrm{Triv}_{X_B}(B))$$

*to get*

$$D_A \oplus D_B \subset T(X_A \times X_B) \oplus T^*(X_A \times X_B) \oplus \mathrm{Triv}_{X_A \times X_B}(A \oplus B)$$

*Confusingly enough, note that that this is* not *the direct product in* $\mathsf{DiracRel}(X)$ *for some* $X$*, because because* $D_A$ *and* $D_B$ *are bundles over different manifolds. However, the fibers of this bundle are in fact the direct product of the fibers of* $D_A$ *and* $D_B$*.*

*Now, we are only halfway done here, because* $\mu_{A,B}$ *must be a* functor*. It is quite easy to define* $\mu_{A,B}$ *on morphisms, however, namely*

$$\mu_{A,B}(f\colon X_A \to X_A', g\colon X_B \to X_B') = f \times g\colon X_A \times X_B \to X_A' \times X_B'$$

*This is easily shown to be a forwards morphism of port-Hamiltonian systems.*

Proposition 8.4. *$\mu$ is a natural transformation.*

**Proof.** This proof is difficult because there are a lot of definitions to unpack, but once we have laid everything out it is fairly straightforward. We want to show that $\mu$ is a natural transformation between two functors that go from $\mathsf{DiracRel} \times \mathsf{DiracRel}$ to $\mathsf{Cat}$. Thus, to show the naturality square, we must consider a morphism in $\mathsf{DiracRel}\times\mathsf{DiracRel}$. Consider four bond spaces $A$, $A'$, $B$, $B'$, and two Dirac relations $R_A\colon A \nrightarrow A'$ and $R_B\colon B \nrightarrow B'$. $(R_A, R_B)\colon (A, B) \nrightarrow (A', B')$ is a morphism in $\mathsf{DiracRel} \times \mathsf{DiracRel}$. We must show that the following square commutes:

$$\begin{array}{ccc} \mathrm{pH}(A) \times \mathrm{pH}(B) & \xrightarrow{\mu_{A,B}} & \mathrm{pH}(A \oplus B) \\ {\scriptstyle \mathrm{pH}(R_A)\times\mathrm{pH}(R_B)}\downarrow & & \downarrow{\scriptstyle \mathrm{pH}(R_A \oplus R_B)} \\ \mathrm{pH}(A') \times \mathrm{pH}(B') & \xrightarrow[\mu_{A',B'}]{} & \mathrm{pH}(A' \oplus B') \end{array}$$

Now, as the action of pH on morphisms is the identity, and the action of $\mu$ is cartesian product, this diagram commutes on morphisms. Thus, we must only check objects. Consider two port-Hamiltonian systems $(X_A, H_A, D_A)$ and $(X_B, H_B, D_B)$. Following the upper path, we get

$$\mu_{A,B}((X_A, H_A, D_A), (X_B, H_B, D_B)) = (X_A \times X_B, H_A + H_B, D_A \oplus D_B)$$

And then $\mathrm{pH}(R_A \oplus R_B)$ sends this to

$$(X_A \times X_B, H_A + H_B, \mathrm{Triv}_{X_A \times X_B}(R_A \oplus R_B) \circ (D_A \oplus D_B))$$

Following the lower path, we get

$$\begin{aligned} & \mathrm{pH}(R_A) \times \mathrm{pH}(R_B)((X_A, H_A, D_A), (X_B, H_B, D_B)) \\ = \ & ((X_A, H_A, \mathrm{Triv}_{X_A}(R_A) \circ D_A), (X_B, H_B, \mathrm{Triv}_{X_B}(R_B) \circ D_B)) \end{aligned}$$

And then

$$\begin{aligned} & \mu_{A',B'}((X_A, H_A, \mathrm{Triv}_{X_A}(R_A) \circ D_A), (X_B, H_B, \mathrm{Triv}_{X_B}(R_B) \circ D_B)) \\ = \ & (X_A \times X_B, H_A + H_B, (\mathrm{Triv}_{X_A}(R_A) \circ D_A) \oplus (\mathrm{Triv}_{X_B}(R_B) \circ D_B)) \end{aligned}$$

Thus, the commutativity of the square reduces to checking that

$$\mathrm{Triv}_{X_A \times X_B}(R_A \oplus R_B) \circ (D_A \oplus D_B) = (\mathrm{Triv}_{X_A}(R_A) \circ D_A) \oplus (\mathrm{Triv}_{X_B}(R_B) \circ D_B)$$

This can be checked fiberwise, using the interchange law from the fact that $\oplus$ is a functor from $\mathsf{DiracRel} \times \mathsf{DiracRel}$ to $\mathsf{DiracRel}$. We are done. ☐

In the next proposition, we sketch the proof that $(\mathrm{pH}, \varepsilon, \mu)$ satisfies the coherence conditions.

CONJECTURE 8.5. *$\mu$ and $\varepsilon$ satisfy the coherence conditions for $(\mathrm{pH}, \varepsilon, \mu)$ to be a lax symmetric monoidal pseudofunctor.*

**Proof.** Recall the coherence conditions from section 2.4. The associativity diagram is

$$\begin{array}{ccc}
(\mathrm{pH}(A) \times \mathrm{pH}(B)) \times \mathrm{pH}(C) & \xrightarrow{a_{\mathrm{pH}(A),\mathrm{pH}(B),\mathrm{pH}(C)}} & \mathrm{pH}(A) \times (\mathrm{pH}(B) \times \mathrm{pH}(C)) \\
{\scriptstyle \mu_{A,B} \times 1_{\mathrm{pH}(z)}} \downarrow & & \downarrow {\scriptstyle 1_{\mathrm{pH}(A)} \times \mu_{B,C}} \\
\mathrm{pH}(A \oplus B) \times \mathrm{pH}(C) & & \mathrm{pH}(A) \times \mathrm{pH}(B \oplus C) \\
{\scriptstyle \mu_{A\oplus B,C}} \downarrow & & \downarrow {\scriptstyle \mu_{A,B\oplus C}} \\
\mathrm{pH}((A \oplus B) \oplus C) & \xrightarrow[\mathrm{pH}(a_{A,B,C})]{} & \mathrm{pH}(A \oplus (B \oplus C))
\end{array}$$

the diagrams for the left and right unitors are

$$\begin{array}{ccc}
1 \times \mathrm{pH}(B) & \xrightarrow{\epsilon \times 1_{\mathrm{pH}(B)}} & \mathrm{pH}(\mathbb{B}^0) \times \mathrm{pH}(B) \\
{\scriptstyle l_{\mathrm{pH}(B)}} \downarrow & & \downarrow {\scriptstyle \mu_{\mathbb{B}^0,B}} \\
\mathrm{pH}(B) & \xleftarrow[\mathrm{pH}(l_B)]{} & \mathrm{pH}(\mathbb{B}^0 \oplus B)
\end{array}$$

and

$$\begin{array}{ccc}
\mathrm{pH}(B) \times 1 & \xrightarrow{1_{\mathrm{pH}(B)} \times \epsilon} & \mathrm{pH}(B) \times \mathrm{pH}(\mathbb{B}^0) \\
{\scriptstyle r_{\mathrm{pH}(B)}} \downarrow & & \downarrow {\scriptstyle \mu_{B,\mathbb{B}^0}} \\
\mathrm{pH}(B) & \xleftarrow[\mathrm{pH}(r_B)]{} & \mathrm{pH}(B \oplus \mathbb{B}^0)
\end{array}$$

and finally, the diagram for symmetry is

$$\begin{array}{ccc}
\mathrm{pH}(A) \times \mathrm{pH}(B) & \xrightarrow{B_{\mathrm{pH}(A),\mathrm{pH}(B)}} & \mathrm{pH}(B) \times \mathrm{pH}(A) \\
{\scriptstyle \mu_{A,B}} \downarrow & & \downarrow {\scriptstyle \mu_{B,A}} \\
\mathrm{pH}(A \oplus B) & \xrightarrow[\mathrm{pH}(B_{A,B})]{} & \mathrm{pH}(B \oplus A)
\end{array}$$

To sketch the proof that $(\mathrm{pH}, \varepsilon, \mu)$ is a lax symmetric monoidal pseudofunctor, we must argue that these diagrams commute up to natural isomorphism, and moreover that the natural isomorphisms filling these diagrams themselves satisfy certain coherence conditions. We will argue this as follows. Recall the definition of the laxator $\mu_{A,B}$:

$$((X_A, H_A, D_A), (X_B, H_B, D_B)) \mapsto (X_A \times X_B, H_A + H_B, D_A \oplus D_B)$$

This is composed of three operations: cartesian product of manifolds, sum of functions, and direct product of relations. Intuitively, these operations are associative, unital with respect to $\varepsilon$, and symmetric, and this is what the diagrams are saying. However, note that the cartesian product of manifolds is only weakly associative, unital, and symmetric.

We go through exactly what this means in the case of the diagram for symmetry. If we chase the symmetry diagram starting with an element of $\mathrm{pH}(A) \times \mathrm{pH}(B)$, we get the following:

$$\begin{array}{ccc}
(X_A, H_A, D_A), (X_B, H_B, D_B) & \longmapsto & (X_B, H_B, D_B), (X_A, H_A, D_A) \\
\downarrow & & \downarrow \\
 & & (X_B \times X_A, H_B + H_A, D_B \oplus D_A) \\
 & & \cong \\
(X_A \times X_B, H_A + H_B, D_A \oplus D_B) & \longmapsto & (X_A \times X_B, H_A + H_B, \mathrm{Triv}_{X_A \times X_B}(B_{A,B}) \circ D_A \oplus D_B)
\end{array}$$

Recall that $B_{A,B}\colon A \oplus B \nrightarrow B \oplus A$ is the *braiding*. The isomorphism in the above diagram is the function $X_A \times X_B \to X_B \times X_A$ sending $(x_A, x_B)$ to $(x_B, x_A)$. This is a forwards morphism of port-Hamiltonian systems on $B \oplus A$.

The cases of associativity and unitality are directly analogous; when the diagrams are chased, they need to be completed with isomorphisms, i.e. $(X_A \times X_B) \times X_C \cong X_A \times (X_B \times X_C)$ and $\mathbb{R}^0 \times X_B \cong X_B \cong X_B \times \mathbb{R}^0$.

We then expect that coherence of the isomorphisms making these squares commute will be a consequence of the fact that $\times$ is coherently associative, unital, and symmetric, but we have not as of yet supplied a detailed proof for this. $\square$

At this point, we would like to cite an analogue of Corollary 7.17, and show that this defines some sort of operad algebra of $\mathrm{Op}(\mathsf{DiracRel})$. However, we as of yet do not know the higher notion of operad algebra that this would work with. Thus, in lieue of this, we give some examples of how the composition with this operad works.

## 8.2. Examples of composition

**Example 8.6.** In this example, we model a motor and flywheel system, pictured in Figure 8.1. In this system, a flywheel is connected to a motor. There is an input port that carries electrical power, and an output port that carries rotational power. Recall that "input" here refers to a sign convention, not a causality convention; this same system could be "run backwards" to produce electrical power from rotation power.

We now describe the parts of the system. The motor transforms electric power, in the form of voltage $V$ and current $I$, into rotational power, in the form of torque $\tau$ and rotational velocity $\omega$. Mathematically, this is given by the relation

$$\begin{aligned} \omega &= \kappa V \\ \kappa \tau &= I \end{aligned}$$

In other words, the rotational velocity is proportional to the input voltage, and the torque is proportional to the current. The electric motor has no state.

Then, we attach the electric motor to a flywheel, whose state is its rotational momentum. That is, the state is $X = \mathbb{R}$, and the Hamiltonian is

$$H(L) = \frac{1}{2 I_{\mathrm{rot}}} L^2$$

where $I_{\mathrm{rot}}$ is the moment of inertia of the flywheel (it is unfortunate that there are so few letters in the English language, and that different domains of physics have conflicting uses for those letters; $I$ is used for current and also moment of inertia). Then the kinematic equations for the flywheel are

$$\begin{aligned} \dot{L} &= \tau' \\ \omega &= \frac{\partial H}{\partial L} = \frac{1}{I_{\mathrm{rot}}} L \end{aligned}$$

The electric motor, the flywheel, and the output are all on one shaft, so their rotational velocities are all the same; this is enforced by a flow-matching, effort-summing junction (torque is effort, rotational velocity is flow).

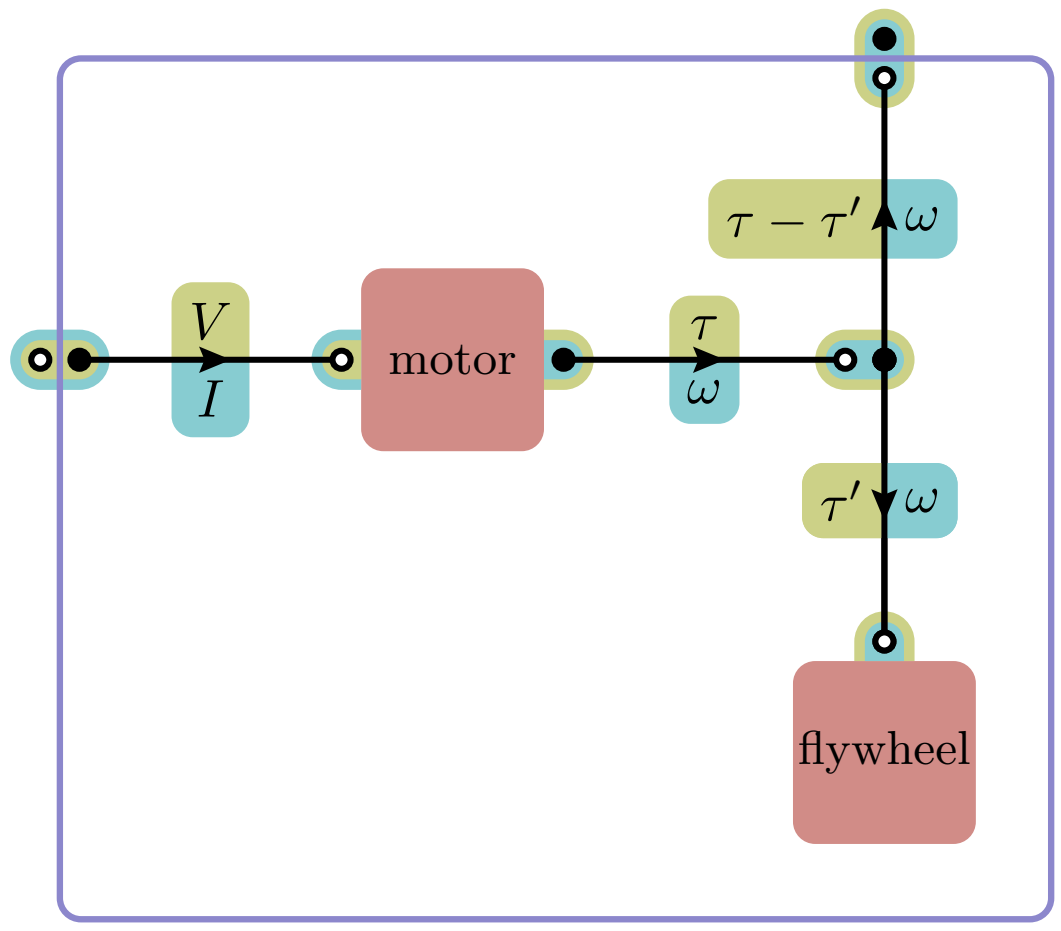

**Figure 8.1.** Dirac diagram for a motor and flywheel.

**Example 8.7.** In this example, we consider a pumping system. The main component of the pumping system is a pump, which uses rotational power from an external port (perhaps from the electrical motor in Example 8.6) to move water from a tank at pressure $p_0$ to a tank at pressure $p_1$. This is governed by the equations

$$p_1 - p_0 = \rho\,\tau$$

$$f = \frac{1}{\rho}\omega$$

for some constant $\rho$. That is, the difference in pressure is proportional to the torque exerted, and the flow through is proportional to the rotational velocity.

The tanks themselves are simple storage systems, analogous to capacitors. The tanks can also be filled and emptied from the outside, using the right and left ports.

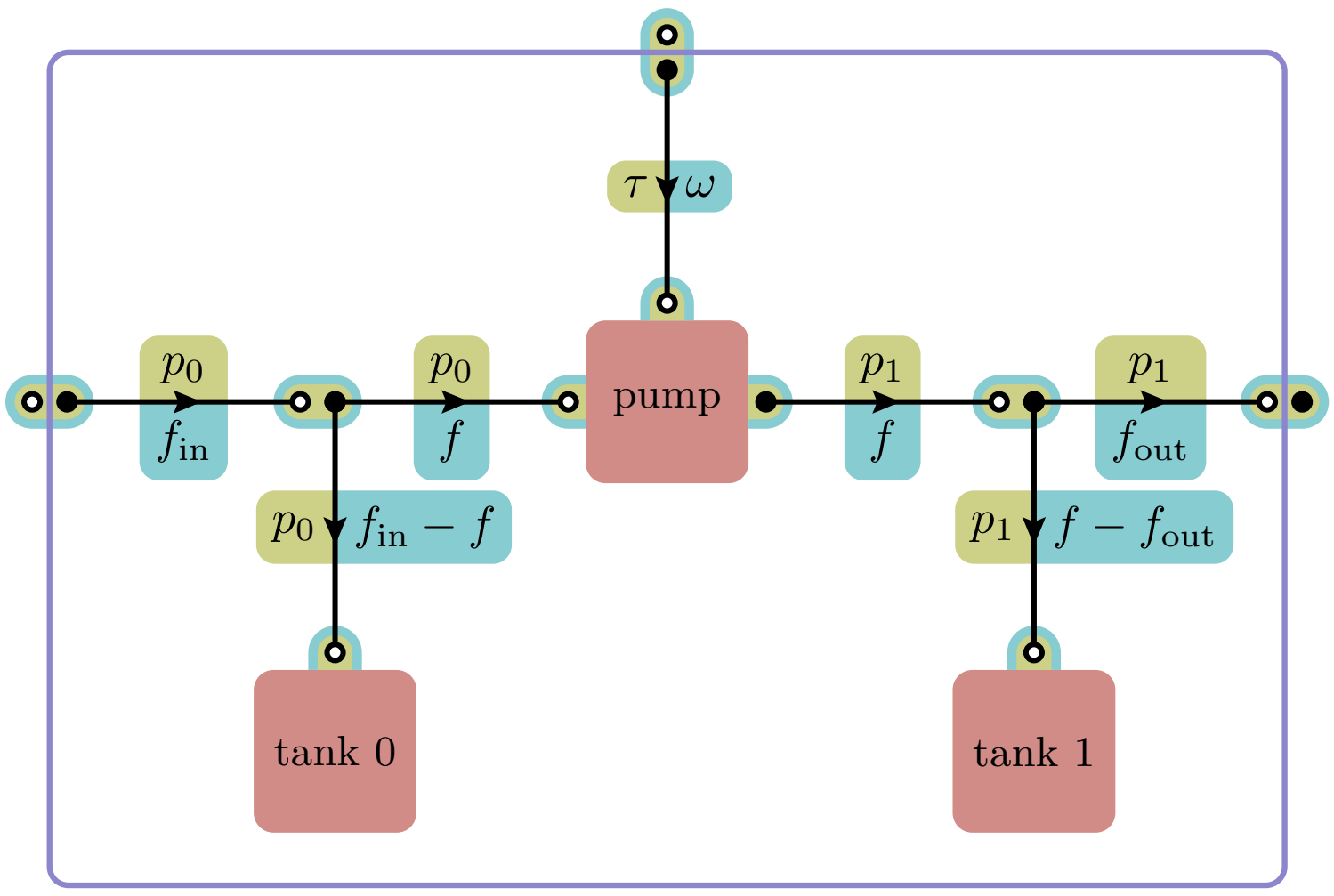

**Figure 8.2.** Dirac diagram for a pumping system.

**Example 8.8.** In this example, we consider a toy hovercar. The hovercar consists of a body with mass, gliding on a frictionless surface with two fans mounted on it, with power to the fans supplied by external wires, and is pictured in Figure 8.3. The point of this example is to demonstrate a non-linear state-dependency in a port-Hamiltonian system.

The configuration space of the hovercar is $Q = \mathbb{R}^2 \times S^1$, i.e. a position $r$ and orientation $\theta$. The state space then is $X = T^* Q$. The Hamiltonian is

$$H(r, \theta, p, L) = \frac{1}{2\,m}|p|^2 + \frac{1}{2\,I}L^2$$

where $m$ is the mass of the hovercar, and $I$ is the moment of intertia. The interface is $\mathbb{B}^2$, with coordinates $(v_x, v_y, F_x, F_y)$. Then the kinematic equations for the hovercar are the following:

$$\begin{aligned}
\dot{r} &= \frac{\partial H}{\partial p} = \frac{1}{m}p \\
\dot{\theta} &= \frac{\partial H}{\partial L} = \frac{1}{I}L \\
\dot{p} &= \begin{bmatrix} F_x \cos\theta \\ F_x \sin\theta \end{bmatrix} \\
\dot{L} &= R\,F_y \\
v_x &= \dot{r}_x \cos\theta + \dot{r}_y \sin\theta \\
v_y &= R\,\dot{\theta}
\end{aligned}$$

Thus, force in the $x$ direction causes the hovercar to move forwards or backwards in the current direction that the hovercar is facing, and force in the $y$ direction causes the hovercar to rotate.

We then supply this force to the hovercar via the fans. We will leave it to the reader to give a specification for the fans; one could include the rotational inertia of the fan blades here, or simply have a stateless conversion from the electrical domain to the mechanical domain.

The interconnection of these parts is then pictured in Figure 8.4.

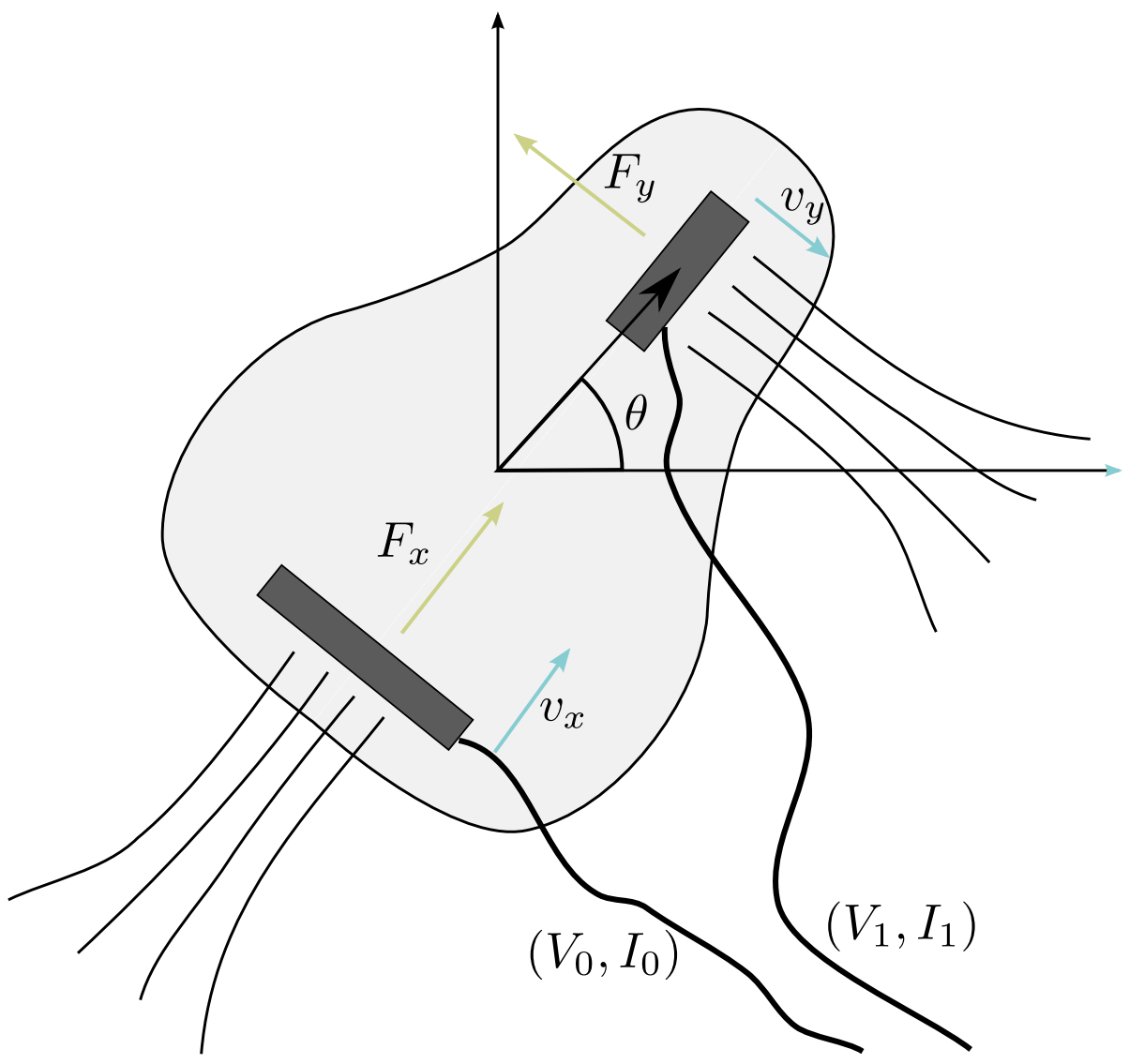

**Figure 8.3.** Schematic of a hovercar.

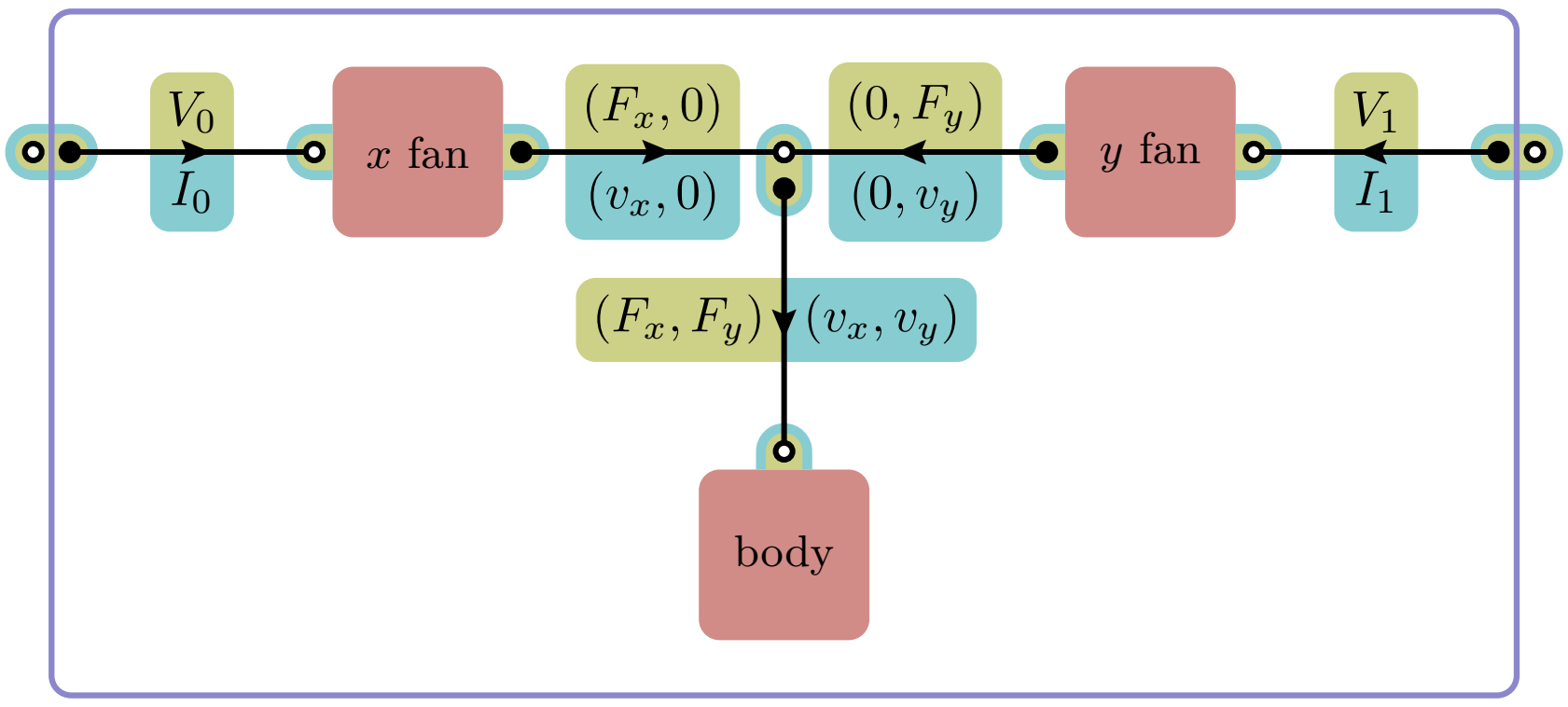

**Figure 8.4.** Dirac diagram for a hovercar.

# 9. Convex Spaces

## 9.1. Convex spaces

This chapter stands alone from most of the previous chapters in this thesis, and only depends on Chapter 2. If the reader has gotten lost in operads and manifolds and Dirac diagrams, we welcome the reader into something quite different, and hopefully refreshing. However, it is necessary background material for the later chapter on thermostatics, where we use convex spaces as state spaces for systems and we use the theory of operad algebras to compose these systems.

The material in this chapter and in the next chapter can be found in Baez, Lynch, and Moeller [1], but we do not go over it in as much detail as in that paper. Rather, we give a "sparse" presentation, showing that we can with not too much effort apply the theory of operads and operad algebras to an entirely different domain.

We start by introducing a generalized conception of "convex space." Recall that if $V$ is a vector space, a convex subset of $V$ is a subset $U \subset V$ such that for all $u_1, u_2 \in U$, and for all $\lambda \in [0,1]$, $\lambda u_1 + (1-\lambda) u_2 \in U$. The following notion treats this operation abstractly, without an underlying vector space. As the properties of a convex space are listed, the reader should note that they hold for the special case of a convex subset of a vector space.

DEFINITION 9.1. *A **convex space** is a set $X$ along with a function $c_\lambda: X \times X \rightarrow X$ for every $\lambda \in [0,1]$ such that the following properties hold.*

- *For all $x, y \in X$, $c_1(x,y) = x$*
- *For all $x \in X$, $\lambda \in [0,1]$ $c_\lambda(x,x) = x$*
- *For all $x, y \in X$, $\lambda \in [0,1]$, $c_\lambda(x,y) = c_{1-\lambda}(y,x)$In this chapter,*
- *For all $x, y, z \in X$, $\lambda, \gamma, \eta \in [0,1]$ such that $\lambda(1-\gamma) = (1-\lambda\gamma)\eta$,*

$$c_\lambda(c_\gamma(x,y), z) = c_{\lambda\gamma}(x, c_\eta(y,z))$$

**Example 9.2.** Any vector space $V$ is a convex space, with

$$c_\lambda(v_1, v_2) = \lambda v_1 + (1-\lambda) v_2$$

for $v_1, v_2 \in V$, $\lambda \in [0,1]$.

**Example 9.3.** If $X$ is a convex space, and $U \subset X$ is such that whenever $x, y \in U$, $c_\lambda(x,y) \in U$, then $c_\lambda$ puts a convex structure on $U$, and we call $U$ a **convex subpace** of $X$. As a sub-example of this, any subset of a vector space closed under convex combinations is a convex space.

**Example 9.4.** As a further subexample of Example 9.3, the ***n*-simplex**

$$\Delta^n = \left\{ p \in \mathbb{R}^{n+1} \middle| p_i \geq 0, \sum_{i=1}^{n} p_i = 1 \right\}$$

is a convex subset of $\mathbb{R}^{n+1}$, representing the set of discrete probability distributions on $n+1$ elements.

**Example 9.5.** As a further subexample of Example 9.3, the **positive orthant**

$$\mathbb{R}^n_{>0} = \{x \in \mathbb{R}^n \mid x_i > 0\}$$

is a convex subset of $\mathbb{R}^n$.

**Example 9.6.** Suppose that $(A, \vee)$ is a **commutative semilattice**, that is $\vee: A \times A \rightarrow A$ is a binary operation such that

- $a \vee a = a$
- $a \vee b = b \vee a$
- $a \vee (b \vee c) = (a \vee b) \vee c$

Then we can put a convex structure on $A$ by defining

$$c_\lambda(a,b) = \begin{cases} a & \text{if } \lambda = 0 \\ b & \text{if } \lambda = 1 \\ a \vee b & \text{otherwise} \end{cases}$$

Examples of this type show that convex spaces are more general than "convex subsets of a vector space."

**Example 9.7.** The extended reals $\bar{\mathbb{R}} = \mathbb{R} \cup \{+\infty, -\infty\}$ have a convex structure uniquely specified by the following equations

$$\begin{array}{rcll} c_\lambda(x, y) & = & \lambda x + (1 - \lambda) y & \text{for } x, y \in \mathbb{R} \\ c_\lambda(x, +\infty) & = & +\infty & \text{for } x \in \mathbb{R} \\ c_\lambda(x, -\infty) & = & -\infty & \text{for } x \in \mathbb{R} \\ c_\lambda(-\infty, \infty) & = & -\infty & \end{array}$$

for $\lambda \in (0, 1)$. This convex structure is a "blend" of Examples 9.3 and 9.6. Note that we have "biased" our definition towards $-\infty$; we could have also biased towards $+\infty$ but for what we do in Chapter 10, this is the physically-relevant definition.

## 9.2. Morphisms of convex spaces

Definition 9.8. *A* ***convex-linear map*** *from a convex space $X$ to a convex space $Y$ is a function $f: X \to Y$ such that*

$$f(c_\lambda(x, x')) = c_\lambda(f(x), f(x'))$$

*for all $x, x' \in X$ and all $\lambda \in [0, 1]$.*

Definition 9.9. *We define* Conv *to be the category of convex spaces and convex-linear maps.*

Proposition 9.10. Conv *is a regular category.*

**Proof.** Convex spaces are models of an algebraic theory [77], [78]. □

Definition 9.11. *As* Conv *is a convex category, we can make its category of relations, via Construction 2.49. We call this* ConvRel*, and it has a natural symmetric monoidal structure given by Proposition 2.51.*

We expand this definition of ConvRel for the sake of clarity.

Definition 9.12. *A* ***convex relation*** *between a convex space $X$ and a convex space $Y$ is a convex subspace $R \subset X \times Y$.* ConvRel *is the category of convex spaces and convex relations, where we use the standard definition for composition of relations.*

In addition to convex-linear maps and convex relations, we consider a third type of map: concave maps.

Definition 9.13. *Suppose that $X$ is a convex space, and $Y$ is a convex space with a partial ordering on it. Then a* ***concave map*** *$f: X \to Y$ is a function such that for all $x, x' \in X$,*

$$f(c_\lambda(x, x')) \geq c_\lambda(f(x), f(x'))$$

We only really consider concave maps into $\bar{\mathbb{R}}$ with the natural ordering on it.

**Example 9.14.** The natural logarithm, $\log: (0, +\infty) \to \bar{\mathbb{R}}$, is a concave map. In fact, we can extend log to $[0, \infty]$ by defining it to be $-\infty$ at 0 and $+\infty$ at $+\infty$, and it is still concave.

We see plenty more examples of concave maps in the next chapter.

# 10. Thermostatics

## 10.1. Thermostatic systems

In this chapter, we give an overview of what we call *thermostatic systems*. These are models of thermodynamical systems that allows us to answer questions about equilibrium (i.e. static points).

In thermostatics, questions about equilibrium are be answered by *constrained maximization of entropy*. We talked about this in Section 1.3; the reader is advised to reread this section if it has been a while. The following definition formalizes exactly what we need to answer questions of equlibrium for thermodynamical systems.

DEFINITION 10.1. *A* ***thermostatic system*** *consists of a convex space $X$ along with a concave function $S: X \to \bar{\mathbb{R}}$, which we call the* ***entropy function****.*

**Example 10.2.** The state space of a classical chemical system is $X = \mathbb{R}_{>0}^{n+2}$, with two coordinates $U$ and $V$ representing total energy and volume, and then $n$ coordinates $N_1, \ldots, N_n$ representing amounts of various chemicals. The entropy function $S: X \to \bar{\mathbb{R}}$ makes a connection between the state of the system and the temperature ($T$), pressure ($p$) and chemical potentials ($\mu_i$), via the equations

$$\begin{aligned} \frac{1}{T} &= \frac{\partial S}{\partial U}(U, V, N_1, \ldots, N_n) \\ -\frac{p}{T} &= \frac{\partial S}{\partial V}(U, V, N_1, \ldots, N_n) \\ \frac{\mu_i}{T} &= \frac{\partial S}{\partial N_i}(U, V, N_1, \ldots, N_n) \end{aligned}$$

Thus, although entropy cannot be measured directly, it plays a crucial role in connecting properties of a chemical system. We call these quantities $\frac{1}{T}$, $-\frac{p}{T}$, and $\frac{\mu_i}{T}$ the **conjugate quantities** to the variables $U, V, N_i$.

**Example 10.3.** The state space of a purely thermal system (like a tank of water that can be heated or cooled, but the volume or amount of water cannot change) is $X = \mathbb{R}_{>0}$. If the system has a constant heat capacity, i.e.

$$U = CT$$

where $U$ is the energy, $C$ is the heat capacity, and $T$ is the temperature, then the entropy function that models this is

$$S(U) = C \log U$$

as then the equation

$$\frac{1}{T} = \frac{\partial S}{\partial U}$$

implies

$$U = CT$$

**Example 10.4.** Often in thermodynamics, it is useful to consider a system with an arbitrarily large heat capacity, a system with a heat capacity so large that the amounts of heat transfered over the course of the relevant process does not change the temperature of the system appreciably. We call this system "the heat bath at temperature $T$", and it has state space $\mathbb{R}$ and entropy function

$$S(\Delta U) = T\,\Delta U$$

where $\Delta U$ is the change in energy from some arbitrary starting point. We use the change in energy rather than the total energy, because the totally energy is effectively infinite (this is captured in the fact that $\Delta U$ can be arbitrarily negative).

**Example 10.5.** Statistical mechanics can be modeled using a state space $X = \mathcal{P}(\Omega)$, the convex space of probability distributions on a finite set $\Omega$, which is the same convex space as the simplex in Example 9.4. Then we use Shannon entropy as the entropy function, i.e.

$$S(p) = -\sum_{i=1}^{n} p_i \log p_i$$

As is normal, we use the convention that $0 \log 0 = 0$.

The conjugate quantities to the $p_i$ are

$$\frac{\partial S}{\partial p_i} = -\log p_i - 1$$

This quantity $-\log p_i$ is called the **surprisal** of the probability distribution at $i$, which roughly measures how unlikely the event $i$ is.

**Example 10.6.** We can also model quantum systems as thermostatic systems. For simplicity, we consider a finite-dimensional Hilbert space $H = \mathbb{C}^n$ with orthonormal basis $e_1, \ldots, e_n$. A **mixed state** on $H$ is a matrix $\rho: H \to H$, with non-negative real eigenvalues $p_1, \ldots, p_n$ such that $p_1 + \cdots + p_n = 1$. The set of mixed states forms a convex space, which we call $X$. We then define $S: X \to \bar{\mathbb{R}}$ by

$$S(\rho) = -\sum_{i=1}^{n} p_i \log p_i$$

where $p_1, \ldots, p_n$ are again the eigenvalues of $\rho$.

We now move on to how one composes thermostatic systems.

## 10.2. Composition of thermostatic systems

We compose thermostatic systems using the same operad-theoretic machinery that we used for port-Hamiltonian systems. The "interface" for a thermostatic system is the state space, so the set of thermostatic systems on a state space is simply the set of entropy functions.

The general idea for how to compose thermostatic systems is the following. Given thermostatic systems $(X_1, S_1), \ldots, (X_n, S_n)$, we first form the thermostatic system $(X_1 \times \cdots \times X_n, S_1 + \cdots + S_n)$. This represents the "independent composition" of the $n$ systems. Then, to model interactions between the systems, we find a convex space $Y$ which models the quantities conserved in the interaction (like total energy, total volume, etc.). We think of $Y$ as a "coarser description" of the whole system. Then, we take a convex relation $R \subset (X_1 \times \cdots \times X_n) \times Y$ that expresses when a state $(x_1, \ldots, x_n) \in X_1 \times \cdots \times X_n$ is compatible with the coarser description $y \in Y$. For instance, it could be compatible only when the conserved quantities in $(x_1, \ldots, x_n)$ add up to the totals in $y$. Finally, we make a thermostatic system on state space $Y$ via the concave map

$$S(y) = \sup_{((x_1, \ldots, x_n), y) \in R} S_1(x_1) + \cdots + S_n(x_n)$$

This process of independent composition followed using a relation to make a system on a new interface is identical on an abstract level to the process of forming a new port-Hamiltonian system, and thus is formalized in the same way: by making a lax symmetric monoidal functor and then turning that into an operad algebra. This time however, we do not need higher operads; this is simply be a lax symmetric monoidal functor into $\mathsf{Set}$ rather than a lax symmetric monoidal pseudofunctor into $\mathsf{Cat}$.

Construction 10.7. *We construct a functor* $\mathrm{Ent}: \mathsf{ConvRel} \to \mathsf{Set}$ *in the following manner. We define* $\mathrm{Ent}$ *on objects by*

$$\mathrm{Ent}(X) = \{S: X \to \bar{\mathbb{R}} \mid S \text{ is concave}\}$$

*We then define* $\mathrm{Ent}$ *on morphisms* $R: X \nrightarrow Y$ *by*

$$\mathrm{Ent}(R)(S \in \mathrm{Ent}(X)) = (y \in Y) \mapsto \sup_{(x,y) \in R} S(x)$$

*We also use* $R_* S$ *to refer to* $\mathrm{Ent}(R)(S)$. *Intuitively,* $R_* S(y)$ *is the maximum entropy of a state compatible with* $y$.

Construction 10.8. *We then define the laxator*

$$\mu_{X,Y}: \mathrm{Ent}(X) \times \mathrm{Ent}(Y) \to \mathrm{Ent}(X \times Y)$$

*by sending* $S \in \mathrm{Ent}(X)$, $T \in \mathrm{Ent}(Y)$ *to* $S + T \in X \times Y$, *where*

$$(S + T)(x, y) = S(x) + T(y)$$

*and the unitor*

$$\epsilon: 1 \to \mathrm{Ent}(1)$$

*by sending the unique element of* $1$ *to the zero function* $1 \to \bar{\mathbb{R}}$*. That this is a valid lax symmetric monoidal functor is proven in Baez, Lynch, and Moeller [1, Theorem 33].*

THEOREM 10.9. $\mathrm{Op}(\mathrm{Ent})$ *is an operad algebra of* $\mathcal{CR} = \mathrm{Op}(\mathsf{ConvRel})$*, where*

$$\mathrm{Op}(\mathrm{Ent})(X) = \{S: X \to \bar{\mathbb{R}} \mid S \textit{ is concave}\}$$

*and the action of a relation*

$$R \in \mathcal{CR}(X_1, \ldots, X_n; Y)$$

*on entropy functions* $S_1: X_1 \to \bar{\mathbb{R}}, \ldots, S_n: X_n \to \bar{\mathbb{R}}$ *is given by*

$$S(y) = \sup_{(x_1, \ldots, x_n, y) \in R} S_1(x_1) + \cdots + S_n(x_n)$$

**Proof.** This follows from Corollary 7.17 applied to the symmetric monoidal functor $(\mathrm{Ent}, \mu, \epsilon)$. □

We now give some examples of composition using this operad algebra.

**Example 10.10.** In this example, we formalize the story about the pie and ice cream given in Section 1.3. We start with two thermostatic systems $(X_1, S_1)$ and $(X_2, S_2)$ with $X_i = \mathbb{R}_{>0}$ and

$$S_i(U) = C_i \log U$$

where $C_i$ is the heat capacity of system $i$. We describe this in terms of the total energy, which takes values in $Y = \mathbb{R}_{>0}$. The relation modeling the conservation of total energy is $R \subset (X_1 \times X_2) \times Y$ defined by

$$R = \{((U_1, U_2), U) \mid U = U_1 + U_2\}$$

This relation is an operation in $\mathrm{Op}(\mathsf{ConvRel})(X_1, X_2; Y)$. We apply this operation to $(S_1, S_2) \in \mathrm{Ent}(X_1) \times \mathrm{Ent}(X_2)$ in the following way. First, we apply the laxator to get

$$S_1 + S_2 \in \mathrm{Ent}(X_1 \times X_2)$$

Second, we apply the functor Ent to $R$ to get

$$R_*(S_1 + S_2)(U) = \sup_{((U_1, U_2), U) \in R} S_1(U_1) + S_2(U_2) = \sup_{U_1 + U_2 = U} S_1(U_1) + S_2(U_2)$$

which is a concave function on $Y$, i.e. an element of $\mathrm{Ent}(Y)$. As discussed in Section 1.3, the supremum over this constraint corresponds to a state where the temperatures of the pie and ice cream have equilibriated. Morever, the explicit form of the entropy function $R_*(S_1 + S_2)$ is

$$R_*(S_1 + S_2)(U) = (C_1 + C_2) \log U + K$$

for some constant $K$, so the system $(Y, R_*(S_1 + S_2))$ acts as a system with constant heat capacity $C_1 + C_2$.

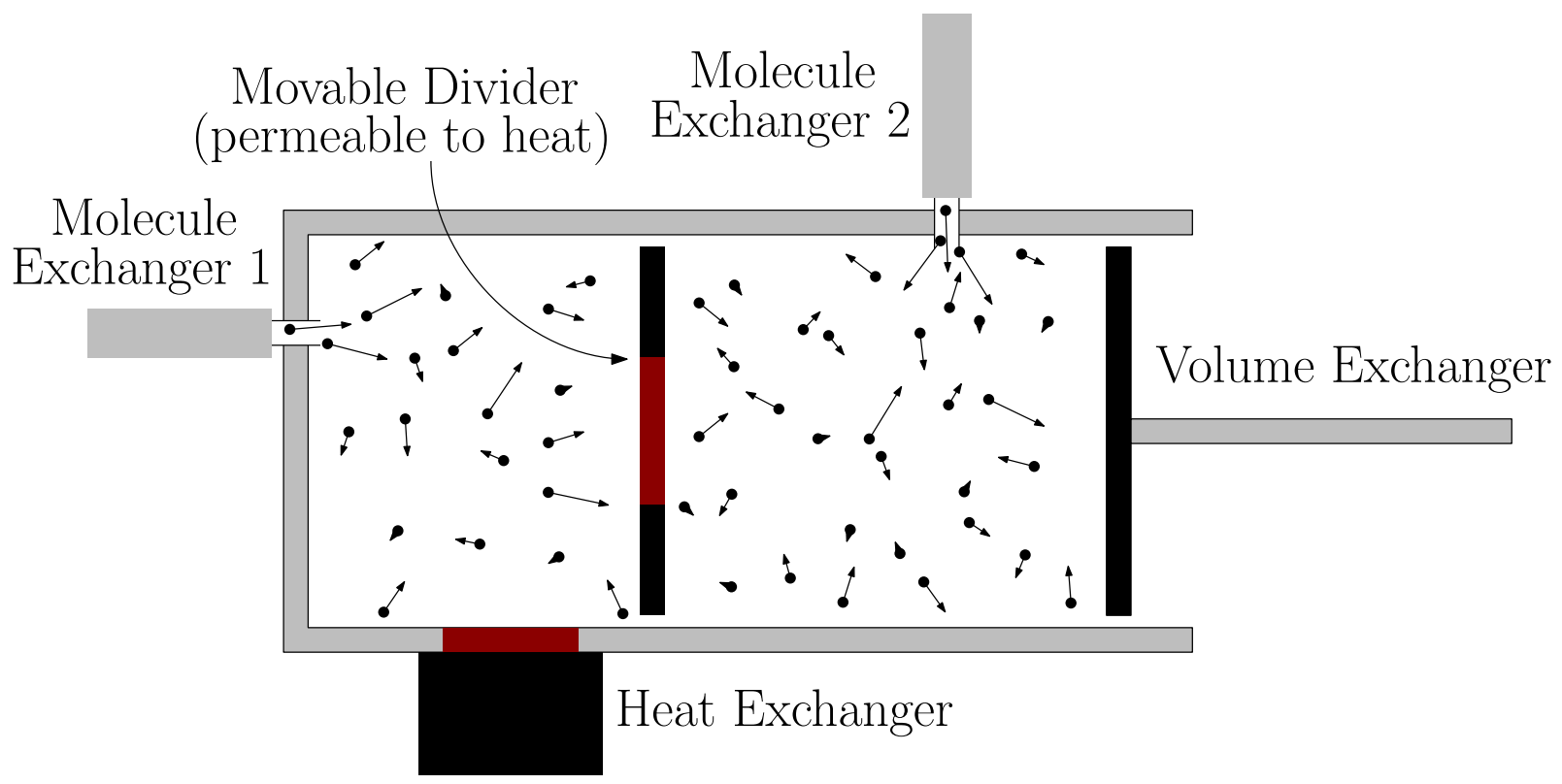

**Figure 10.1.** A setup where two gases can equilibriate temperature and pressure.

**Example 10.11.** We can generalize beyond the previous example to not just equilibriating temperature, but also equlibriating pressure. Consider a system of two ideal gases, as pictured in Figure 10.1. This is a composition of two systems $(X_1, S_1)$, $(X_2, S_2)$, which are each as in Example 10.2 but only one type of molecule. That is, $X_i = \mathbb{R}^3_{>0}$ with coordinates $(U_i, V_i, N_i)$ of energy, volume, and particle number, and the entropy functions are empirically derived from measurement of the quantities

$$\begin{aligned}\frac{1}{T_i} &= \frac{\partial S_i}{\partial U_i}(U_i, V_i, N_i)\\ -\frac{p_i}{T_i} &= \frac{\partial S_i}{\partial V_i}(U_i, V_i, N_i)\\ \frac{\mu_i}{T} &= \frac{\partial S_i}{\partial N_i}(U_i, V_i, N_i)\end{aligned}$$

Then the quantities that are conserved in the interaction of the two systems with each other are $U_1 + U_2$, $V_1 + V_2$, $N_1$, and $N_2$, as the two systems can exchange heat and volume, but not particles. Thus, we coarse grain the system in terms of $Y = \mathbb{R}^4_{>0}$ with coordinates $(U, V, N_1^{\text{ext}}, N_2^{\text{ext}})$, using the relation defined by the equations

$$\begin{aligned}U_1 + U_2 &= U\\ V_1 + V_2 &= V\\ N_1 &= N_1^{\text{ext}}\\ N_2 &= N_2^{\text{ext}}\end{aligned}$$

When we maximize entropy with respect to these constraints, by similar logic as before we end up equilibriating the conjugate variables to the quanties that are allowed to flow between systems. That is, we allow heat to flow between systems, so we end up with

$$\frac{1}{T_1} = \frac{1}{T_2}$$

at equilibrium. We also allow volume to flow between systems (by moving the divider), so we end up with

$$-\frac{p_1}{T_1} = -\frac{p_2}{T_2}$$

at equilibrium. Thus, the temperature and the pressure equilibriate.

**Example 10.12.** Suppose that we take a chemical system $(X_1 = \mathbb{R}^{m+2}, S_1)$ and let it freely exchange heat with a heat bath $\left(X_2 = \mathbb{R}, S_2(U) = \frac{1}{T}U\right)$ at fixed temperature $T$ (see Example 10.4), assuming that the total energy was fixed to 0 (i.e., all of the energy in the system came out of the heat bath originally). The result is a system expressed in terms of volume and particle numbers as

$$S(V, N_1, \ldots, N_m) = \sup_{U + U' = 0} S_1(U, V, N_1, \ldots, N_m) + \frac{1}{T}U' = \sup_U S_1(U, V, N_1, \ldots, N_m) - \frac{1}{T}U$$

This is the Legendre transform of $S_1$, and is known as "Helmholtz free entropy". This quantity is used to investigate systems at constant temperature; our formalism shows that this is a natural consequence of composing a system with a heat bath.

Note that we end up describing the system in terms of $Y = \mathbb{R}^{m+1}$, which is not a complete description of conserved quantities because we also have energy conserved. But because we are dealing with a heat bath, it matters not the total energy as long as it is fixed, so we might as well just assume the total energy is 0 and drop an unnecessary variable from $Y$.

We can generalize the previous example to fix both the pressure and temperature of a system by attaching two thermostatic systems $\left(X_2 = \mathbb{R}, S_2(U) = \frac{1}{T}U\right)$ and $\left(X_3 = \mathbb{R}, S_3(V) = \frac{-p}{T}V\right)$.

## 10.3. Chemical reactions

In this section we model a question that is frequently asked by chemists: what will the equilibrium of a reaction (or collection of reactions) be, at a certain temperature and pressure? To pose this question, we must start by giving the data of a collection of reactions in a mathematical format. To get a sense for what this data must look like, we give some examples of reactions.

**Example 10.13.** Water is formed by reacting hydrogen and oxygen.

$$2\,H_2 + O_2 \rightleftharpoons 2\,H_2O$$

**Example 10.14.** The Haber process reacts atmospheric nitrogen and hydrogen to get ammonia:

$$N_2 + 3\,H_2 \rightleftharpoons 2\,NH_3$$

In general, to specify a reaction, we first give a collection of types of molecules (which we call **species**) present in the reaction, which we call $\Sigma$. Then we have a function $i: \Sigma \to \mathbb{N}$ which counts how many times each species appears on the left hand side of the reaction, and a function $o: \Sigma \to \mathbb{N}$ which counts how many times each species appears on the right hand side.

DEFINITION 10.15. *A **reaction** is a tuple $(\Sigma, i: \Sigma \to \mathbb{N}, o: \Sigma \to \mathbb{N})$ where $\Sigma$ is a finite set. If $s \in \Sigma$, we call $i(s)$ the **input stoichiometric coefficient** in the reaction, and $o(s)$ the **output stoichiometric coefficient**.*

**Example 10.16.** For Example 10.13, we have $\Sigma = \{H_2, O_2, H_2O\}$, and

$$\begin{aligned} i(H_2) &= 2 \\ i(O_2) &= 1 \\ i(H_2O) &= 0 \\ o(H_2) &= 0 \\ o(O_2) &= 0 \\ o(H_2O) &= 2 \end{aligned}$$

**Example 10.17.** For Example 10.14, we have $\Sigma = \{N_2, H_2, NH_3\}$ and

$$\begin{aligned} i(N_2) &= 1 \\ i(H_2) &= 3 \\ i(NH_3) &= 0 \\ o(N_2) &= 0 \\ o(H_2) &= 0 \\ o(H_2O) &= 2 \end{aligned}$$

In general, we might have multiple reactions sharing the same set of species.

DEFINITION 10.18. *A **reaction network** is a tuple $(\Sigma, R, \{i_r: \Sigma \to \mathbb{N}\}_{r \in R}, \{o_r: \Sigma \to \mathbb{N}\}_{r \in R})$, where $R$ is the set of **reactions**.*

Note that a reaction can be identified with a reaction network with just one reaction.

Chemists measure molecules in moles. This is not a unit, it is simply a very large number, $N_A = 6.0221415 \times 10^{23}$ to be precise. Chemists will say "4 moles of hydrogen" to mean $4N_A$ molecules of hydrogen.

DEFINITION 10.19. *A **state** of a reaction network $(\Sigma, R, \{i_r: \Sigma \to \mathbb{N}\}_{r \in R}, \{o_r: \Sigma \to \mathbb{N}\}_{r \in R})$ is an element of $\mathbb{R}_{\geq 0}^{\Sigma}$.*

The meaning of a state $x \in \mathbb{R}_{\geq 0}^{\Sigma}$ is that we model some situation in the world where there are $x_s$ moles of species $s$, for each $s \in \Sigma$. Chemists also measure reaction progress in moles. When a reaction $(\Sigma, i, o)$ has progressed for $x$ moles in a solution, it means that for each species $s \in \Sigma$, $i(s)\,x$ moles of that species have been used, and $o(s)\,x$ moles have been produced, for a total change of $(o(s) - i(s))\,x$ moles. We can represent this more succinctly with the following definition.

DEFINITION 10.20. *For a reaction $(\Sigma, i, o)$, the **stoichiometric number** of a species $s \in S$ is $o(s) - i(s)$. The **stoichiometric vector** for the reaction is $\nu \in \mathbb{N}^S$ defined by $\nu_s = o(s) - i(s)$. If $(\Sigma, R, \{i_r: \Sigma \to \mathbb{N}\}_{r \in R}, \{o_r: \Sigma \to \mathbb{N}\}_{r \in R})$ is a reaction network, then the corresponding **stoichiometric matrix** is defined by $\nu_{sr} = o_r(s) - i_r(s)$.*

The use of the stoichiometric matrix is that if reaction $r$ has progressed for $\xi_r$ moles (i.e. $\xi \in \mathbb{R}^R$), then the total change in all of the species given by the matrix multiplication $\nu\xi$, as $\nu$ is a linear map from $\mathbb{R}^R$ to $\mathbb{R}^S$.

Definition 10.21. *For a reaction network* $(\Sigma, R, \{i_r: \Sigma \to \mathbb{N}\}_{r\in R}, \{o_r: \Sigma \to \mathbb{N}\}_{r\in R})$*, the* ***stoichiometric subspace*** $\mathrm{Stoch} \subset \mathbb{R}^\Sigma$ *is defined by* $\mathrm{Stoch} = \mathrm{im}\,\nu$*.*

Now, suppose that $x: [a, b] \to \mathbb{R}^\Sigma_{\geq 0}$ represents the state over time of an experiment in which the reactions in $(S, R, \{i_r\}, \{o_r\})$ are taking place. Then it must be true that for any $t, s \in [a, b]$, $x(t) - x(s) \in \mathrm{Stoch}$. Thus, the projection $p: \mathbb{R}^\Sigma_{\geq 0} \to \mathbb{R}^\Sigma_{\geq 0}/\mathrm{Stoch}$ must be conserved, where $\mathbb{R}^\Sigma_{\geq 0}/\mathrm{Stoch}$ is the convex space given by quotienting $\mathbb{R}^\Sigma_{\geq 0}$ by the equivalence relation $x \sim x'$ iff $x - x' \in \mathrm{Stoch}$. That is, $p(x(t)) = p(x(s))$ for all $t, s \in [a, b]$. Moreover, if $p(x) = p(x')$, then it is possible to reach $x'$ from $x$ using the reactions, as there exists some $\xi$ such that $x - x' = \nu\xi$, so if we run reaction $r$ for $\xi_r$ moles, then we will get to $x'$ starting at $x$.

We can now apply the methods of thermostatics to investigate questions of equilibrium subject to the constraint that the projection to $\mathbb{R}^\Sigma_{\geq 0}/\mathrm{Stoch}$ is preserved. The typical setup that might be considered in chemistry is the following.

We start with a thermostatic system $\left(X_1 = \mathbb{R}^{2+|\Sigma|}_{\geq 0}, S_1\right)$, a state of which consists of a vector $(U, V, \{N_s\}_{s\in\Sigma})$ of energy, volume, and particle numbers (measured in moles). Reactions in a lab typically take place under constant pressure and constant temperature, so we attach this thermostatic system to a heat bath $\left(X_2 = \mathbb{R}, S_2(U) = \frac{1}{T}U\right)$, and a pressure bath $\left(X_3 = \mathbb{R}, S_3(V) = \frac{-p}{T}V\right)$. As noted before, we describe the conserved quantities using the space $Y = \mathbb{R}^3_{\geq 0}/\mathrm{Stoch}$. Volume and energy are also conserved, but as noted in Example 10.12, because we are connecting to heat/pressure baths, the total energy/volume can just be set to 0 and thus we do not need a variable in $Y$ to account for it.

More formally, the relation $R$ between $X_1 \times X_2 \times X_3$ and $Y$ can be given as follows, using coordinates $(U, V, \vec{N})$ for $X_1$, $U'$ for $X_2$, $V'$ for $X_3$, and $y$ for $Y$.

$$\begin{aligned} U + U' &= 0 \\ V + V' &= 0 \\ p(\vec{N}) &= y \end{aligned}$$

where $p: \mathbb{R}^\Sigma_{\geq 0} \to Y$ is the projection discussed earlier. When we compute the entropy $S = R_*(S_1 + S_2 + S_3)$, we end up getting

$$S(y) = \sup_{U, V, p(\vec{N}) = y} S_1(U, V, \vec{N}) - \frac{1}{T}U + \frac{p}{T}V$$

The quantity

$$\Xi = \sup_{U, V} S_1(U, V, \vec{N}) - \frac{1}{T}U + \frac{p}{T}V$$

is known as Gibbs free entropy or the Planck potential, and maximizing it over all $\vec{N}$ compatible with $y$ is known to give the equilibrium of a collection of reactions; for more information on this see Callen [79, Section 6-7], where $\Xi$ is referred to as a Massieu function. Thus, we see that our framework rederives chemical equilibria.

## 10.4. Statistical mechanics

Statistical mechanics is the study of systems in physics with many orders of magnitude more degrees of freedom than would be feasible to model explicitly. Fortunately, it can be productive to model these systems using probability and statistics, as there are so many degrees of freedom that laws of large numbers apply very strongly.

Statistical mechanics is a very broad subject. Here we give just a brief taste of how thermostatics and statistical mechanics can intersect. It is well-known that the canonical distribution, a common distribution for modeling systems within statistical mechanics, can be derived by maximizing entropy. We give an overview of why this is so, and how this fits into the thermostatic framework.

We start with a probabilistic system, as in Example 10.5. That is, we have a state space $X_1 = P(\Omega)$, for $\Omega$ a finite set, where $P(\Omega)$ is the convex space of probability distributions on $\Omega$, and the entropy function is Shannon entropy:

$$S(p) = -\sum_{\omega \in \Omega} p(\omega) \log p(\omega)$$

We then describe this system in terms of the expected values of a macroscale observable. Specifically, we assume that there is some function $H\colon \Omega \to \mathbb{R}$ that gives the energy of each $\omega \in \Omega$, and that measuring the energy of a general state $p \in P(\Omega)$ gives

$$\langle H \rangle_p = \sum_{\omega \in \Omega} H(\omega)\, p(\omega)$$

Now, suppose that we were to measure the energy of the system $U$. The mathematician E.T. Jaynes would argue that the only reasonable choice of probability distribution that expresses our knowledge of the system is the probability distribution $p$ that maximizes $S(p)$, subject to the constraint that $\langle H \rangle_p = U$; see Jaynes [18, Chapter 11] for more details on this. There are also theorems in probability that state that under some reasonable assumptions, the "most likely" probability distribution producing our outcome is in fact this entropy-maximizing distribution; see Cover and Thomas [80, Section 11.4].

We can model this situation via the relation $R$ between $X = P(\Omega)$ and $Y = \mathbb{R}$ given by the graph of the function $p \mapsto \langle H \rangle_p$. We get an entropy function

$$R_* S(U) = \sup_{p \in P(\Omega), \langle H \rangle_p = U} S(p)$$

This supremum in fact has a unique maximizer, which we calculate by means of Lagrange multipliers. We have two constraints, one of which comes from the fact that $p$ must be a probability distribution

$$\sum_{\omega \in \Omega} p(\omega) \;=\; 1$$

and the other which comes from our observation $U$, i.e.

$$\sum_{\omega \in \Omega} H(\omega)\, p(\omega) \;=\; U$$

We will thus maximize the quantity

$$S(p) - \beta(\langle H \rangle_p) - \gamma\left(\sum_{\omega \in \Omega} p(w)\right) = -\sum_{\omega \in \Omega} p(w) \log(p(\omega)) + \beta H(\omega)\, p(\omega) + \gamma p(\omega)$$

We do this by setting the derivatives with respect to each $p(\omega)$ to 0. We then get

$$1 + \log(p(\omega)) + \beta H(\omega) + \gamma = 0$$

Thus

$$p(\omega) = \mathrm{e}^{-1 - \beta H(\omega) - \gamma}$$

We normalize this to a probability distribution by tweaking $\gamma$ so that

$$\mathrm{e}^{-1-\gamma} = \sum_{\omega \in \Omega} \mathrm{e}^{-\beta H(\omega)}$$

This is called the **partition function** and is typically called $Z(\beta)$. We then have a family of probability distributions given by

$$p_\beta(\omega) = \frac{\mathrm{e}^{-\beta H(\omega)}}{Z(\beta)}$$

This is called the **canonical distribution**. We can now tweak $\beta$ until we find one such that $\langle H \rangle_{p_\beta} = U$. The partition function and the canonical distribution are the basic ingredients in statistical mechanics [81].

The end result here is that the derivation of the canonical distribution fits into the same framework as macroscale thermostatics; that of the operad algebra of thermostatic systems. It is widely appreciated that the maximizing entropy subject to constraints is an important procedure in science; we have shown here that this procedure can be captured within a formal system general enough to handle a wide range of applications, and we hope that this will lay the groundwork for further formalization of thermodynamics and statistical mechanics.